\input pictex \input miniltx \input graphicx.sty \input color.sty \input amssym \font \bbfive = bbm5 \font \bbeight = bbm8 \font \bbten = bbm10 \font \rs
= rsfs10 \font \rssmall = rsfs10 scaled 833 \font \eightbf = cmbx8 \font \eighti = cmmi8 \skewchar \eighti = '177 \font \fouri = cmmi5 scaled 800 \font
\eightit = cmti8 \font \eightrm = cmr8 \font \eightsl = cmsl8 \font \eightsy = cmsy8 \skewchar \eightsy = '60 \font \eighttt = cmtt8 \hyphenchar \eighttt
= -1   \font \sixi = cmmi6 \skewchar \sixi = '177 \font \sixrm = cmr6 \font \sixsy = cmsy6 \skewchar \sixsy =
'60 \font \tensc = cmcsc10   \scriptfont \bffam = \bbeight \scriptscriptfont \bffam = \bbfive \textfont
\bffam = \bbten \newskip \ttglue \def \eightpoint {\def \rm {\fam 0 \eightrm }\relax \textfont 0 = \eightrm \scriptfont 0 = \sixrm \scriptscriptfont 0
= \fiverm \textfont 1 = \eighti \scriptfont 1 = \sixi \scriptscriptfont 1 = \fouri \textfont 2 = \eightsy \scriptfont 2 = \sixsy \scriptscriptfont 2 =
\fivesy \textfont 3 = \tenex \scriptfont 3 = \tenex \scriptscriptfont 3 = \tenex \def \it {\fam \itfam \eightit }\relax \textfont \itfam = \eightit \def
\sl {\fam \slfam \eightsl }\relax \textfont \slfam = \eightsl \def \bf {\fam \bffam \eightbf }\relax \textfont \bffam = \bbeight \scriptfont \bffam =
\bbfive \scriptscriptfont \bffam = \bbfive \def \tt {\fam \ttfam \eighttt }\relax \textfont \ttfam = \eighttt \tt \ttglue = .5em plus.25em minus.15em
\normalbaselineskip = 9pt \def \MF {{\manual opqr}\-{\manual stuq}}\relax \let \sc = \sixrm \let \big = \eightbig \let \rs = \rssmall \setbox \strutbox
= \hbox {\vrule height7pt depth2pt width0pt}\relax \normalbaselines \rm } \def \setfont #1{\font \auxfont =#1 \auxfont } \def \withfont #1#2{{\setfont
{#1}#2}}   \def
\mathbox #1{{\mathchoice {\hbox {\rm #1}} {\hbox {\rm #1}} {\hbox {\eightrm #1}} {\hbox {\sixrm #1}}}} \def \TRUE {Y} \def \FALSE {N} \def \ifundef
#1{\expandafter \ifx \csname #1\endcsname \relax } \def \undefrule {\kern 2pt \vrule width 2pt height 5pt depth 0pt \kern 2pt} \def \UndefLabels {}
\def \possundef #1{\ifundef {#1}\undefrule {\eighttt #1}\undefrule \global \edef \UndefLabels {\UndefLabels #1\par } \else \csname #1\endcsname \fi }
\newcount \secno \secno = 0 \newcount \stno \stno = 0 \newcount \eqcntr \eqcntr = 0 \ifundef {showlabel} \global \def \showlabel {\FALSE } \fi \ifundef
{auxwrite} \global \def \auxwrite {\TRUE } \fi \ifundef {auxread} \global \def \auxread {\TRUE } \fi \def \define #1#2{\global \expandafter \edef \csname
#1\endcsname {#2}} \long \def \error #1{\medskip \noindent {\bf ******* #1}} \def \fatal #1{\error {#1\par Exiting...}\end } \def \advseqnumbering
{\global \advance \stno by 1 \global \eqcntr =0} \def \current {\ifnum \secno = 0 \number \stno \else \number \secno \ifnum \stno = 0 \else .\number
\stno \fi \fi } \def \rem #1{\vadjust {\vbox to 0pt{\vss \hfill \raise 3.5pt \hbox to 0pt{ #1\hss }}}} \font \tiny = cmr6 scaled 800 \def \deflabel
#1#2{\relax \if \TRUE \showlabel \rem {\tiny #1}\fi \ifundef {#1PrimarilyDefined}\relax \define {#1}{#2}\relax \define {#1PrimarilyDefined}{#2}\relax
\if \TRUE \auxwrite \immediate \write 1 {\string \newlabe l {#1}{#2}}\fi \else \edef \old {\csname #1\endcsname }\relax \edef \new {#2}\relax \if \old
\new \else \fatal {Duplicate definition for label ``{\tt #1}'', already defined as ``{\tt \old }''.}\fi \fi } 
\def \label #1 {\deflabel {#1}{\current }} \def \equationmark #1 {\ifundef {InsideBlock} \advseqnumbering \eqno {(\current )} \deflabel {#1}{\current
} \else \global \advance \eqcntr by 1 \edef \subeqmarkaux {\current .\number \eqcntr } \eqno {(\subeqmarkaux )} \deflabel {#1}{\subeqmarkaux } \fi \if
\TRUE \showlabel \hbox {\tiny #1}\fi } \def \lbldeq #1 $$#2$${\ifundef {InsideBlock}\advseqnumbering \edef \lbl {\current }\else \global \advance \eqcntr
by 1 \edef \lbl {\current .\number \eqcntr }\fi $$ #2 \deflabel {#1}{\lbl }\eqno {(\lbl )} $$} \def \split #1.#2.#3.#4;{\global \def \parone {#1}\global
\def \partwo {#2}\global \def \parthree {#3}\global \def \parfour {#4}} \def \NA {NA} \def \ref #1{\split #1.NA.NA.NA;(\possundef {\parone }\ifx \partwo
\NA \else .\partwo \fi )}  \newcount \bibno \bibno = 0  \def \Bibitem #1 #2; #3;
#4 \par {\smallbreak \global \advance \bibno by 1 \item {[\possundef {#1}]} #2, {``#3''}, #4.\par \ifundef {#1PrimarilyDefined}\else \fatal {Duplicate
definition for bibliography item ``{\tt #1}'', already defined in ``{\tt [\csname #1\endcsname ]}''.} \fi \ifundef {#1}\else \edef \prevNum {\csname
#1\endcsname } \ifnum \bibno =\prevNum \else \error {Mismatch bibliography item ``{\tt #1}'', defined earlier (in aux file ?) as ``{\tt \prevNum }''
but should be ``{\tt \number \bibno }''.  Running again should fix this.}  \fi \fi \define {#1PrimarilyDefined}{#2}\relax \if \TRUE \auxwrite \immediate
\write 1 {\string \newbi b {#1}{\number \bibno }}\fi } \def \jrn #1, #2 (#3), #4-#5;{{\sl #1}, {\bf #2} (#3), #4--#5} \def \Article #1 #2; #3; #4 \par
{\Bibitem #1 #2; #3; \jrn #4; \par } \def \references {\begingroup \bigbreak \eightpoint \centerline {\tensc References} \nobreak \medskip \frenchspacing }
\catcode `\@=11 \def \citetrk #1{{\bf \possundef {#1}}} \def \c@ite #1{{\rm [\citetrk {#1}]}} \def \sc@ite [#1]#2{{\rm [\citetrk {#2}\hskip 0.7pt:\hskip
2pt #1]}} \def \du@lcite {\if \pe@k [\expandafter \sc@ite \else \expandafter \c@ite \fi } \def \cite {\futurelet \pe@k \du@lcite } \catcode `\@=12 \def
\Headlines #1#2{\nopagenumbers \headline {\ifnum \pageno = 1 \hfil \else \ifodd \pageno \tensc \hfil \lcase {#1} \hfil \folio \else \tensc \folio \hfil
\lcase {#2} \hfil \fi \fi }} \def \title #1{\medskip \centerline {\withfont {cmbx12}{\ucase {#1}}}}  \long \def \Quote #1\endQuote {\begingroup \leftskip 35pt \rightskip 35pt \parindent 17pt
\eightpoint #1\par \endgroup } \long \def \Abstract #1\endAbstract {\vskip 1cm \Quote \noindent #1\endQuote }   \def \Note #1{\footnote {}{\eightpoint #1}} \def \Date #1 {\Note
{\it Date: #1.}} \newcount \auxone \newcount \auxtwo \newcount \auxthree \def \currenttime {\auxone =\time \auxtwo =\time \divide \auxone by 60 \auxthree
=\auxone \multiply \auxthree by 60 \advance \auxtwo by -\auxthree \ifnum \auxone <10 0\fi \number \auxone :\ifnum \auxtwo <10 0\fi \number \auxtwo } \def
\today {\ifcase \month \or January\or February\or March\or April\or May\or June\or July\or August\or September\or October\or November\or December\fi
{ }\number \day , \number \year }  \def \hojeExtenso {\number \day \ de \ifcase \month \or
janeiro\or fevereiro\or mar\c co\or abril\or maio\or junho\or julho\or agosto\or setembro\or outubro\or novembro\or decembro\fi \ de \number \year }  \def \part #1#2{\vfill \eject \null \vskip 0.3\vsize \withfont {cmbx10 scaled 1440}{\centerline {PART #1}
\vskip 1.5cm \centerline {#2}} \vfill \eject }   \def \fix {\smallskip \noindent $\blacktriangleright $\kern 12pt} 
 \def \ucase #1{\edef \auxvar {\uppercase {#1}}\auxvar } \def \lcase #1{\edef \auxvar {\lowercase {#1}}\auxvar } \def
\emph #1{{\it #1}\/} \def \section #1 \par {\global \advance \secno by 1 \stno = 0 \goodbreak \bigbreak \noindent {\bf \number \secno .\enspace #1.}
\nobreak \medskip \noindent } \def \state #1 #2\par {\begingroup \def \InsideBlock {} \medbreak \noindent \advseqnumbering {\bf \current .\enspace
#1.\enspace \sl #2\par }\medbreak \endgroup } \def \definition #1\par {\state Definition \rm #1\par } \newcount \CloseProofFlag   \long \def \Proof #1\endProof {\begingroup \def
\InsideBlock {} \global \CloseProofFlag =0 \medbreak \noindent {\it Proof.\enspace }#1 \ifnum \CloseProofFlag =0 \hfill $\endproofmarker $ \looseness =
-1 \fi \medbreak \endgroup } \def \quebra #1{#1 $$$$ #1} \def \explica #1#2{\mathrel {\buildrel \hbox {\sixrm #1} \over #2}} \def \explain #1#2{\explica
{\ref {#1}}{#2}}  \def \=#1{\explain {#1}{=}} \def \pilar #1{\vrule height #1 width 0pt}
 \newcount \fnctr \fnctr = 0 \def \fn #1{\global \advance \fnctr by 1 \edef \footnumb {$^{\number \fnctr
}$}\relax \footnote {\footnumb }{\eightpoint #1\par \vskip -10pt}} \def \text #1{\hbox {#1}} \def \bool #1{[{#1}]}  \def
\item #1{\par \noindent \kern 1.1truecm\hangindent 1.1truecm \llap {#1\enspace }\ignorespaces } \def \iItem {\smallskip } \def \Item #1{\smallskip \item
{{\rm #1}}} \newcount \zitemno \zitemno = 0 \def \izitem {\global \zitemno = 0} 
\def \iItemize {\izitem } \def \zitemplus {\global \advance \zitemno by 1 \relax } \def \rzitem {\romannumeral \zitemno } \def \rzitemplus {\zitemplus
\rzitem } \def \zitem {\Item {{\rm (\rzitemplus )}}} \def \iItem {\zitem }  \newcount \nitemno \nitemno = 0  \def \nitem {\global \advance \nitemno by 1 \Item {{\rm (\number \nitemno )}}} \newcount \aitemno \aitemno = -1 \def \boxlet #1{\hbox to
6.5pt{\hfill #1\hfill }} \def \iaitem {\aitemno = -1} \def \aitemconv {\ifcase \aitemno a\or b\or c\or d\or e\or f\or g\or h\or i\or j\or k\or l\or m\or
n\or o\or p\or q\or r\or s\or t\or u\or v\or w\or x\or y\or z\else zzz\fi } \def \aitem {\global \advance \aitemno by 1\Item {(\boxlet \aitemconv )}}
 \def \deflabeloc #1#2{\deflabel {#1}{\current .#2}{\def \showlabel {\FALSE }\deflabel {Local#1}{#2}}} \def \lbldzitem
#1 {\zitem \deflabeloc {#1}{\rzitem }} \def \lbldaitem #1 {\aitem \deflabeloc {#1}{\aitemconv }} \def \aitemmark #1 {\deflabel {#1}{\aitemconv }}
\def \iItemmark #1 {\zitemmark {#1} } \def \zitemmark #1 {\deflabel {#1}{\current .\rzitem }{\def \showlabel {\FALSE }\deflabel {Local#1}{\rzitem
}}} \def \Case #1:{\medskip \noindent {\tensc Case #1:}} \def \<{\left \langle \vrule width 0pt depth 0pt height 8pt } \def \>{\right \rangle } \def
\({\big (} \def \){\big )} \def \ds {\displaystyle } \def \and {\mathchoice {\hbox {\quad and \quad }} {\hbox { and }} {\hbox { and }} {\hbox { and }}}
   \def \IFF
{\kern 7pt\Leftrightarrow \kern 7pt} \def \IMPLY {\kern 7pt \Rightarrow \kern 7pt} \def \for #1{\mathchoice {\quad \forall \,#1} {\hbox { for all } #1}
{\forall #1}{\forall #1}} \def \endproofmarker {\square } \def \"#1{{\it #1}\/} \def \umlaut #1{{\accent "7F #1}} \def \inv {^{-1}} \def \*{\otimes }
\def \caldef #1{\global \expandafter \edef \csname #1\endcsname {{\cal #1}}} \def \mathcal #1{{\cal #1}} \def \bfdef #1{\global \expandafter \edef
\csname #1\endcsname {{\bf #1}}} \bfdef N \bfdef Z \bfdef C \bfdef R  \def \exists {\mathchar "0239\kern 1pt } \if \TRUE
\auxread \IfFileExists {\jobname .aux}{\input \jobname .aux}{\null } \fi \if \TRUE \auxwrite \immediate \openout 1 \jobname .aux \fi \def \close {\if
\empty \UndefLabels \else \message {*** There were undefined labels ***} \medskip \noindent ****************** \ Undefined Labels: \tt \par \UndefLabels
\fi \if \TRUE \auxwrite \closeout 1 \fi \par \vfill \supereject \end }  \def \Caixa #1{\setbox 1=\hbox {$#1$\kern 1pt}\global \edef \tamcaixa {\the \wd 1}\box 1} \def \caixa #1{\hbox to \tamcaixa {$#1$\hfil }} \def
\medprod {\mathop {\mathchoice {\hbox {$\mathchar "1351$}}{\mathchar "1351}{\mathchar "1351}{\mathchar "1351}}}  \def \medcup {\mathop {\mathchoice {\raise 1pt \hbox {$\mathchar
"1353$}}{\mathchar "1353}{\mathchar "1353}{\mathchar "1353}}} \def \medcap {\mathop {\mathchoice {\raise 1pt \hbox {$\mathchar "1354$}}{\mathchar
"1354}{\mathchar "1354}{\mathchar "1354}}}        \def
\clauses #1{\def \crr {\vrule width 0pt height 10pt depth 5pt\cr }\left \{ \matrix {#1}\right .} \def \cl #1 #2 #3 {#1, & \hbox {#2 } #3\hfill \crr }
\def \paper #1#2#3{ \hsize #1truemm   \advance \hsize by -#3truemm  \advance \hsize by -#3truemm \vsize #2truemm
\advance \vsize by -#3truemm  \advance \vsize by -#3truemm \hoffset =-1truein \advance \hoffset by #3truemm \voffset =-1truein \advance \voffset
by #3truemm }  \def \FPCovar {2.6} \def \FPIntroEpsilon {12.23}  \def \FPBackInvarCorol {13.4}
\def \FPBigNewTightResult {16.18} \def \FPBoundedSets {13.2} \def \FPCovarPhi {17.4} \def \FPDefIdempSLA {3.11} \def \FPDeltaMatchCosntru {10.13}
\def \FPEquality {7.21} \def \FPFRFstarRTwo {10.12} \def \FPFormOfHull {7.13} \def \FPFormOfSubSets {7.15} \def \FPForwardInvariance {10.19} \def
\FPIntroPhiSigmaNonDeg {16.1} \def \FPIntroRelTight {15.8} \def \FPMapForISG {12.22} \def \FPNonOpenUltra {17.11} \def \FPOpenStringsInvariant {12.24}
\def \FPRangeInsv {7.19} \def \FPRelTightClosed {15.9} \def \FPStringDirected {10.1.iii} \def \FPStringDivisors {10.2} \def \FPSigmaPhiOnOpen {16.15}
\def \FPBirthOfString {17.7} \def \FPDualBackOnStrings {17.6} \def \FPFRFstarR {10.10} \def \FPIntroStarAction {10.9} \def \FPPropsOpenString {11.2}
\def \FPOpenInvarUnderTheta {11.4} \def \FPFactorlemmaTwo {15.11}     \def \MaxOpenOrDeadEnd {11.3} \def \first #1{\cite [#1]{ESOne}} \def \supp {\mathbox {supp}} \def \Ker {\hbox {Ker}} \def \sla
{semilattice} \def \interior #1{\mathaccent '27#1} \def \longhookrightarrow {\lhook \joinrel \longrightarrow } \def \Imply {\ \mathrel {\Rightarrow }\ }
\def \arw #1{\ {\buildrel #1 \over \longrightarrow }\ } \def \kl #1 / #2 / #3/{#1, & \hbox {#2 } #3\hfill \crr } \def \|{\mathrel {|}} \def \itmproof #1
{\medskip \noindent #1\enspace } \def \dualrep {\hat \theta } \def \Iff {\ \mathrel {\Leftrightarrow }\ } \def \Bool #1{\left [#1\right ]} \def \id {\hbox
{id}} \def \ac #1#2{\dualrep _{#1}(#2)} \def \acinv #1#2{\dualrep _{#1}\inv (#2)} \def \G {{\cal G}} \def \I {{\cal I}} \def \S {{\cal S}} \def \P {{\cal
P}} \def \E {{\cal E}} \def \tS {\tilde S} \font \smallgothicfont = eufb10 scaled 833 \font \gothicfont = eufb10 \def \goth #1{\hbox {\gothicfont #1}}
\def \hull {\goth H} \def \ehull {\goth E} \def \Ehat {\goth E\kern -5pt \widehat {\vrule height 6.5pt width 5pt depth -10pt}} \def \sinfShort {\Ehat
_\infty } \def \stightShort {\Ehat \sub {tight}} \def \smaxShort {\Ehat \sub {max}} \def \sessShort {\Ehat \sub {ess}} \def \X {{\cal X}} \def \SX {S_\X
} \def \tSX {\tilde \SX } \def \TX {{\cal T}_\X } \def \tildeTX {\tilde {\cal T}_\X } \def \MX {{\cal M}_\X } \def \OX {{\cal O}_\X } \def \LX {L_\X }
\def \sub #1{_{\hbox {\sixrm #1}}} \def \gothEhat {\goth E\kern -5pt \widehat {\vrule height 6.5pt width 5pt depth -10pt}} \newcount \currentgrouplevel

\centerline {\bf SUBSHIFT SEMIGROUPS} \medskip \centerline {\tensc R. Exel\footnote {$^*$}{\eightrm Partially supported by CNPq.} and B. Steinberg\footnote
{$^{**}$}{\eightrm Partially supported by the Fulbright Commission.}} \vskip 1cm \footnote {}{\eightrm August 21,  2019.}

\def \fin {{\hbox {\sixrm fin}}}

\Abstract Given a one-sided subshift $\X $ on a finite alphabet, we consider the semigroup $\SX =\LX \cup \{0\}$, where $\LX $ is the language of $\X $,
equipped with the multiplication operation given by concatenation, when allowed, and set to vanish otherwise.  We then study the inverse hull $\hbox
{\smallgothicfont H}(\SX )$, relating it with C*-algebras that have been discussed in the literature in association with subshifts.  \endAbstract

\section Introduction

Since Cuntz and Krieger associated $C^*$-algebras to Markov shifts in \cite {CK}, there has been a fascinating interaction between operator theory
and symbolic dynamics.  Often dynamical aspects of the system are reflected to some extent in the operator algebra (perhaps considered jointly with
its associated diagonal subalgebra) and many natural dynamical invariants, like the Bowen-Franks group of a Markov shift, have operator theoretic
interpretations (in this case as the $K_0$-group of the Cuntz-Krieger algebra).

In the list of references at the end of this paper the reader will find a (non-exhaustive) list of papers, too numerous to be cited here, mostly by
T.~M.~Carlsen, K.~Matsumoto and several co-authors, dealing with this rich interplay.

There have been a number of attempts over the years to associate $C^*$-algebras to general subshifts.  Recall that if $\Sigma $ is a finite alphabet,
then a subshift $\X $ is a shift-invariant closed subspace of $\Sigma ^{\bf N}$.  Matsumoto made a first attempt to define a $C^*$-algebra associated
to an arbitrary subshift \cite {MatsuOri} and proved a number of results about these algebras \cite {MatsumotoDimension}, \cite {MatsuAuto}, \cite
{MatsumotoStabilized}, before some mistakes were discovered and it became clear that there were two different algebras: $\MX $ (which we call the
Matsumoto algebra) and $\OX $ (which is usually called the Carlsen-Matsumoto algebra).  The latter is generally accepted as the ``correct'' algebra for
shift spaces.  See \cite {MatsuCarl} for more on the history of these algebras.  It should be noted that for Markov shifts, both the Matsumoto and the
Carlsen-Matsumoto algebras are just the usual Cuntz-Krieger algebra associated to the shift.

One of the difficulties in finding the correct algebra for shifts has been due to the lack of a systematic approach to the construction.  For a shift
of finite type, like a Markov shift, the shift map is a local homeomorphism and so one can build a Deaconu-Renault groupoid from the dynamics \cite
{cuntzlike}.  The $C^*$-algebra of this groupoid is then the Cuntz-Krieger algebra of the shift.  For more general subshifts, the shift map is not a
local homeomorphism and so if one wants to build a Deaconu-Renault groupoid, one leaves the realm of etale groupoids for that of semi-etale groupoids
\cite {CarlsenThomsen}, \cite {Thomsen}, which is less familiar.

This paper bears to fruition a more systematic approach for constructing both the Matsumoto and Carlsen-Matsumoto algebras using the scheme of building
from a combinatorial object an inverse semigroup and then an etale groupoid as per  \cite {actions}.  Here our basic combinatorial object is the semigroup
of admissible factors associated to a shift.  If $\X $ is a shift, then $\LX $ is the language of all finite words that appear in an element of $\X $.
There is a natural semigroup structure on $\SX =\LX \cup \{0\}$ where the product is concatenation when the result is an admissible word, and otherwise
the product is $0$.  This semigroup was first introduced long ago by Hedlund and Morse \cite {HeadlundMorse}, who used the semigroup associated to the
Thue-Morse dynamical system to give the first example of a finitely generated, infinite semigroup in which each element is nilpotent.

Over the past several years, a number of mathematicians have associated $C^*$-algebras to left-cancellative semigroups, or more generally categories
\cite {BedosSpielberg},  \cite {SpielbergA}, \cite {SpielbergB}.  For example, graph $C^*$-algebras (which for finite graphs are basically the same
thing as Cuntz-Krieger algebras) are the $C^*$-algebras associated to the path category (or free category) of the graph. However, our semigroup $\SX $
is far from cancellative because of $0$ and except in the case of Markov shifts it fails the kind of partial associativity you would want for it to come
from some category-like structure.  Our fundamental observation, which led to our previous paper \cite {ESOne}, is that if a semigroup $S$ with $0$ is
$0$-left cancellative, meaning that (like in an integral domain) $xy=xz\neq 0$ implies $y=z$, then one can associate a natural inverse semigroup $\hull
(S)$ to $S$, called its inverse hull, which generalizes the constructions commonly used for left cancellative semigroups and categories.  Moreover,
the semigroup $\SX $ coming from a subshift is $0$-cancellative (it satisfies the cancellation condition on both sides) and so we can build the desired
inverse semigroup.

We should point out that, for left-cancellative semigroups (without zero), the inverse hull has already been considered in relation to C*-algebras by X. Li
\cite [Chapter 5]{CELY}, this being in fact a significant source of inspiration for our work both in \cite {ESOne} and in the present paper.  In relation
to Li's work, the main difficulty we face here is related to the fact that our semigroups are allowed to have zero (a crucial requirement if one is to
include language semigroups), and hence left-cancellativity never holds. We must therefore  make do with the weaker assumption of $0$-left-cancellativity.

It turns out that the universal groupoid of the inverse hull $\hull (\SX )$ of $\SX $ has a number of natural closed invariant subspaces.  Previous to our
work, the so-called tight spectrum (see \cite {actions}) has been the most important invariant subspace of the universal groupoid of an inverse semigroup.
But it turns out that the tight groupoid of $\hull (\SX )$ leads to an algebra that has not yet been studied.  The Carlsen-Matsumoto algebra $\OX $
comes from considering the reduction to a closed invariant subspace of the tight spectrum coming from the closure of the ultrafilters corresponding
to infinite words in the shift.  The Matsumoto algebra comes from a larger closed invariant subspace, which is incomparable with the tight spectrum,
that we call the essentially tight spectrum.

One advantage of using the inverse semigroup approach is that one easily checks that $\hull (\SX )$ is a strongly $0$-$E$-unitary inverse semigroup with
universal group the free group. The general machinery of \cite {MilanSteinberg} then gives you for free that the $C^*$-algebra is a partial crossed
product of the free group with a commutative $C^*$-algebra.  Moreover, there is a natural way to define an analogue of the shift on the full spectrum
of $\hull (\SX )$, which can then be used to identify the universal groupoid with a Deaconu-Renault groupoid, yielding amenability for this groupoid
and all its reductions.  In fact, we show more generally that for a class of free group partial actions, called semi-saturated and orthogonal, the
corresponding partial transformation groupoid is always a Deaconu-Renault groupoid, which may be of interest in its own right.

The systematic approach via $0$-left-cancellative semigroups also leaves open the possibility of studying higher rank shift spaces, built from forbidding
patterns in higher rank graphs.

Some of the main results given here were announced in \cite {announce}.

\section Essentially tight representations

Let, for the time being, $\E $ be a {\sla } with zero and let $\Omega $ be any set.  Also let $$ \pi :\E \to \P (\Omega ) $$ be a representation of $\E
$ in the Boolean algebra $\P (\Omega )$ formed by all subsets of $\Omega $.  Recall from \cite [11.6]{actions} that $\pi $ is said to be \emph {tight}
if, for every finite subsets $X,Y\subseteq \E $, and for every finite cover $Z$ for $$ \E ^{X,Y} := \{z\in \E : z\leq x,\ \forall x\in X,\hbox { and }
z\perp y,\ \forall y\in Y\}, $$ one has that $$ \medcup _{z\in Z}\pi (z) = \pi (X, Y), $$ where $$ \pi (X, Y):= \Big (\medcap _{x\in X} \pi (x) \Big )
\cap \Big ( \medcap _{y\in Y} \Omega \setminus \pi (y)\Big ).  $$

Notice that, according to the definition of covers \cite [11.5]{actions}, one has that $Z\subseteq \E ^{X,Y}$, so if follows from \cite [11.3]{actions}
that $$ \medcup _{z\in Z}\pi (z) \subseteq \pi (X, Y), $$ regardless of whether or not $\pi $ is tight.  We may thus think of $$ \pi (X, Y) \ \setminus
\ \medcup _{z\in Z}\pi (z) \equationmark DefectSet $$ as the \emph {defect} set relative to the cover $Z$ of $\E ^{X,Y}$, so that tightness means that
all possible defect sets are empty.

\definition Let $\pi $ be a representation of $\E $ in $\P (\Omega )$.  We will say that $\pi $ is \emph {essen\-tial\-ly tight} if, for every finite
subsets $X,Y\subseteq \E $, and for every finite cover $Z$ for $\E ^{X,Y}$, one has that the corresponding defect set \ref {DefectSet} is finite.  In case
$\E $ is itself a sub-{\sla } of $\P (\Omega )$, we will say that $\E $ is an \emph {(essentially) tight sub-{\sla }} if the identity representation
$\E \to \P (\Omega )$ is (essentially) tight.

A similar concept, which we will not use here, could be defined if instead of $\P (\Omega )$, our representation $\pi $ took values in an arbitrary
Boolean algebra equipped with an ideal whose elements one sees as negligible.

\def \Q {{\cal Q}}

Let $\P _\fin (\Omega )$ be the Boolean ideal of $\P (\Omega )$ formed by all finite subsets of $\Omega $.  The quotient Boolean algebra $$ \Q (\Omega
) := \P (\Omega )/\P _\fin (\Omega ) $$ is defined as the quotient of $\P (\Omega )$ by the equivalence relation $$ X\sim Y \iff X\mathop {\Delta }Y
\hbox { is finite}, $$ where the \emph {symmetric difference} between $X$ and $Y$ is defined by $$ X\mathop {\Delta }Y = (X\setminus Y)\cup (Y\setminus X).  $$

The quotient mapping will henceforth be denoted by $$ q:\P (\Omega )\to \Q (\Omega ).  \equationmark IntroQuotientBoole $$ Employing this terminology
notice that a representation $\pi :\E \to \P (\Omega )$ is essentially tight if and only $q\circ \pi $ is a tight representation of $\E $ in $\Q (\Omega )$.

As an example, observe that the identity mapping $\P (\Omega )\to \P (\Omega )$ is tight by \cite [11.9]{actions}, so $\P (\Omega )$ is a tight sub-{\sla
} of itself.  On the other hand, when $\Omega =\{0, 1,2\}$, the sub-{\sla } consisting of the subsets $$ \emptyset ,\ \{0\},\ \{1\},\ \{0,1,2\}, $$
is not tight, since $\big \{\{0\}, \{1\}\big \}$ is a cover for $\{0,1,2\}$, with nonempty defect set $\{2\}$.  This, and all other defect sets are
evidently finite, so the above is an example of an essentially tight sub-{\sla }.

An example of a sub-{\sla } which is not essentially tight may easily be built based on the above idea: take $\Omega ={\bf N}$, and consider the
sub-{\sla } consisting of the subsets $$ \emptyset ,\ \{0\},\ \{1\},\ {\bf N}.  $$ Then $\big \{\{0\}, \{1\}\big \}$ is a cover for ${\bf N}$, with
infinite defect set $\{2,3,\ldots \}$.

The following concept is a variation of the notion of relatively tight characters introduced in \first {\FPIntroRelTight }.

\definition \label IntroEssTightChar Let $\E $ be a {\sla } and let $\pi $ be a representation of $\E $ on a set $\Omega $.  Given a character $\varphi
$ on $\E $, we will say that $\varphi $ is \emph {essentially tight relative to $\pi $}, or simply \emph {$\pi $-essentially tight}, if $\varphi $
is tight relative to the representation $q\circ \pi $ of\/ $\E $ in $\Q (\Omega )$.

Combining an essentially tight representation with an essentially tight character, one gets a (truly) tight character:

\state Proposition \label TwoEssTight Let $\E $ be a {\sla } and let $\pi $ be an essentially tight representation of\/ $\E $ on a set $\Omega $.
Then every $\pi $-essentially tight character on $\E $ is tight.

\Proof Let $\varphi $ be any $\pi $-essentially tight character.  Since $\varphi $ is nonzero, it satisfies condition (i) of \cite [Lemma 11.7]{actions}
so, in order to prove that $\varphi $ is tight, it suffices to verify the simplified tightness condition of \cite [Proposition 11.8]{actions}.  We thus
pick $x$ in $\E $, and a finite covering $\{y_1, y_2,\ldots , y_n\}$ of $x$, and our task is to prove that $$ \varphi (x) = \bigvee _{i=1}^n \varphi
(y_i).  \equationmark SimplifiedPhiTight $$ By the essential tightness of $\pi $, the defect set $$ \pi (x)\setminus \medcup _{i=1}^n\pi (y_i) $$
is finite and consequently $$ q\big (\pi (x)\big ) = \bigvee _{i=1}^nq\big (\pi (y_i)\big ), $$ so \ref {SimplifiedPhiTight} follows from the $\pi
$-essential tightness of $\varphi $.  \endProof

We would now like to prove the converse of the above result, but there are issues we need to discuss first.  The main outstanding question is whether
or not the simplified tightness condition of \cite [Proposition 11.8]{actions} is enough to check tightness of $q\circ \pi $.

Since the proof of the simplified condition is already within reach, let us tackle it now, postponing the question of essential tightness of $\pi $
for a while.

\state Lemma \label LemmaEssTight Let $\E $ be a {\sla } and let $\pi $ be a representation of $\E $ on a set $\Omega $.  Suppose that every $\pi
$-essentially tight character on $\E $ is tight.  Then, for every $x$ in $\E $, and for every finite covering $\{y_1, y_2,\ldots , y_n\}$ of $x$,
one has that $$ q\big (\pi (x)\big ) = \bigvee _{i=1}^nq\big (\pi (y_i)\big ).  $$

\Proof Suppose by way of contradiction that the equality displayed above is false, whence the defect set $$ D:= \pi (x) \ \setminus \ \medcup _{i=1}^n\pi
(y_i) $$ is infinite.  Since $D$ is not in $\P _\fin (\Omega )$, it follows that $q(D)$ is nonzero, so there exists a character $\chi $ of $\Q (\Omega
)$, preserving meets and joins, such that $\chi \big (q(D)\big )=1$.

Considering the composition $$ \varphi : \E \arw \pi \P (\Omega ) \arw q \Q (\Omega ) \arw \chi \{0, 1\}, $$ notice that all arrows involved are {\sla
} homomorphisms, and hence so is $\varphi $.  In order for $\varphi $ to qualify as a character we must only show that $\varphi $ is nonzero, but this
follows from the fact that $D\subseteq \pi (x)$, and hence $$ \varphi (x) = \chi (q\big (\pi (x))\big ) \geq \chi \big (q\big (D)\big ) =1.  $$

We next claim that $\varphi $ is $\pi $-essentially tight.  For this, suppose that $z,w_1,\ldots ,w_n$ are elements of $\E $ such that $ q\big (\pi (z)\big
) = \bigvee _{i=1}^n q\big (\pi (w_i)\big ).  $ Because $\chi $ preserves joins we deduce that $$ \varphi (z)= \chi \big (q\big (\pi (z)\big )\big )=
\bigvee _{i=1}^n \chi \big (q\big (\pi (w_i)\big )\big ) = \bigvee _{i=1}^n \varphi (w_i), $$ thus proving that $\varphi $ is indeed $\pi $-essentially tight.

The hypothesis therefore implies that $\varphi $ is tight, but we will reach a contradiction by showing that this is not so, hence concluding the proof
of the statement.

Observing that $D\cap \pi (y_i)=\emptyset $, for all $i$, we have that $$ 0 = \chi \big (q(D\cap \pi (y_i))\big ) = \chi \big (q\big (D\big )\big )
\chi \big (q\big (\pi (y_i)\big )\big ) = \varphi (y_i), $$ so $\bigvee _{i=1}^n \varphi (y_i)=0$, but since $\varphi (x)=1$, as shown above, we see
that $\varphi $ fails to verify the tightness condition relative to the cover $\{y_1, y_2,\ldots , y_n\}$ of $x$.  This concludes the proof.  \endProof

Adding a few assumptions to our representation $\pi $, we may now obtain a necessary and sufficient condition for essential tightness.

\state Proposition \label NecSufEssTight Let $\E $ be a {\sla } and let $\pi $ be a representation of $\E $ on a set $\Omega $.  Suppose that at least
one of the following two conditions hold: \iaitem \aitem there is a finite subset $\{x_1, \ldots , x_n\}\subseteq \E $, such that $\Omega \setminus
\medcup _{i=1}^n\pi (x_i)$ is finite, or \aitem $\E $ admits no finite cover.  \medskip \noindent Then the following are equivalent: \iItemize \iItem
$\pi $ is an essentially tight representation, \iItem every $\pi $-essentially tight character on $\E $ is tight.

\Proof That (i)$\Imply $(ii) is precisely the content of \ref {TwoEssTight}, while the converse follows immediately from \cite [Proposition 11.8]{actions}
and \ref {LemmaEssTight}.  \endProof

It is our next main goal in this section to study $0$-left-cancellative semigroups $S$ for which the sub-{\sla } $\ehull (S) \subseteq \P (S')$ is
essentially tight.

\state Lemma \label CoversAndLengths Let $S$ be a $0$-left-cancellative semigroup admitting least common multiples and let $\ell $ be a homogeneous
\first {13.2.ii} $N$-valued length function for $S$, where $N$ is a given totally ordered set.  Given $X$ in $\ehull (S) $, let $\{Y_1,Y_2,\ldots ,Y_n\}$
be a cover for $X$.  Then there exists some $n_0$ in $N$ such that $$ \interior X\cap \big \{s\in S : \ell (s)>n_0\}\subseteq \medcup _{i=1}^n Y_i, $$
where the \emph {interior} of $X$ is defined by $$ \interior X = \{s\in X: \exists p\in S,\ sp\in X\}.  $$

\Proof By \first {\FPFormOfHull } we may write $X = uF_\Lambda $, where $\Lambda $ is a finite subset of $\tS $, with $\Lambda \cap S\neq \emptyset $,
and $u\in \Lambda $.  In addition, since each $Y_i\subseteq X$, we may use \first {\FPFormOfSubSets } to write $Y_i =ux_iF_ {\Delta _i}$, where each
$\Delta _i$ is a finite subset of $\tS $, each $x_i$ is in $\tS $, and $\Lambda x_i\subseteq \Delta _i$.

Setting $\ds m_0=\max _{1\leq i\leq n}\ell (x_i)$, let $$ A=\{t\in F_u: \ell (t)\leq m_0\}, $$ so that $A$ is evidently a bounded subset of $F_u$,
and hence $uA$ is bounded because $\ell $ is homogeneous.  Therefore there is some $n_0$ in $N$ such that $\ell (ut)\leq n_0$, for every $t$ in $A$.

Picking any $s$ in $\interior X$, with $\ell (s)>n_0$, we will conclude the proof by showing that $s$ lies in $\medcup _{i=1}^n Y_i.$ Given that
$X=uF_\Lambda $, we may write $s=ur$, with $r\in F_\Lambda $.  Notice that $$ \ell (r)>m_0, $$ since otherwise $r$ would lie in $A$, and this would
imply that $\ell (s) = \ell (ur) \leq n_0$.  We next claim that $$ F_{\Lambda r}\neq \emptyset .  $$ In order to see this, choose $p$ in $S$ such that
$sp\in X$, and write $sp=uq$, for some $q$ in $F_\Lambda $.  We then have that $uq=sp=urp$, whence $q=rp$, and then for every $t$ in $\Lambda $, $$
0 \neq tq = trp, $$ proving that $p\in F_{\Lambda r}$, and verifying our claim.  It then follows that $urF_{\Lambda r}$ is a nonempty subset of $X$.
Since $\{Y_1,Y_2,\ldots ,Y_n\}$ covers $X$, there is some $i$ such that $urF_{\Lambda r}$ intersects $Y_i$.  So, picking $t$ in $urF_{\Lambda r}\cap
Y_i$, we may write $$ t =urx=ux_iy, $$ with $x\in F_{\Lambda r}$, and $y\in F_ {\Delta _i}$.  This implies that $ rx=x_iy, $ so $x_i\| rx$, and since
$$ \ell (x_i)\leq m_0<\ell (r), $$ we have that $x_i\|r$.  Consequently we may write $r=x_iz$, for some $z$ in $S$, and then $$ x_izx = rx = x_iy, $$
so $ zx = y\in F_ {\Delta _i}.  $ Since $F_ {\Delta _i}$ is easily seen to be hereditary with respect to the division relation, we see that $z\in F_
{\Delta _i}$, as well, whence $$ s = ur = ux_iz \in ux_iF_ {\Delta _i} = Y_i.  $$ This concludes the proof. \endProof

\definition \label DefLocFin We will say that an $N$-valued length function $\ell $ is \emph {locally finite} provided $$ \{s\in S': \ell (s)\leq n\}
$$ is a finite set for every $n\in N$.

The following is the achievement of the goal we set ourselves earlier.

\state Theorem \label GeneralEssTight Let $S$ be a $0$-left-cancellative semigroup admitting least common multiples and suppose that: \iItemize \iItem
$S$ possesses a locally finite, homogeneous length function $\ell $, \iItem there is a finite collection $s_1, s_2,\ldots ,s_n$ of elements of $S$
such that $ S'\setminus \medcup _{i=1}^n E^\theta _{s_i} $ is finite, \iItem the \emph {boundary} $$ \partial X:= X\setminus \interior X $$ of every
constructible set $X$ is finite.  \medskip \noindent Then $\ehull (S) $ is an essentially tight sub-{\sla } of $\P (S')$.

\Proof Considering the ideal $\P _\fin (S')$ of $\P (S')$ formed by all finite subsets, we must prove that the composition $$ \beta \ :\ \ehull (S)
\ \longhookrightarrow \ \P (S') \ \longrightarrow \ \P (S')/\P _\fin (S'), $$ is tight, where the arrow in the left-hand-side is the inclusion map,
and the right-hand-side arrow is the quotient map.

By (ii) one has that $\beta $ satisfies \cite [11.7.i]{actions}, so we may use \cite [11.8]{actions} to reduce our task to proving that $$ \beta
(X)=\bigvee _{i=1}^n\beta (Z_i), \equationmark SimplifiedTight $$ whenever $X\in \ehull (S) $, and $\{Z_1,Z_2,\ldots ,Z_n\}$ is a cover for $X$.
Employing \ref {CoversAndLengths}, let $n_0$ be such that $$ \interior X\cap \big \{s\in S : \ell (s)>n_0\}\subseteq \medcup _{i=1}^n Z_i, $$ so the
corresponding defect set satisfies $$ X\setminus \medcup _{i=1}^n Z_i \ \subseteq \ \partial X \cup \Big (\interior X\setminus \medcup _{i=1}^n Z_i\Big )
\ \subseteq \ \partial X \cup \big \{s\in S : \ell (s)\leq n_0\}, $$ which is a finite set because $\partial X$ is finite and $\ell $ is locally finite.
This evidently implies \ref {SimplifiedTight}, so the proof is concluded.  \endProof

\section Subshift semigroups

\label SubshiftSgrSection By a \emph {subshift} on a finite alphabet $\Sigma $ one means a subset $\X \subseteq \Sigma ^{\bf N}$, which is closed
relative to the product topology, and invariant under the left shift map $$ x_1x_2x_3\ldots \ \mapsto \ x_2x_3x_4\ldots $$

Given any such $\X $, the \emph {language} of $\X $ is the set $\LX $ formed by all finite words appearing as a block in some infinite word belonging
to $\X $.  We will not allow the empty word in $\LX $, as sometimes done in connection with subshifts, so all of our words have strictly positive length.

Summarizing, a finite sequence $a_1a_2\ldots a_n$, with $a_i\in \Sigma $, lies in $\LX $, if and only if $n\geq 1$, and there exists an element
$x=x_1x_2x_3\ldots \in \X $, and some $k\geq 0$, such that $$ a_i=x_{k+i}, \for i=1,\ldots ,n.  $$

Throughout this paper we will be concerned with the semigroup $$ \SX =\LX \cup \{0\}, \equationmark IntroSX $$ equipped with the multiplication operation
given by $$ \def \quad { } \mu \nu =\left \{\matrix { \mu \cdot \nu , & \hbox { if $\mu ,\nu \neq 0$, and $\mu \nu \in L$}, \cr \pilar {12pt} 0, &
\hbox { otherwise, }\hfill }\right .  $$ where $\mu \cdot \nu $ stands for the concatenation of $\mu $ and $\nu $.

Incidentally, this is the semigroup described in \first {6.2}, where we did not insist that $\X $ be closed relative to the product topology of $\Sigma
^{\bf N}$.  However in our present development it will be crucial that $\X $ be closed, and that property will manifest itself through the following condition.

\state Proposition \label ClosedShifts Let $\Sigma $ be a finite alphabet and let $\X \subseteq \Lambda ^{\bf N}$ be any subset.  Let $P$ be the
collection of all finite prefixes\fn {When $\X $ is invariant under the left shift notice that $P=\LX $.} of the members of $\X $.  Then the following
conditions are equivalent: \iItemize \iItem $\X $ is closed relative to the product topology of\/ $\Sigma ^{\bf N}$, \iItem for every infinite word
$\omega \in \Lambda ^{\bf N}$, if all prefixes of $\omega $ lie in $P$, then $\omega \in \X $.

\Proof (i)$\Imply $(ii) \enspace Given that all prefixes of $\omega $ lie in $P$, we may choose, for each $n$ in ${\bf N}$, an element $\omega _n\in \X $
sharing with $\omega $ the prefix of length $n$.  It therefore follows that $\omega _n\to \omega $ in the product topology, so $\omega $ lies in $\X $
because $\X $ is closed.

\itmproof {(ii)$\Imply $(i)} Let $\omega $ be a point in the closure of $\X $.  Given any integer $n$, let $\alpha $ be the prefix of $\omega $ of
length $n$.  The set $A$ formed by all members of $\Sigma ^{\bf N}$ having $\alpha $ as a prefix is know as a \emph {cylinder}, and it is well known
that $A$ is an open set in the product topology.

Since $\omega $ lies in $A$, we have that $A\cap \X \neq \emptyset $.  Picking $\eta $ in that intersection we see that $\alpha $ is a prefix of $\eta $,
whence $\alpha $ lie in $P$.  This shows that all prefixes of $\omega $ lie in $P$, whence $\omega \in \X $, by (ii).  \endProof

\fix Throughout this chapter we will let $\X $ be a fixed subshift. \medskip

When considering the unitized semigroup $\tilde \SX = \SX \cup \{1\}$, the added unit may sometimes be thought of as the empty word $\varnothing $ but,
as already mentioned, $\varnothing $ will not be allowed in the language $\LX $.

Given $\mu $ and $\nu $ in $\SX $, with $\nu $ nonzero, notice that $\mu \|\nu $ if and only if $\mu $ is a prefix of $\nu $.  Given that divisibility
is also well defined for the unitized semigroup $\tilde \SX $ \first {5.4}, and given that $1$ divides any $\mu \in \SX $, we will also say that $1$
is a prefix of $\mu $.

Some of the special properties of $\SX $ of easy verification are listed below:

\state Proposition \label EasyProp \iItemize \iItem $\SX $ is $0$-left-cancellative and $0$-right-cancellative, \iItem for every $s$ and $t$ in $\SX $,
one has that either $s\|t$, or $t\|s$, or else there is no common multiple of $s$ and $t$ besides $0$, \iItem $\SX $ has no nonzero idempotent elements.

A further special property of $\SX $ is a very strong uniqueness of the normal form \first {7.13 \& 7.21} for elements in $\hull (\SX ) $, as we shall
see next:

\state Proposition \label SubshiftUniqueness For $i=1,2$, let $\Lambda _i$ be a finite subset of $\tSX $ intersecting $\SX $ non-trivially, and let
$u_i,v_i\in \Lambda _i$ be such that $$ \theta _{u_1}f_{\Lambda _1} \theta _{v_1}\inv = \theta _{u_2}f_{\Lambda _2} \theta _{v_2}\inv \neq 0.  $$
Then $u_1=u_2$, $v_1=v_2$, and $F_{\Lambda _1} = F_{\Lambda _2}$.

\Proof Letting $x_1$ and $x_2$ be as in \first {\FPEquality }, notice that by \first {\FPEquality .i}, the range of $\theta _{u_i}f_{\Lambda _i} \theta
_{v_i}\inv $ is contained in $E_{u_ix_i}$.  Using \first {\FPRangeInsv } we then deduce that $F_{\Lambda _i}\subseteq E_ {x_i}$.

By hypothesis $F_{\Lambda _i}$ is nonempty so we may pick some $\mu $ there.  Observing that $F_{\Lambda _i}$ is hereditary, we have that the first
letter of $\mu $, say $\mu _1$, is also in $F_{\Lambda _i}$, and hence in $E_ {x_i}$.  However, unless $x_i=1$, every element of $E_ {x_i}$ has length
at least $2$.  Therefore $x_i=1$, and then we have by \first {\FPEquality .ii} that $u_1=u_2$, and $F_{\Lambda _1} = F_{\Lambda _2}$.

Observing that since \first {\FPEquality .b\&c} speak of nonzero idempotents, the only valid alternative is \first {\FPEquality .a}, from where it
follows that $v_1=v_2$.  \endProof

Given the importance of strings \first {10.1}, let us give an explicit description of these in the present context.

\state Proposition \label DescribeStringsAsWords Given a (finite or infinite) word $$ \omega = \omega _1\omega _2\omega _3\ldots , $$ on the alphabet
$\Sigma $, assume $\omega $ to be \emph {admissible} (meaning that $\omega $ belongs to $\LX $, if finite, or to $\X $, if infinite) and consider the
set $\sigma _\omega $ formed by all prefixes of $\omega $ having positive length, namely $$ \sigma _\omega =\{\omega _1,\ \omega _1\omega _2,\ \omega
_1\omega _2\omega _3,\ \ldots \ \}.  $$ Then: \iItemize \iItem $\sigma _\omega $ is a string, \iItem $\sigma _\omega $ is an open string if and only if
$\omega $ is an infinite word, \iItem $\sigma _\omega $ is a maximal string if and only if $\omega $ is an infinite word, \iItem for any string $\sigma $
in $\SX $, there exists a unique admissible word $\omega $ such that $\sigma =\sigma _\omega $.

\Proof Before we begin we should note that, in case $\omega $ is a finite admissible word, one has that $\omega $ represents an element of $\SX $,
in which case $\sigma _\omega $ has already been defined in \first {\FPStringDivisors }.  However there is no conflict since the above definition of
$\sigma _\omega $ agrees with the meaning given to this notation by \first {\FPStringDivisors }.

We leave the easy proofs of (i), (ii) and (iii) to the reader and concentrate on proving (iv).  Observe that, given any $\mu $ and $\nu $ in $\SX $,
admitting a nonzero common multiple, and such that $\ell (\mu )\leq \ell (\nu )$, one necessarily has that $\mu \|\nu $.  Thus, if we order the elements
of a given string according to length, we will have also ordered them according to divisibility.  The conclusion is then evident.  \endProof

Strings are one of our best instruments to provide characters on $\ehull (\SX ) $. Now that we have a concrete description of strings in terms of
admissible words, let us give an equally concrete description of the characters induced by strings.  In the result below we will focus on infinite words,
but a similar result for finite words is easy to deduce from \first {\FPDeltaMatchCosntru }.

\state Proposition \label DescribeCharsAsWords Let $\omega $ be a given infinite admissible word, and let $X$ be any $\theta $-constructible set, written
in normal form, namely $X=uF^\theta _\Lambda $, where $\Lambda $ is a finite subset of $\tSX $, intersecting $\SX $ nontrivially, and $u\in \Lambda $.
Regarding the string $\sigma _\omega $, and the associated character $\varphi _{\sigma _\omega }$ \first {16.1.ii}, the following are equivalent:
\iItemize \iItem $\varphi _{\sigma _\omega }(X)=1$, \iItem $u$ is a prefix of $\omega $, and upon writing $\omega =u\eta $, for some infinite word
$\eta $, one has that $t\eta $ is admissible (i.e.~lies in $\X $), for every $t$ in $\Lambda $.

\Proof Assuming that $\varphi _{\sigma _\omega }(X)=1$, we have that $$ \sigma _\omega \in \varepsilon (X) = \varepsilon (uF^\theta _\Lambda ) = \theta
^\star _u(F^\star _\Lambda ), $$ which is equivalent to $$ \emptyset \neq \sigma _\omega \cap E^\theta _u \subseteq \theta _u(F^\theta _\Lambda ), $$
by \first {\FPFRFstarRTwo }.  Any element in $\sigma _\omega \cap E^\theta _u$ has $u$ as a prefix and is itself a prefix of $\omega $, so one sees that
$u$ is necessarily a prefix of $\omega $.  Writing $\omega =u\eta $, where $\eta $ is an infinite word, observe that $$ \{u\eta _1\eta _2\ldots \eta _k:
k\geq 1\} = \sigma _\omega \cap E^\theta _u \subseteq \theta _u(F^\theta _\Lambda ).  $$ Applying $\theta _u \inv $ to the left and rightmost terms above
then leads to $$ \{\eta _1\eta _2\ldots \eta _k: k\geq 1\} \subseteq F^\theta _\Lambda , $$ from where we conclude that $t\eta _1\eta _2\ldots \eta _k$
is admissible for every $t$ in $\Lambda $, and for every $k\geq 1$, so that all prefixes of the infinite word $t\eta $ are admissible, and hence $t\eta $
lies in $\X $ by \ref {ClosedShifts}.

The converse is proven by essentially running the above argument backwards.  \endProof

In case the set $X$ of the above result coincides with $F^\theta _\mu $, for some $\mu $ in $\LX $, we get the following simplification:

\state Proposition \label DescribeCharsAsWordsOnFollower Let $\omega $ be an infinite admissible word, and let $\mu \in \LX $.  Then $$ \varphi _{\sigma
_\omega }(F^\theta _\mu )=1 \IFF \mu \omega \in \X .  $$

\Proof This follows at once by letting $u=1$ and $\Lambda =\{1, \mu \}$ in \ref {DescribeCharsAsWords}.  \endProof

Let us now study a situation in which infinite words provide all ultra-characters.

\state Proposition \label GroundForSubshifts The following are equivalent: \iItemize \iItem for every finite subset $\Lambda \subseteq \tSX $, such that
$\Lambda \cap \SX \neq \emptyset $, one has that $F^\theta _\Lambda $ is either empty or infinite, \iItem every nonempty constructible set is infinite,
\iItem $\{E^\theta _a: a\in \Sigma \}$ is a cover for $\ehull (\SX ) $, \iItem $\ehull (\SX ) $ admits no ground ultra-characters \first {17.8}, \iItem
for every ultra-character $\varphi $ on $\ehull (\SX ) $, there exists an infinite admissible word $\omega $, such that $\varphi =\varphi _{\sigma _\omega }$.

\Proof (i)$\Imply $(ii).\enspace Given any nonempty constructible set, write it as $uF^\theta _\Lambda $, where $\Lambda \subseteq \tSX $ is finite,
$\Lambda \cap \SX \neq \emptyset $, and $u\in \Lambda $, by \first {\FPFormOfHull }.  The conclusion then follows easily by observing that $uF^\theta
_\Lambda $ has the same number of elements as $F^\theta _\Lambda $.

\itmproof (ii)$\Imply $(iii). Notice that, for every $a$ in $\Sigma $, one has that $E^\theta _a$ consists of all words in $\LX $ beginning with the
letter ``$a$'', and having length at least $2$.  Therefore the words of length $1$ do not belong to any $E^\theta _a$.  More precisely $$ \LX \setminus
\bigcup _{a\in \Sigma } E^\theta _a \subseteq \Sigma .  $$

If $Y$ is any nonempty constructible subset of $\SX '=\LX $, then $Y$ is infinite by (ii) so it cannot be contained in the finite alphabet $\Sigma $,
and hence $Y$ must have a nonempty intersection with some $E^\theta _a$, proving (iii).

\itmproof (iii)$\Imply $(iv).  If $\varphi $ is any ultra-character on $\ehull (\SX ) $, then $\varphi $ is tight by \cite [12.7]{actions} and hence $$
\bigvee _{a\in \Sigma } \varphi (E^\theta _a)=1, $$ so there exists $a$ in $\Sigma $ such that $\varphi (E^\theta _a)=1$, whence $a\in \sigma _\varphi $,
proving $\varphi $ not to be a ground character.

\itmproof (iv)$\Imply $(i).  Suppose by way contradiction, that there exists a finite subset $\Lambda \subseteq \tSX $, such that $\Lambda \cap \SX \neq
\emptyset $, and $F^\theta _\Lambda $ is finite and nonempty.  What matters to us is that $F^\theta _\Lambda $ is an example of a finite constructible
set, so there exists a smallest nonempty constructible set $Y\subseteq F^\theta _\Lambda $, which is therefore a minimal element of the {\sla } $\ehull
(\SX ) $.  The principal filter generated by $Y$ is then an ultra-filter, and the associated character, say $\varphi $, is an ultra-character.

In order to compute $\sigma _\varphi $, recall that a given $t$ in $\SX $ belongs to $\sigma _\varphi $ if and only if $$ 1=\varphi (E^\theta _t) =
\bool {Y\subseteq E^\theta _t}, $$ which is to say that $Y\subseteq E^\theta _t$.  Supposing that indeed $Y\subseteq E^\theta _t$, notice that $$ \ell
(y)\geq \ell (t)+1, \for y\in Y, $$ so if we fix any $y_0$ in $Y$, we deduce that $$ \ell (y_0)-1\ \geq \ \max _{t\in \sigma _\varphi }\ \ell (t), $$
so $\sigma _\varphi $ is bounded in length.  A moment's reflection is enough to allow one to see that $\sigma _\varphi $ cannot be an open string,
whence $\varphi $ is not an open character.  We then deduce from \first {\FPNonOpenUltra } that $\varphi =\dualrep _u(\psi )$, where $u\in \tSX $,
and $\psi $ is a ground ultra-character, but these are banned by (iv).  This is therefore a contradiction, and hence (i) is proved.

\itmproof (iv)$\Imply $(v). Given an ultra-character $\varphi $ we have by the hypothesis and \first {\FPNonOpenUltra } that $\varphi =\varphi _\sigma
$, for some open quasi-maximal string $\sigma $.  Using \ref {DescribeStringsAsWords.iv} we have that $\sigma =\sigma _\omega $, for some admissible
word $\omega $, and since $\sigma $ is open, $\omega $ must be infinite by \ref {DescribeStringsAsWords.ii}.

\itmproof (v)$\Imply $(iv).  Arguing by contradiction let $\varphi $ be a ground ultra-character.  By hypothesis we may write $\varphi =\varphi _\sigma
$, where $\sigma =\sigma _\omega $, and $\omega $ is an infinite admissible word.  Using \ref {DescribeStringsAsWords.ii} we have that $\sigma _\omega
$ is an open string, so it follows from \first {\FPSigmaPhiOnOpen .iii} that $$ \sigma _\omega = \sigma = \Sigma \big (\Phi (\sigma )\big ) = \Sigma
(\varphi ) = \sigma _\varphi = \emptyset , $$ a contradiction, thus proving (iv).  \endProof

Let us conclude this section with an example to show that nonempty finite constructible sets may indeed exist and hence the equivalent conditions of
\ref {GroundForSubshifts} do not always hold.  Consider the alphabet $\Sigma =\{a, b, c\}$ and let $\X $ be the subshift on $\Sigma $ consisting of
all infinite words $\omega $ such that, in any block of $\omega $ of length three, there are no repeated letters.  Alternatively, a set of forbidden
words defining $\X $ is the set of all words of length three with some repetition.

It is then easy to see that the language $\LX $ of $\X $ is formed by all finite words on $\Sigma $ with the same restriction on blocks of length three
described above.

Notice that $c\in F_{\{a, b\}}$, because both $ac$ and $bc$ are in $\LX $.  However there is no element in $F_{\{a, b\}}$ other than $c$, because it is
evident that neither $a$ nor $b$ lie in $F_{\{a,b\}}$, and for any $x$ in $\Sigma $, either $acx$ or $bcx$ will involve a repetition. So \emph {voil\`a}
the finite constructible set: $$ F_{\{a,b\}}=\{c\}.  \equationmark FiniteConstructibleSet $$

\section A non essentially tight {\sla } of constructible sets

\label CounterExEssTight In this section we shall provide an example of a subshift whose associated semigroup $\SX $ contains constructible subsets with
infinite boundary.  We will also see that $\ehull (\SX ) $ is not an essentially tight sub-{\sla } of $\P (\SX ')$.  This will show that condition \ref
{GeneralEssTight.iii} cannot be removed from that result.

\def \opt #1#2#3#4{[#1, #2, #3, #4]} \def \any {\kern 1pt \star }

From now on we shall adopt a convention, loosely based on the notion of \emph {regular expressions}\fn {\tt https://en.wikipedia.org/wiki/Regular\_expression},
designed to simplify the specification of sets of words on a given alphabet $\Sigma $.  Given $a$ in $\Sigma $, and $n\geq 0$, we will write $$ a^n $$
to refer to the word $aa\ldots a$, where the letter $a$ repeats $n$ times.  In case $n=0$, then $a^n$ will stand for the empty word.  We will also write
$$ a^* $$ for any word of the form $a^n$, where $n\geq 0$, and in case we wish to refer to $a^n$ only for strictly positive $n$'s, we will write $$ a^+ .  $$

Given any finite collection of members of $\Sigma $, say $a_1,a_2,\ldots ,a_n$, we will write $$ \opt {a_1}{a_2}{\ldots }{a_n} $$ to refer to any word of
length 1, formed by any one of the given $a_i$.  For example, choosing the alphabet $\Sigma =\{0, 1,2, 3, 4\}$, when we refer to a word of the form $$
10^+ 4\opt 0234, $$ we mean any member of the set $$ \{1\underbrace {00\ldots 0}_{n\hbox { \eightrm times}}4a: n\geq 1, \ a\in \{0, 2, 3, 4\}\big \}.
$$ Furthermore, if $\mu $ is any finite word, then by $$ \mu \any \equationmark WordsAny $$ we shall mean any finite word of the form $\mu \nu $,
where $\nu $ is any word, including the empty one.

For example, given $a$ in $\Sigma $, the words of the form $a\any $ are precisely the words beginning in the letter $a$, including the single letter
word $a$ itself.

\bigskip From now on we will fix the alphabet $\Sigma =\{0, 1,2, 3, 4\}$, and we will consider the subshift $\X $ given by specifying the following
set of forbidden words: $$ 10^+ 4\opt 0234, \qquad 20^+ 4\opt 0134, \and 30^+ 4.  $$

Notice that, for any $n\geq 1$, the infinite words $$ 10^n411111\ldots , \and 20^n422222\ldots $$ lie in $\X $, so, in particular, the word $0^n4$
belongs to the language $\LX $ of $\X $, and in fact $$ 0^n4\in F_{\{1, 2\}}, $$ although $0^n4\notin F_3$, because $30^n4$ is forbidden.

Still asuming that $n\geq 1$, notice that if $a$ is any member of $\Sigma $, then either $$ 10^n4a,\quad \hbox {or} \quad 20^n4a $$ is a forbidden
word, according to whether $a$ is of the form $\opt 0234$ or $\opt 0134$.  As a consequence we deduce that $0^n4a$ does not belong to $F_{\{1, 2\}}$,
for any $a$ in $\Sigma $, whence $0^n4$ lies in the boundary $\partial F_{\{1, 2\}}$.  We thus see that $F_{\{1, 2\}}$ has an infinite boundary.

It is our intention to prove that, if $\SX $ is the semigroup associated to $\X $, as in section \ref {SubshiftSgrSection}, then $\ehull (\SX ) $ is not
an essentially tight sub-{\sla } of $\P (\SX ')$.  This will be accomplished by proving that $$ \big \{E_1, \ E_2, \ E_3, \ E_4, \ 0F_{\{10,20,30\}}\big
\} \equationmark HereIsCover $$ is a cover for $F_{\{1,2\}} $, with infinite defect set $$ D:= F_{\{1,2\}} \setminus \big (E_1 \cup E_2 \cup E_3\cup
E_4 \cup 0F_{\{10,20,30\}}\big ).  $$

By far the easiest thing to check is that $D$ is infinite: simply notice that $0^n4$ belongs to $D$, for every $n\geq 1$.  However proving that \ref
{HereIsCover} is a cover for $F_{\{1,2\}} $ is not so straightforward.  In order to do this, let us first describe all elements of $D$.

It is elementary to see that $\Sigma $ is contained in $D$, among other reasons because for every $a$ in $\Sigma $, the words in $E_a$ have length at
least $2$.  In fact we claim that $$ D=\Sigma \cup \{0^n4: n\geq 1\}.  \equationmark HereIsDefectSet $$

Having already checked the inclusion "$\supseteq $", we now focus on proving "$\subseteq $".  For this let us pick any word $\mu $ in $D$ of length $2$
or more.  Since $\mu $ is not in $E_1 \cup E_2 \cup E_3\cup E_4$ we see that the $\mu $ must begin with $0$, so we may write $\mu =0\nu $, where $\nu $ is
a word in $\LX $ of length 1 or more.  Moreover, since $\mu \in F_{\{1,2\}} $, we have that $\nu $ lies in $F_{\{10,20\}}$, so necessarily $\nu \notin F_{30}$.

The word $30\nu $ must therefore involve a forbidden word $\beta $, which is surely not a subword of $\nu $, hence $\beta $ must begin in either
$3$ or $0$.  Given the lack of forbidden words beginning in $0$, we see that $\beta $ must begin in $3$, whence $\beta =30^n4$, for some $n\geq 1$.
From this it follows that $0^{n-1}4$ is a prefix of $\nu $, and hence $0^n4$ is a prefix of $\mu $.  Having already seen that $0^n4$ is a maximal word
in $F_{\{1,2\}} $, we deduce that $\mu =0^n4$.  This shows that any word of length 2 or more in $F_{\{1,2\}} $ lies in the set $\{0^n4: n\geq 1\}$,
and hence \ref {HereIsDefectSet} is proved.

Showing our stated goal that \ref {HereIsCover} is a cover for $F_{\{1,2\}} $ consists in showing that every nonempty constructible subset of $F_{\{1,2\}} $
must intercept one of the members of \ref {HereIsCover}, which is the same as showing that $D$ contains no nonempty constructible set.  We must therefore
get a good handle on constructible sets.

Given $\mu $ in $\LX $, let us first analyze $F_\mu $.  If $\nu $ is a member of $\LX $ such that $\nu \notin F_\mu $, then $\mu \cdot \nu =0$,
which means that $\mu \nu $ contains a forbidden subword, say $\beta $.  Evidently $\beta $ is not a subword of either $\mu $ or $\nu $, so $\beta
$ must straddle $\mu $ and $\nu $, meaning that there is a suffix $\mu '$ of $\mu $ and a prefix $\nu '$ of $\nu $, such that $\mu '\nu '=\beta $.
Consequently, if no suffix of $\mu $ combine with any prefix of $\nu $ to form a forbidden word, we must have that $\nu \in F_\mu $.

We summarize this as follows:

\state Proposition Given $\mu $ and $\nu $ in $\LX $, one has that $\nu \notin F_\mu $, if and only if there is a forbidden word of the form $\mu '\nu
'$, where $\mu '$ is a suffix of $\mu $, and $\nu '$ is a prefix of $\nu $.

Suppose, for example that $10^n4$ is a suffix of $\mu $, where $n\geq 1$.  Then, for any $\nu $ beginning in $0$, $2$, $3$ or $4$, the concatenated word
$\mu \nu $ will involve a forbidden word of the form $10^+ 4\opt 0234$.  On the other hand, if $\nu $ begins in $1$, then $\mu \nu $ will not contain
any forbidden word because $10^n41$ cannot be extended towards the left or right-hand-side in such a way as to form a forbidden word involving the \emph
{bridge} ``$41$''.  As a consequence we deduce that $ F_\mu $ consists precisely of the words in $\LX $ of the form $1\any $ (as defined in \ref {WordsAny}).

Arguing as above we may determine every $F_\mu $, as follows:

\state Proposition Let $\mu \in \LX $.  \iaitem \aitem If $\mu $ ends in $10^+ 4$, then $F_\mu $ consists of all admissible words of the form $1\any $.
\aitem If $\mu $ ends in $20^+ 4$, then $F_\mu $ consists of all admissible words of the form $2\any $.  \aitem If $\mu $ ends in $4$, and it does
not fit either (a) or (b) above, then $F_\mu =\LX $.  \aitem If $\mu $ ends in $10^+ $, then $F_\mu $ consists of all admissible words, except those
of the form $0^* 4\opt 0234\any $.  \aitem If $\mu $ ends in $20^+ $, then $F_\mu $ consists of all admissible words, except those of the form $0^*
4\opt 0134\any $.  \aitem If $\mu $ ends in $30^+ $, then $F_\mu $ consists of all admissible words, except those of the form $0^* 4\any $.  \aitem If
$\mu $ ends in $0$, and it does not fit either (d), (e) or (f) above, then $F_\mu =\LX $.  \aitem If $\mu $ ends in $1$, then $F_\mu $ consists of all
admissible words, except those of the form $0^+ 4\opt 0234\any $.  \aitem If $\mu $ ends in $2$, then $F_\mu $ consists of all admissible words, except
those of the form $0^+ 4\opt 0134\any $.  \aitem If $\mu $ ends in $3$, then $F_\mu $ consists of all admissible words, except those of the form $0^+ 4\any $.

All possible alternatives for $\mu $ being taken care of, we have described all possible sets $F_\mu $.  But, despite this detailed description, the
only relevant information we care, easily deducible from the above, is as follows:

\state Proposition \label DescribeFMu Let $\mu \in \LX $.  Then one and only one of the following properties hold: \iItemize \iItem $F_\mu $ consists
of all admissible words of the form $1\any $.  \iItem $F_\mu $ consists of all admissible words of the form $2\any $.  \iItem $F_\mu $ contains all
words mentioned in (i) and (ii), above.

Once we know this much we may also say something relevant about sets of the form $F_\Lambda $:

\state Proposition \label DescribeFLambda Let $\Lambda $ be any subset of $\tS $ and suppose that $F_\Lambda $ is nonempty.  Then, substituting $F_\Lambda
$ for $F_\mu $, one has that exactly one among the properties \ref {DescribeFMu.i--iii} hold.

\Proof This follows immediately from \ref {DescribeFMu}, since $F_\Lambda $ is the intersection of finitely many $F_\mu $'s.  \endProof

Given any nonempty constructible subset $X$ of $\SX '$, we then have by \first {\FPFormOfHull } that $X=uF_\Lambda $, where $\Lambda \subseteq \tS $
is a finite subset, $\Lambda \cap \SX \neq \emptyset $, and $u\in \Lambda $.  Combining this with \ref {DescribeFLambda} we then deduce that either $X$
contains all admissible words of the form $\mu 1\any $, or else $X$ contains all admissible words of the form $\mu 2\any $.

Regarding the defect set $D$ mentioned in \ref {HereIsDefectSet}, it is therefore evident that no non\-empty constructible subset of $\SX '$ may be
found inside $D$, whence \ref {HereIsCover} is indeed a cover for $F_{\{1,2\}} $.

Since we have already determined that $D$ is infinite we have completed the verification that $\ehull (\SX ) $ is not an essentially tight sub-{\sla }
of $\P (\SX ')$.

The big conclusion, which motivated the present section, is that Theorem \ref {GeneralEssTight} is rendered false if its third condition is removed.

\section Subsets of the spectrum of the {\sla } of constructible sets

In this section we will briefly return to the general situation of a $0$-left-cancellative semigroup $S$ in order to point out certain noteworthy
subsets of $\Ehat (S)$.

Once $\Ehat (S)$ is seen as the unit space of the universal groupoid of $\hull (S)$, one may reduce the latter to these subsets obtaining various related
groupoids.  When applied to the semigroup $\SX $ for a given subshift $\X $, and upon considering the C*-algebra of such subgroupoids, we will have
produced most of the C*-algebras that have been studied in connection to the subshift $\X $, such as the Matsumoto and the Carlsen-Matsumoto C*-algebras.
By studying conditions under which these subsets agree one can obtain conditions for the different versions of Matsumoto's algebras to agree.

\fix Throughout this section we therefore fix a $0$-left-cancellative semigroup $S$ admitting least common multiples.  Among the well known subsets of
$\Ehat (S)$ one has \def \myalign #1=#2;{\hbox to 1.5cm{\hfill $#1$} = \hbox to 6cm{$#2$\hfill }} $$ \myalign \Ehat _\infty (S)= \{\varphi \in \Ehat (S):
\varphi \hbox { is an ultra-character}\},; $$ and $$ \myalign \Ehat \sub {tight}(S) = \{\varphi \in \Ehat (S): \varphi \hbox { is a tight character}\},;
$$ to which we would like to add a few more, beginning with $$ \myalign \Ehat \sub {max}(S) = \{\varphi _\sigma : \sigma \hbox { is a maximal string}\}.; $$

Since $\ehull (S) $ is a sub-{\sla } of $\P (S')$, we may view the inclusion map $$ \iota : \ehull (S) \to \P (S'), \equationmark InclusionRep $$
as a representation and we may then consider the essentially tight characters (see \ref {IntroEssTightChar}) relative to $\iota $.

\definition \label DefSess We will denote by $\Ehat \sub {ess}(S) $ the set of all characters on $\ehull (S) $ which are essentially tight relative to
the representation $\iota $ described in \ref {InclusionRep}.

In other words, a character $\varphi $ is in $\Ehat \sub {ess}(S) $ if and only if, whenever $X, Y_1,\ldots ,Y_n$ are in $\ehull (S) $, and the symmetric
difference $$ X\mathop {\Delta }\big (\bigcup _{i=1}^nY_i\big ) $$ is finite, then $\varphi (X)=\bigvee _{i=1}^n\varphi (Y_i).$ In particular, given
any finite set $X$, we necessarily have that $$ \varphi (X)=0, \equationmark SessKillFinite $$ since one could take all of the $Y_i$ above to be empty.

\state Proposition \label SessClosed $\Ehat \sub {ess}(S) $ is a closed subset of\/ $\Ehat (S)$.

\Proof This follows from \first {\FPRelTightClosed } by observing that $\Ehat \sub {ess}(S) $ coincides with the set of all $\pi $-tight characters,
where $\pi $ is the representation of $\ehull (S) $ obtained as the composition $$ \ehull (S) \arw \iota \P (S')\arw q\Q (S') $$ where $q$ is the
quotient mapping mentioned in \ref {IntroQuotientBoole}. \endProof

It is our next major goal to understand the relationship between all of the above subsets of $\Ehat (S)$, which we will henceforth denote by $$
\sinfShort , \ \stightShort ,\ \smaxShort \and \sessShort , $$ for simplicity.

It is well known that $\sinfShort $ is a dense subset of $\stightShort $.  In case all maximal strings are open, such as for semigroups possessing right
local units, or for semigroups associated to subshifts (see \ref {DescribeStringsAsWords}), we have by \first {\FPBigNewTightResult } that $\varphi
_\sigma $ is an ultra-character for every maximal string $\sigma $, whence $ \smaxShort \subseteq \sinfShort , $ and then the first three of the above
sets are related by $$ \smaxShort \subseteq \sinfShort \subseteq \stightShort .  \equationmark FirstThreeRel $$

We next want to study a situation in which $\sessShort $ may be included in the above picture.

\state Lemma \label FiniteAlternative Assuming that $S$ possesses a length function $\ell $, let $\sigma $ be a string in $S^\star $ and let $X\in
\ehull (S) $.  Then \iItemize \iItem $\sigma \in \varepsilon (X) \Imply \sigma \setminus X$ is bounded (relative to $\ell $, as defined in \first
{\FPBoundedSets }), \iItem $\sigma \notin \varepsilon (X) \Imply \sigma \cap X$ is bounded.  \medskip \noindent Moreover, if $\sigma $ is unbounded,
then the converse of (i) and (ii) also hold.

\Proof We begin by writing $X=\theta _u(F_\Lambda )$, where $\Lambda $ is a finite subset of $\tS $, intersecting $S$ non-trivially, and $u\in \Lambda $.
In order to prove (i), assume that $$ \sigma \in \varepsilon (X)\explica {\first {\FPIntroEpsilon }}=\theta ^\star _u(F^\star _\Lambda ), $$ so $$
\emptyset \neq \sigma \cap E_u \subseteq \theta _u(F_\Lambda ) = X, \equationmark FRFstarRTwoRepeated $$ by \first {\FPFRFstarRTwo }.  Let us now split
the argument according to whether $u=1$ or $u\in S$.  Under the former assumption, $E_u=S'$, so $\sigma \subseteq X$, by \ref {FRFstarRTwoRepeated},
whence $\sigma \setminus X$ is actually empty.

If $u\in S$, then the fact that $\sigma \cap E_u$ is nonempty implies that there exists some $t$ in $S$ such that $ut\in \sigma $.  We then claim that $$
\{s\in \sigma :\ell (s)>\ell (u)\} \subseteq E_u.  \equationmark LongHasPrefix $$

In order to prove it, pick any $s$ in $\sigma $ with $\ell (s)>\ell (u)$, and use \first {\FPStringDirected } to obtain $v$ and $w$ in $\tS $ such that
$sv=utw\in \sigma $.  It follows that $u\|sv$, and by length considerations we conclude that $u\|s$, so we may write $s=uz$, for some $z$ in $\tS $.
Clearly $z\neq 1$, since otherwise $\ell (s)=\ell (u)$, so $z$ lies in $S$, and hence $s\in E_u$, proving the claim.

We then have that $$ \sigma \setminus X \explain {FRFstarRTwoRepeated}\subseteq \sigma \setminus E_u \explain {LongHasPrefix}\subseteq \{s\in \sigma
:\ell (s)\leq \ell (u)\}, $$ so $\sigma {\setminus } X$ is bounded, as desired.  This proves (i).

In order to prove (ii), suppose that $\sigma \notin \theta ^\star _u(F^\star _\Lambda )$, so that \ref {FRFstarRTwoRepeated} fails by \first
{\FPFRFstarRTwo }.  Should the failure be due to the fact that $\sigma \cap E_u$ is empty, then $$ \sigma \subseteq S'\setminus E_u \subseteq S'\setminus
X, $$ from where we see that $\sigma \cap X$ is actually empty.  On the other hand, suppose that $\sigma \cap E_u$ is nonempty, so the failure of \ref
{FRFstarRTwoRepeated} must be due to the fact that $\sigma \cap E_u$ is not a subset of $\theta _u(F_\Lambda )$.  Picking any $$ s \in (\sigma \cap
E_u) \setminus \theta _u(F_\Lambda ), $$ write $s=ur$, for some $r$ in $S$, and observe that $r$ is not in $F_\Lambda $, so there exists some $t$ in
$\Lambda $ such that $tr=0$.  We next claim that, for every $p$ in $\sigma $, one has that $$ \ell (p)>\ell (s) \Imply p\notin \theta _u(F_\Lambda ).
\equationmark BigNonThere $$ In order to prove this claim, choose any $p$ in $\sigma $ with $\ell (p)>\ell (s)$, and observe that since both $s$ and $p$
lie in $\sigma $, by \first {\FPStringDirected } there are $v$ and $w$ in $\tS $, such that $sv=pw\in \sigma $.  We then have that $s\|p$, by comparing
lengths, so $p=sz$, for some $z$ in $\tS $, whence $$ p=sz=urz.  $$ From this it follows that $p$ lies in $E_u$, so to show that $p$ is not in $\theta
_u(F_\Lambda )$, it is enough to check that $\theta _u\inv (p)$ is not in $F_\Lambda $, but this follows easily from $$ t\kern 1pt\theta _u\inv (p)=trz=0.  $$

Claim \ref {BigNonThere} is therefore verified, so we deduce that $$ \sigma \cap X\subseteq \{p\in \sigma :\ell (p)\leq \ell (s)\}, $$ thus showing
that $\sigma \cap X$ is bounded.

Assuming that $\sigma $ is unbounded, let us now prove the converse of (i), so suppose that $\sigma \setminus X$ is bounded.  Arguing by contradiction,
suppose that $\sigma $ is not in $\varepsilon (X)$, so (ii) implies that $\sigma \cap X$ is bounded, and since $$ \sigma = (\sigma \cap X)\cup (\sigma
\setminus X), $$ we deduce that $\sigma $ is bounded, a contradiction.  This completes the proof of the converse of (i), while the converse of (ii)
may be proved by following a similar pattern.  \endProof

It is interesting to notice that, as a consequence of the above result, a string cannot be split into an unbounded part inside, and an unbounded part
outside a given $\theta $-constructible set.

In case $\sigma $ is a bounded string then both \ref {FiniteAlternative.i-ii} are obviously true, so the relevance of the above result is really
restricted to unbounded strings.

When $\sigma $ is unbounded, it is necessarily non-degenerate, so, according to \first {\FPIntroPhiSigmaNonDeg }, one has that $\varphi _\sigma $ is a
character on $\ehull (S) $.  Recalling that $$ \varphi _\sigma (X)=\bool {\sigma \in \varepsilon (X)}, $$ we may rephrase \ref {FiniteAlternative} as follows:

\state Corollary \label FiniteAlternativeChar Assuming that $S$ possesses a length function $\ell $, let $\sigma $ be an unbounded string in $S^\star $
and let $X\in \ehull (S) $. Then \iItemize \iItem $\varphi _\sigma (X)=1 \Iff \sigma \setminus X$ is bounded, \iItem $\varphi _\sigma (X)=0 \Iff \sigma
\cap X$ is bounded.

Recall from \ref {DefLocFin} that a length function is said to be locally finite if $$ \{s\in S': \ell (s)\leq n\}, $$ is a finite set for every $n$.
In this case every bounded set is finite and it is obvious that, conversely, every finite set is bounded.  In the conclusions of the last two results
above, one could therefore replace each occurrence of the word ``bounded" by the word ``finite".

\state Proposition \label UnboundedIsEssential Assuming that $S$ possesses a locally finite length function $\ell $, one has that $\varphi _\sigma \in
\pilar {11pt}\sessShort $, for every unbounded string $\sigma $.

\Proof Suppose that $Y_1,\ldots ,Y_n,X$ are in $\ehull (S) $ and $$ q(X)=\bigvee _{i=1}^nq(Y_i), $$ where $q$ is the quotient map introduced in \ref
{IntroQuotientBoole}.  This means that $X$ coincides with $\bigcup _{i=1}^nY_i$, up to a finite set, in the sense that the symmetric difference is finite.
We must then prove that $$ \varphi _\sigma (X)=\bigvee _{i=1}^n\varphi _\sigma (Y_i).  \equationmark PhiSigmaSup $$ For this, let us first assume that
$\varphi _\sigma (X)=1$, which \ref {FiniteAlternativeChar} says is equivalent to the fact that $$ \sigma \setminus X \equationmark FirstFinite $$
is bounded, hence finite.  Incidentally this means that $\sigma $ is mostly contained in $X$ in the sense that the part of $\sigma $ which is not
contained in $X$ is finite.

If we assume that, contrary to what is required, $\varphi _\sigma (Y_i)=0$, for all $i$, then, again by \ref {FiniteAlternativeChar} we would have that
$\sigma \cap Y_i$ is bounded hence finite, so $$ \sigma \cap \big (\bigcup _{i=1}^nY_i\big ) \equationmark SecondFinite $$ is also finite.  Removing from
$\sigma $ the finite sets \ref {FirstFinite} and \ref {SecondFinite}, we would then be left with an infinite set contained in $$ X\setminus \bigcup
_{i=1}^nY_i, $$ hence contradicting the fact that the symmetric difference between $X$ and $\bigcup _{i=1}^nY_i$ is finite.  It then follows that
$\varphi _\sigma (Y_i)=1$, for some $i$, hence proving \ref {PhiSigmaSup}.

On the other hand, suppose that $\varphi _\sigma (X)=0$, and, again arguing by contradiction, that $\varphi _\sigma (Y_i)=1$, for some $i$.  Then by the
local finiteness of $\ell $, and by \ref {FiniteAlternativeChar}, we have that $\sigma \cap X$ and $\sigma \setminus Y_i$ are finite sets, and if we remove
these from $\sigma $ we will be left with an infinite set contained in $Y_i\setminus X$, and hence also in $$ \bigcup _{i=1}^nY_i\setminus X, $$ once more
contradicting the finiteness of the symmetric difference between $X$ and $\bigcup _{i=1}^nY_i$.  This concludes the proof of \ref {PhiSigmaSup}.  \endProof

In \first {\FPBackInvarCorol } we have already seen the relevance of the hypothesis that maximal strings be unbounded.  Another such situation is the
following immediate consequence of \ref {UnboundedIsEssential}:

\state Theorem \label SmaxInSess Let $S$ be a $0$-left-cancellative semigroups admitting least common multiples, and possessing a locally finite length
function $\ell $.  Assuming that every maximal string is unbounded, one has that $\pilar {11pt}\smaxShort \subseteq \sessShort $.

If a string $\sigma $ is not open, we have seen in \first {\FPPropsOpenString .ii} that necessarily $\sigma =\delta _r$, for some $r$ in $S$.  In this
case it is obvious that $\ell (s)\leq \ell (r)$, for each $s$ in $\sigma $, and in particular we see that $\sigma $ is bounded.  This says that every
unbounded string is open, so, under the hypotheses of the above result (every maximal string is unbounded), we have that every maximal string is open,
in which case \ref {FirstThreeRel} holds.  This proves the following:

\state Corollary \label InclusionsSpec Under the assumptions of \ref {SmaxInSess} one has that $$ \def \quad {\kern 8pt} \matrix { \smaxShort & \subseteq
& \sinfShort & \subseteq & \stightShort \cr \kern 3pt \rotatebox {-90}{$\subseteq $}\hfill \cr \pilar {16pt} \sessShort } $$

We should remark that, for every subshift $\X $, the associated semigroup $\SX $ satisfies the assumptions of \ref {SmaxInSess} by \ref
{DescribeStringsAsWords.iii}.

Let us now present a few examples to illustrate that certain inclusions mentioned above may be proper.  In all of the examples below we will consider
a specific subshift $\X $ and we will always refer to the semigroup $\SX $ introduced in \ref {IntroSX}.

\state Proposition For a suitably chosen subshift $\X $, there exists an ultra-character $\psi _1$ not belonging to $\sessShort $.  Consequently $\psi
_1$ is not in $\smaxShort $ either, and hence $\smaxShort \varsubsetneq \sinfShort $.

\Proof Considering the subshift defined at the end of section \ref {SubshiftSgrSection}, recall from \ref {FiniteConstructibleSet} that $F_{\{a,b\}}=\{c\}$.
We then have that $\{c\}$ is a minimal element of $\ehull (S) $, whence the character $\psi _1$ defined by $$ \psi _1(X) = \bool {\{c\}\subseteq X}
$$ is an ultra-character.  Since $ \psi _1\big (\{c\}\big ) = 1, $ we see that $\psi _1$ assigns a nonzero value to a finite set, whence $\psi _1$
is not in $\sessShort $ by \ref {SessKillFinite}.  \endProof

As seen in \ref {SessClosed}, we have that $\sessShort $ is closed, hence the closure of $\smaxShort $ is contained in $\sessShort $.  Therefore $\psi
_1$ is not in the closure of $\smaxShort $, whence $\smaxShort $ is not even dense in $\sinfShort $.

It is interesting to observe that under the equivalent conditions of \ref {GroundForSubshifts}, one has that $\smaxShort $ coincides with $\sinfShort $.

In the diagram below we picture Venn diagrams for the four subsets of $\Ehat (S)$ under analysis, highlighting the inclusions already mentioned in
\ref {InclusionsSpec}.  Sets we know are always closed are drawn with a solid line.  We also illustrate the character $\psi _1$ mentioned in the above result.

\bigskip \bigskip \bigskip \null \hfill 
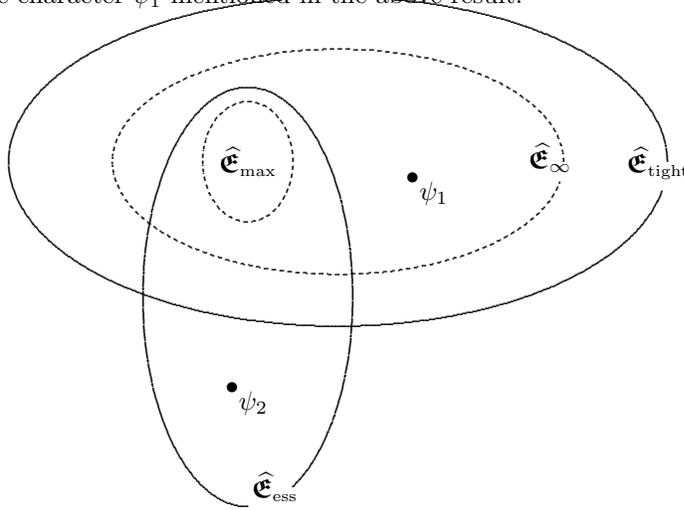
\beginpicture \setcoordinatesystem units <0.0060truecm, -0.0090truecm> point at -1300 -3000 \ellipticalarc
axes ratio 2:1 350 degrees from 730 0 center at 00 0 \put {$\stightShort $} at 710 2 \ellipticalarc axes ratio 1:2 -335 degrees from -200 510 center
at -200 200 \put {$\sessShort $} at -140 480 \put {$\bullet _{ _{\textstyle \psi _1}}$} at 200 40 \put {$\bullet _{ _{\textstyle \psi _2}}$} at -200
350 \setdashes <1.5pt> \ellipticalarc axes ratio 2:1 350 degrees from 500 0 center at 0 0 \put {$\sinfShort $} at 470 -6 \ellipticalarc axes ratio 3:4
360 degrees from -300 0 center at -200 0 \put {$\smaxShort $} at -200 00 \endpicture \hfill \null \bigskip \bigskip

As already pointed out, $\sinfShort $ is dense in $\stightShort $, but unfortunately Venn diagrams cannot depict fine topological features such as density.

\state Proposition For a suitably chosen subshift, there exists a character $\psi _2$ in $\sessShort $ which is not tight.

\Proof Consider the subshift described in section \ref {CounterExEssTight}.  There we showed that $\ehull (S) $ is not an essentially tight sub-{\sla }
of $\P (\LX )$. In other words, the inclusion representation $\iota $ of \ref {InclusionRep} is not essentially tight.

We next plan to use \ref {NecSufEssTight}, so we must adapt ourselves to the hypotheses required there and we will do so by verifying \ref {NecSufEssTight.a},
namely that $\LX $ can be written, up to a finite set, as the union of finitely many $\theta $-constructible sets.  This is in fact easily checked,
since the only elements of $\LX $ not in $$ \bigcup _{a\in \Sigma }E_a $$ are the words of length one of which there are finitely many.

The application of \ref {NecSufEssTight} is thus legitimized, and hence the fact that $\iota $ is not essentially tight implies that \ref {NecSufEssTight.ii}
fails, meaning that there exists a character $\psi _2$ in $\sessShort $ which is not tight, as desired.  \endProof

The character $\psi _2$ of the above result therefore belongs to $\sessShort $ and not to $\smaxShort $, so we see that $$ \smaxShort \varsubsetneq
\sessShort .  $$ Since $\stightShort $ is closed and contains $\smaxShort $, it follows that $\psi _2$ is not in the closure of $\smaxShort $.
This shows that $\smaxShort $ is not even dense in $\sessShort $.

\section Carlsen-Matsumoto Condition

The main goal of this section is to study a situation in which, contrary to the last example of the previous section, $\smaxShort $ is dense in $\sessShort $.

For the time being we will let $S$ be a fixed $0$-left-cancellative semigroup admitting least common multiples.  Recall from \first {7.7} that for each
nonempty finite subset $\Lambda $ of $\tS $, we have $$ F^\theta _\Lambda = \{s\in S': ts\neq 0,\ \forall t\in \Lambda \}.  $$

Here we would like to extend the above notion as follows:

\definition \label DefineDDom Given finite subsets $\Lambda $ and $\Gamma $ of $\tS $, we will denote by $F^\theta _{\Lambda , \Gamma }$ the subset of
$S'$ given by $$ F^\theta _{\Lambda ,\Gamma } = \{s\in S': ts\neq 0,\ \forall t\in \Lambda , \ rs=0,\ \forall r\in \Gamma \}.  $$

An alternative way to express the above is clearly $$ F^\theta _{\Lambda ,\Gamma } = \big (\medcap _{t\in \Lambda } F^\theta _t\big ) \cap \big (\medcap
_{r\in \Gamma } S'\setminus F^\theta _r\big ).  \equationmark AltDescrDDom $$

There is still another way to write $F^\theta _{\Lambda ,\Gamma }$, inspired by covers and their associated defect sets.

\state Proposition \label OtherDescrDDom Let $\Lambda $ and\/ $\Gamma $ be finite subsets of $\tS $.  Then $$ F^\theta _\Lambda \supseteq \medcup _{r\in
\Gamma } F^\theta _{\Lambda \cup \{r\} }, $$ and $$ F^\theta _\Lambda \setminus \big ( \medcup _{r\in \Gamma } F^\theta _{\Lambda \cup \{r\} }\big )
= F^\theta _{\Lambda ,\Gamma }.  $$

\Proof Left for the reader. \endProof

\definition \label CarlMatsuCond We will say that $S$ satisfies Carlsen and Matsumoto's condition $(*)$ if, whenever $\Lambda $ and $\Gamma $ are finite
subsets of $\tS $ such that $F^\theta _{\Lambda ,\Gamma }$ is infinite, there exists a maximal string $\sigma $ such that $$ \sigma \subseteq F^\theta
_t, \and \sigma \not \subseteq F^\theta _r, $$ for all $t$ in $\Lambda $, and for all $r$ in $\Gamma $.

Notice that by \first {\FPFRFstarR .i}, the above condition on $\sigma $ is equivalent to $$ \sigma \in F^\star _t, \and \sigma \notin F^\star _r,
$$ which in turn is the same as saying that $$ \sigma \in \big (\medcap _{t\in \Lambda } F^\star _ t\big ) \cap \big (\medcap _{r\in \Gamma } S^\star
\setminus F^\star _r\big ) =: F^\star _{\Lambda ,\Gamma }, $$ where the above definition of $F^\star _{\Lambda ,\Gamma }$ is of course compatible with
\ref {AltDescrDDom}, provided we think of it as applied to the representation $\theta ^\star $ of $S$ on $S^\star $ given by \first {\FPIntroStarAction .ii}.

Condition $(*)$ may then be phrased in the following very concise way:

\state Proposition \label ConciseCMCond A $0$-left-cancellative semigroup $S$ satisfies condition $(*)$ if and only if, for all finite subsets $\Lambda ,
\Gamma \subseteq \tS $, one has that $$ F^\theta _{\Lambda ,\Gamma } \hbox { is infinite } \Imply F^\star _{\Lambda ,\Gamma }\cap S^\infty \neq \emptyset .  $$

We leave it for the reader to check that, when $S$ is the language semigroup associated to a subshift, then \ref {CarlMatsuCond} is equivalent to
condition $(*)$ introduced by Carlsen and Matsumoto in \cite [Section 3]{MatsuCarl}, but we warn the reader that the description of condition $(*)$
in \cite {MatsuCarl} is incorrectly stated and must be amended by requiring that the sequence $\{\mu _i\}_i$, mentioned there, have an infinite range.

The relevance of condition $(*)$ to our theory is highlighted in the following:

\state Proposition \label StarNecess Let $S$ be a $0$-left-cancellative semigroup admitting least common multiples.  If \ $\smaxShort $ is dense in
$\sessShort $, then condition $(*)$ holds.

\Proof Let $\Lambda $ and $\Gamma $ be finite subsets of $\tS $, such that $F^\theta _{\Lambda ,\Gamma }$ is infinite.  Choose any non principal
ultrafilter $\xi $ on $S'$, such that $F^\theta _{\Lambda ,\Gamma }\in \xi $, and let $\varphi $ be the character of $\P (S')$ associated to $\xi $.
Therefore $\varphi $ preserves meets and joins, $\varphi (C)=0$, for every finite subset $C$ of $S'$, and $\varphi (F^\theta _{\Lambda ,\Gamma })=1$.

Regarding the discussion right after \ref {DefSess}, observe that if $X, Y_1,\ldots ,Y_n$ are in $\ehull (S) $, and the symmetric difference $$ X\mathop
{\Delta }\big (\medcup _{i=1}^nY_i\big ) $$ is finite, then $$ \varphi (X)=\varphi \big (\medcup _{i=1}^nY_i\big ) = \bigvee _{i=1}^n\varphi (Y_i).
$$ This implies that the restriction of $\varphi $ to $\ehull (S) $, which we also denote by $\varphi $ by abuse of language, belongs to $\sessShort $.

Observe that $$ F^\theta _{\Lambda ,\Gamma } \subseteq F^\theta _\Lambda , \and F^\theta _{\Lambda ,\Gamma }\cap F^\theta _{\Lambda \cup \{r\} }=\emptyset
, \for r\in \Gamma .  $$ Therefore the fact that $\varphi (F^\theta _{\Lambda ,\Gamma })=1$ implies that $$ \varphi (F^\theta _\Lambda )=1, \and \varphi
(F^\theta _{\Lambda \cup \{r\} })=0, \for r\in \Gamma .  $$ It follows that $\varphi $ lies in the open subset of $\Ehat (S)$ defined by $$ U=\{\psi
\in \Ehat (S): \psi (F^\theta _\Lambda )= 1, \ \psi (F^\theta _{\Lambda \cup \{r\} })=0, \ \forall r\in \Gamma \}.  $$

The hypothesis that $\smaxShort $ is dense in $\sessShort $ thus leads to the existence of a maximal string $\sigma $ such that $\varphi _\sigma $
lies in $U$.  For every $t$ in $\Lambda $, we have that $F^\theta _\Lambda \subseteq F^\theta _t$, whence $$ 1=\varphi _\sigma (F^\theta _\Lambda )
\leq \varphi _\sigma (F^\theta _t), $$ and hence it follows that $$ 1=\varphi _\sigma (F^\theta _t)=\bool {\sigma \in F^\star _ t}, $$ so $\sigma \in
F^\star _ t$.  On the other hand, for every $r$ in $\Gamma $, we observe that $$ 0=\varphi _\sigma (F^\theta _{\Lambda \cup \{r\} }) = \varphi _\sigma
(F^\theta _{\Lambda }\cap F^\theta _r) = \varphi _\sigma (F^\theta _{\Lambda }) \varphi _\sigma (F^\theta _r) = \varphi _\sigma (F^\theta _r), $$
whence $$ 0=\varphi _\sigma (F^\theta _r)=\bool {\sigma \in F^\star _ r}, $$ so $\sigma \notin F^\star _ r$.  This shows that $$ \sigma \in F^\star
_{\Lambda ,\Gamma }\cap S^\infty , $$ hence proving condition $(*)$ via \ref {ConciseCMCond}.  \endProof

Having shown that condition $(*)$ is necessary for the density of $\smaxShort $ in $\sessShort $, let us now prove that it is also sufficient for
semigroups associated to subshifts.

\state Proposition \label CMConditionForSX Let $\X $ be a subshift, and let $\SX $ be the associated semigroup.  Then the following are equivalent:
\iItemize \iItem $\SX $ satisfies condition $(*)$, \iItem $\sessShort $ coincides with the closure of\/ $\smaxShort $.

\Proof In view of \ref {StarNecess} we need only prove that (i) implies (ii).  As already seen, $\SX $ satisfies the hypotheses of \ref {SmaxInSess},
and hence $$ \smaxShort \subseteq \sessShort .  $$ Since the latter is closed, we see that the closure of the former is contained in the latter, so it
remains to show that $\smaxShort $ is dense in $\sessShort $.

Let $\varphi $ be any character in $\sessShort $ and let $U$ be an open subset of $\Ehat (S)$ containing $\varphi $.  Our task is to find a maximal
string $\sigma $, such that $\varphi _\sigma $ lies in $U$.

By shrinking $U$ a bit if necessary, we may assume that $U$ is a basic open set, meaning that $$ U=\{\psi \in \Ehat (S): \psi (X_i)=1,\ \psi (Y_j)=0,
\ \forall i=1, \ldots , n, \ \forall j=1, \ldots , m\}, $$ where $n,m\geq 0$, and the $X_i$ and the $Y_j$ are suitably chosen elements of $\ehull (S) $.

Since $\varphi $ is nonzero, there exists some $X$ in $\ehull (S) $ such that $\varphi (X)=1$, which in turn may be added to the $X_i$, making $U$
a bit smaller, but still containing $\varphi $, and allowing us to assume that $n\geq 1$.

Notice that for any character $\psi $, one has that $\psi (X_i)=1$, for every $i$, if and only if $\psi (X)=1$, where $$ X=\medcap _{i=1}^nX_i $$ where
we have reset the notation $X$.  So we may replace all of the $X_i$ by the single $X$, in the sense that $$ U=\{\psi \in \Ehat (S): \psi (X)=1,\ \psi
(Y_j)=0, \ \forall j=1, \ldots , m\}.  \equationmark DefnU $$

We will next make a further reduction in order to be able to assume that the $Y_j$ are subsets of $X$, as follows: letting $Z_j=X\cap Y_j$, for every
character $\psi $ one has that $$ \psi (Z_j)=\psi (X)\psi (Y_j), $$ so, assuming that $\psi (X)=1$, we have that $$ \psi (Y_j)=0 \IFF \psi (Z_j)=0.  $$
We may therefore replace the $Y_j$ by the corresponding $Z_j$ in \ref {DefnU} without altering $U$, hence allowing us to proceed under the assumption
that $Y_j\subseteq X$, for every $j$.

By \first {\FPFormOfHull } we may write $X=uF_\Lambda $, for some finite set $\Lambda \subseteq \tS $, with $\Lambda \cap S\neq \emptyset $, and $u\in
\Lambda $.  Since the $Y_j$ are subsets of $X$ we may furthermore use \first {\FPFormOfSubSets } to write $Y_j=ux_jF_{\Delta _j}$, where $x_j\in \tS $,
each $\Delta _j$ is a finite subset of $\tS $, and $\Lambda x_j\subseteq \Delta _j$.

For each letter $a$ in the alphabet $\Sigma $, consider the set\fn {No claim is being made regarding whether or not $ua\neq 0$ (meaning that $ua$ lies
in $\LX $), but when $ua=0$, then $uaF_{\Lambda a}$ is clearly empty.} $$ uaF_{\Lambda a} = \{uax\in S: x\in S, \ tax\neq 0, \ \forall t\in \Lambda \},
$$ and observe that $$ uF_\Lambda \supseteq \medcup _{a\in \Sigma } uaF_{\Lambda a}.  $$

Wondering about the validity of the reverse inclusion, notice that any element $s$ in $uF_\Lambda $, with length at least $\ell (u)+2$, must have the
form $s=ux$, with $\ell (x)\geq 2$, and if the first letter of $x$ is $a$, then necessarily $s\in uaF_{\Lambda a}$.  From this it clearly follows that $$
uF_\Lambda \setminus \medcup _{a\in \Sigma } uaF_{\Lambda a} $$ is a finite set, and since $\varphi $ is essentially tight relative to $\iota $, we deduce
that $$ 1=\varphi (X)=\varphi (uF_\Lambda )= \bigvee _{a\in \Sigma } \varphi (uaF_{\Lambda a}), $$ whence $\varphi (uaF_{\Lambda a})=1$, for some $a$.
Repeating this argument sufficiently many times we may therefore find an arbitrarily long word $x$, such that $$ \varphi (uxF_{\Lambda x})=1, $$ and
$uxF_{\Lambda x}\subseteq X$ (for obvious reasons).  We shall therefore choose an $x$, as above, such that $$ \ell (x) > \ell (x_j), \for j=1, \ldots , m.  $$

We next consider the following sets $$ X'=uxF_{\Lambda x}, $$ $$ Y_j'=X'\cap Y_j, \for j=1,\ldots ,m, $$ and $$ U'=\{\psi \in \Ehat (S): \psi (X')=1,\
\psi (Y'_j)=0, \ \forall j=1, \ldots , m\}.  $$

It is a routine task to show that $$ \varphi \in U'\subseteq U, $$ and if the reader is wondering why is it worth introducing $U'$, the reason is that
it has a nicer description in the sense that the $Y'_j$ are more closely related to $X'$, but we still need a little more work to see that.

Staring at $Y'_j$ in the no-frills form $$ Y'_j = uxF_{\Lambda x}\cap ux_jF_{\Delta _j}, $$ and recalling that $\ell (x)> \ell (x_j)$, it is immediate
that $Y'_j$ is empty unless $x_j$ is a prefix of $x$.  Note we may assume that $ux_j\in F_{\Delta _j}$.

Abandoning all of the empty $Y'_j$, which after all have no effect in the above definition of $U'$, we may assume that, for each $j$, there exists some $y_j$
in $\LX $ such that $x=x_jy_j$.  We may then see $Y'_j$ as the domain of the idempotent element of $\hull (S) $ given by $$ \theta _{ux}f_{\Lambda x}\theta
_{ux}\inv \theta _{ux_j}f_{\Delta _j}\theta _{ux_j}\inv = \theta _{ux_jy_j}f_{\Lambda x}\theta _{y_j}\inv \theta _{ux_j}\inv \theta _{ux_j}f_{\Delta
_j}\theta _{ux_j}\inv = \theta _{ux_jy_j}f_{\Lambda x}\theta _{y_j}\inv f_{ux_j}f_{\Delta _j}\theta _{ux_j}\inv \quebra = \theta _{ux_jy_j}f_{\Lambda
x}\theta _{y_j}\inv f_{\Delta _j}\theta _{ux_j}\inv \explica {\first {\FPCovar .ii}}= \theta _{ux_jy_j}f_{\Lambda x}f_{\Delta _jy_j}\theta _{y_j}\inv
\theta _{ux_j}\inv = \theta _{ux}f_{\Lambda x\cup \Delta _jy_j}\theta _{ux}\inv , $$ from where we deduce that $$ Y'_j = uxF_{\Lambda x\cup \Delta _jy_j}.  $$

Setting $u'=ux$, $\Lambda '=\Lambda x$, and $\Delta _j'=\Lambda x\cup \Delta _jy_j$, we then have that $\Lambda '\subseteq \Delta _j'$, while the
announced nicer description of the ingredients involed in the definition of $U'$ becomes $$ X'=u'F_{\Lambda '}, \and Y'_j=u'F_{\Delta '_j}.  $$

Having finally arrived at a convenient description of our open neighborhood of $\varphi $, we will now reset the notation so far introduced in this
proof by assuming that $$ X=uF_{\Lambda }, \and Y_j=uF_{\Delta _j}, $$ and that $U$ is defined as in \ref {DefnU} in terms of the above $X$ and $Y_j$.
We also assume that $\Lambda $, as well as the $\Delta _j$ are finite subsets of $\tS $ intersecting $S$ nontrivially, that $u\in \Lambda $, and that
$\Lambda \subseteq \Delta _j$, whence $Y_j\subseteq X$.

We next consider the question of whether or not $$ Z:=X\setminus \medcap _{j=1}^mY_j $$ is finite.  Should the answer be affirmative, the fact that
$\varphi $ is essentially tight relative to $\iota $ would imply that $$ 1=\varphi (X)= \bigvee _{j=1}^m\varphi (Y_j)=0, $$ a contradiction.  The answer
to the above question is therefore negative and hence $Z$ is an infinite set.  Letting $W=\theta _u\inv (Z)$, we then have that $W$ is also an infinite
set, while $$ W\subseteq F_\Lambda , \and W\cap F_{\Delta _j}=\emptyset , \for j=1, \ldots , m.  $$

Given any $s$ is $W$, one has that $s\notin F_{\Delta _1}$, which of course means that there is some $r$ in $\Delta _1$, such that $rs=0$.  Since there
are infinitely many elements in $W$, and only finitely many elements in $\Delta _1$, there necessarily exists a single $r_1\in \Delta _1$, and an
infinite subset $W_1\subseteq W$, such that $r_1s=0$, for all $s$ in $W_1$.

Similarly there exists some $r_2\in \Delta _2$, and an infinite subset $W_2\subseteq W_1$, such that $r_2s=0$, for all $s$ in $W_2$.  Continuing in
this way, we thus obtain a vector $$ \big (r_j\big )_{1\leq j\leq m}\in \medprod _{1\leq j\leq m}\Delta _j, $$ and an infinite subset $V\subseteq W$,
such that $$ r_js=0, \for j=1,\ldots ,m,\for s\in V.  $$

Setting $\Gamma =\{r_1, \ldots , r_m\}$, we then have that $V\subseteq F_{\Lambda ,\Gamma }$.  The grand conclusion so far is that $F_{\Lambda ,\Gamma
}$ is infinite, whence by $(*)$ we deduce the existence of a maximal string $\tau $ such that $$ \tau \in F^\star _ t , \for t\in \Lambda , \and \tau
\notin F^\star _ r , \for r\in \Gamma .  $$

Letting $\sigma =\theta ^\star _u(\tau )$, Proposition \first {\FPForwardInvariance } implies that $\sigma $ is a maximal, and $$ 1 = \bool {\tau \in
F^\star _ \Lambda } = \bool {\theta ^\star _u(\tau )\in \theta ^\star _u(F^\star _ \Lambda )} = \Bool {\sigma \in \varepsilon \big (uF^\theta _\Lambda
\big )} = \varphi _\sigma (uF^\theta _\Lambda ) = \varphi _\sigma (X).  $$

In order to show that $\varphi _\sigma $ lies in $U$ we must show that $\varphi _\sigma (Y_j)=0$, for all $j$, so we must compute $$ \varphi _\sigma
(Y_j) = \bool {\sigma \in \varepsilon \big (uF^\theta _{\Delta _j}\big )} = \bool {\theta ^\star _u(\tau )\in \theta ^\star _u(F^\star _ {\Delta _j})}
= \bool {\tau \in F^\star _ {\Delta _j}}.  $$ However notice that, since $r_j$ is in $\Delta _j$, one has that $$ \tau \notin F^\star _ {r_j} \supseteq
F^\star _ {\Delta _j}, $$ so indeed $\varphi _\sigma (Y_j) =0$.  This shows that $\varphi _\sigma $ is in $U$, hence proving that $\varphi $ is in the
closure of $\smaxShort $.  \endProof

\def \H {{\cal H}}

\section Matsumoto's C*-algebra

\label MatsuSection Let $\X $ be a subshift.  Although the empty word, henceforth denoted by $\varnothing $, is not allowed in $\LX $, it will nevertheless
play a role in what follows, so we will often refer to the extended language $$ \tilde \LX = \LX \cup \{\varnothing \}.  $$ Making the empty word act as a
unit for $\SX $, we may identify $\varnothing $ with the unit element of $\tSX $, so we have that $$ \tSX = \SX \cup \{\varnothing \}= \LX \cup \{0\}\cup
\{\varnothing \}.  $$ Removing zero from $\tSX $ we are left with the extended language, meaning that $$ \tSX '= \tSX \setminus \{0\} = \tilde \LX .  $$

As mentioned after \first {5.4}, a unitized semigroup may sometimes fail to be $0$-left-cancellati\-ve even when the original semigroup has this property.
However the present case is not an example of that undesirable situation, meaning that $\tSX $ is a well behaved $0$-left-cancellative semigroup,
as the reader may easily verify.

Consider the Hilbert space $$ \tilde \H = \ell ^2(\tilde \LX ) $$ formed by all square summable sequences $x = (x_\mu )_{\mu \in \tilde \LX }$ of
complex numbers.  The canonical orthonormal basis of $\tilde \H $ will be denoted by $\{\delta _\mu \}_{\mu \in \tilde \LX }$.  Given any $\nu $ in
$\tilde \LX $, consider the unique bounded linear operator $\tilde T _\nu $ on $\tilde \H $ such that $$ \tilde T _\nu (\delta _\mu ) = \left \{\matrix {
\delta _{\nu \mu }, & \hbox { if } \nu \mu \in \tilde \LX , \cr \pilar {12pt} 0, & \hbox { otherwise.} }\right .  $$

It is easy to see that $\tilde T _\varnothing $ is the identity operator and that each $\tilde T _\nu $ is a partial isometry.  In addition, for every
$\nu $ and $\mu $ in $\tilde \LX $, one has that $$ \tilde T _\nu \circ \tilde T _\mu = \left \{\matrix { \tilde T _{\nu \mu }, & \hbox { if } \nu \mu
\in \tilde \LX , \cr \pilar {12pt} 0, & \hbox { otherwise.} }\right .  $$ In particular, if $\nu =\nu _1\nu _2\ldots \nu _n$, with $\nu _i\in \Sigma $,
then $$ \tilde T _\nu = \tilde T _{\nu _1}\tilde T _{\nu _2}\ldots \tilde T _{\nu _n}.  \equationmark TNuInLetters $$

Setting $\tilde T _0=0$, we then obtain a multiplicative mapping $$ \nu \in \tSX \mapsto \tilde T _\nu \in B(\tilde \H ) .  \equationmark MatsuRep $$

Notice that a similar Hilbert space representation could be defined for any $0$-left-cancellative semigroup.

Since we will be working with many Hilbert space representations from now on, let us give the precise definitions below:

\definition \label HSRep Let $\H $ be any Hilbert space.  \iItemize \iItem A representation of a \underbar {semigroup} $S$ on $\H $ is any mapping
$\rho :S\to B(\H )$, whose range consists of partial isometries, satisfying $\rho (0)=0$, and $\rho (st)=\rho (s)\rho (t)$, for every $s$ and $t$
in $S$.  \iItem A representation of an \underbar {inverse semigroup} $\S $ on $\H $ is any mapping $\rho :\S \to B(\H )$, satisfying $\rho (0)=0$, $\rho
(st)=\rho (s)\rho (t)$, and $\rho (s\inv )=\rho (s)^*$, for every $s$ and $t$ in $\S $.  \iItem A representation of a \underbar {\sla } $\E $ on $\H $
is a representation of $\E $ in the sense of (ii), once $\E $ is seen as an inverse semigroup.

It is therefore clear that the mapping referred to in \ref {MatsuRep} is a representation of $\tSX $ on $\tilde \H $, in the sense of \ref {HSRep.i}.

\definition ({\it Cf.~}\cite {MatsuOri}) The \emph {Toeplitz-Matsumoto C*-algebra} $\tildeTX $ is the closed $*$-subal\-gebra of operators on $\tilde
\H $ generated by $\{\tilde T _\nu :\nu \in \tilde \LX \}$.

By \ref {TNuInLetters} one has that $\tildeTX $ may also be described as the smallest unital C*-algebra containing $\{\tilde T _a:a\in \Sigma \}$.

Denoting by $I$ the identity operator on $\tilde \H $, observe that $$ P:=I-\sum _{a\in \Sigma }\tilde T _a\tilde T _a^* $$ is the orthogonal projection
onto the one-dimensional subspace of $\tilde \H $ spanned by the so called \emph {vacuum vector} $\delta _\varnothing $, so it follows that $P$ lies
in $\tildeTX $.  Moreover, given any $\mu ,\nu ,\alpha $ in $\tilde \LX $, one has that $$ \tilde T _\mu P\tilde T _\nu ^*(\delta _\alpha ) = \left
\{\matrix { \delta _{\mu }, & \hbox { if } \alpha =\nu , \hfill \cr \pilar {12pt} 0, & \hbox { otherwise,} }\right .  $$ so $\tilde T _\mu P\tilde T
_\nu ^*$ is the rank one partial isometry from the vector $\delta _\nu $ to $\delta _\mu $.  So we see that the algebra $K(\tilde \H ) $ of all compact
operators on $\tilde \H $ is a sub-algebra of $\tildeTX $.

\definition \label MatsuDefinition ({\it Cf.~}\cite {MatsuOri}) The \emph {Matsumoto C*-algebra} $\MX $ is defined as the quotient $$ \MX = \tildeTX
/K(\tilde \H ) .  $$

Matsumoto's original definition was given in terms of a two-sided subshift, namely a closed subset $\X \subseteq \Sigma ^{\bf Z}$, invariant under
the bilateral shift map.  Compared to Matsumoto's, definition \ref {MatsuDefinition} is slightly more general in the sense that, for every two-sided
subshift $\X \subseteq \Sigma ^{\bf Z}$, the canonical projection of $\X $ on $\Sigma ^{\bf N}$ is a one-sided subshift to which \ref {MatsuDefinition}
attaches precisely the same C*-algebra as that defined by Matsumoto in \cite {MatsuOri}.  On the other hand, Matsumoto's original definition requires
essentially no modifications in order to apply to one-sided subshifts.

The reader is warned that there are numerous C*-algebras associated to a subshift in the literature, sometimes presented with conflicting notation.
See \cite [Section 7]{CarlsenSilvestrov} for a comparative study of these algebras.

\bigskip

As we shall see, the inclusion of the empty word $\varnothing $ in $\tilde \LX $, and the subsequent inclusion of $\delta _\varnothing $ as a basis
vector in $\ell ^2(\tilde \LX )$, causes many technical problems and introduces unwanted relationships between unrelated ingredients.  Consider, for
example, two subshifts $\X _1$ and $\X _2$ on disjoint alphabets $\Sigma _1$ and $\Sigma _2$.  Then clearly $$ \X :=\X _1\mathop {\dot \cup }\X _2 $$
is a subshift on the alphabet $\Sigma :=\Sigma _1\mathop {\dot \cup } \Sigma _2$.  Being made out of two totally unrelated parts, it is evident that
$\X $ is a reducible subshift in the sense that both $\X _1$ and $\X _2$ are invariant under the shift map.  However, letting $L_{\scriptscriptstyle \X
_1}$ and $L_{\scriptscriptstyle \X _2}$ be the languages of $\X _1$ and $\X _2$, respectively, the ``empty word on $\Sigma _1$'' gets confused with the
``empty word on $\Sigma _2$'', and in particular $\ell ^2(L_{\scriptscriptstyle \X _1})$ and $\ell ^2(L_{\scriptscriptstyle \X _2})$ are not invariant
under $\tildeTX $.  This is because if $\mu _1$ and $\mu _2$ are any two words in $L_{\scriptscriptstyle \X _1}$ and $L_{\scriptscriptstyle \X _2}$,
respectively, then the operator $$ \tilde T _{\mu _2} \tilde T _{\mu _1}^* $$ maps $\delta _{\mu _1}$ to $\delta _{\mu _2}$ (via $\delta _\varnothing
$), trespassing a boundary one would expect to exist between the two disjoint pieces of $\X $.

The presence of the empty word causes a further problem in that the representation $\tilde T $ mentioned in \ref {MatsuRep} might not extend to $\hull (\SX
) $.  To see where the problem lies, suppose that the language $\LX $ contains words $\mu _1$ and $\mu _2$ with $F_{\mu _1}\cap F_{\mu _2}=\emptyset $.
Then $$ \theta _{\mu _1} \theta _{\mu _2}\inv = \theta _{\mu _1}f_{\mu _1}f_{\mu _2} \theta _{\mu _2}\inv = \theta _{\mu _1}\id _\emptyset \theta _{\mu
_2}\inv , $$ so we see that $\theta _{\mu _1} \theta _{\mu _2}\inv $ is the empty map.  However, if $\tilde T $ could be extended to a $*$-representation,
say $\pi $, of $\hull (\SX ) $ on $\tilde \H $, then $$ 0 = \pi (\theta _{\mu _1} \theta _{\mu _2}\inv ) = \pi (\theta _{\mu _1})\pi (\theta _{\mu
_2})^* = \tilde T _{\mu _1}\tilde T _{\mu _2}^*, $$ but the latter is nonzero since $$ \tilde T _{\mu _1}\tilde T _{\mu _2}^*(\delta _{\mu _2})= \tilde
T _{\mu _1}(\delta _\varnothing )= \delta _{\mu _1}.  $$

Fortunately there exists an alternative construction of $\MX $ which does not involve the empty word nor the extended language $\tilde \LX $.  The key
observation is that $\ell ^2(\LX )$ is a subspace of $\ell ^2(\tilde \LX )$ of codimension one, and hence the compression to $\ell ^2(\LX )$ of any
operator $U$ on $\ell ^2(\tilde \LX )$ differs from $U$ by a compact perturbation.

To be more precise, let $$ \H =\ell ^2(\LX ), $$ and for each $\nu $ in $\LX $, consider the unique bounded linear operator $T_\nu $ on $\H $ such that
$$ T_\nu (\delta _\mu ) = \left \{\matrix { \delta _{\nu \mu }, & \hbox { if } \nu \mu \in \LX , \cr \pilar {12pt} 0, & \hbox { otherwise.} }\right .
\equationmark DefineMatsuRep $$ Notice that $T_\nu $ may also be viewed as the restriction of $\tilde T _\nu $ to $\H $.

As before, it is easy to see that, upon setting $T_0=0$, we obtain a multiplicative mapping $$ \nu \in \SX \mapsto T_\nu \in B(\H ).  \equationmark
MatsuRepNoUnit $$

\state Proposition \label IntroTX Letting $\TX $ be the closed $*$-subalgebra of operators on $\H $ generated by $\{T_\nu :\nu \in \LX \}$, one has
that $$ \MX \simeq {\TX + K(\H ) \over K(\H )}.  $$

\Proof Consider the $*$-homomorphism $\Phi $ defined as the composition $$ \Phi : \TX + K(\H ) \hookrightarrow B(\H ) \arw \iota B(\tilde \H ) \arw
q B(\tilde \H ) / K(\tilde \H ), $$ where the leftmost arrow is the inclusion of $\TX + K(\H )$ in $B(\H )$, the rightmost arrow is the quotient map,
and $\iota $ is the natural map extending any operator $U$ from $\H $ to $\tilde \H $ by setting $U(\delta _\varnothing )=0$.

Observing that, for every $\mu $ in $\LX $, one has that $$ \iota (T_\mu )|_\H =\tilde T _\mu |_\H , $$ and noticing that the codimension of $\H $ in
$\tilde \H $ is one, we see that $\iota (T_\mu )-\tilde T _\mu $ has rank at most one, and hence is a compact operator.  Consequently $$ \Phi (T_\mu )
= q(\tilde T _\mu ), \for \mu \in \LX .  $$

We wish to use the above conclusion to deduce that the range of $\Phi $ coincides with $q(\tildeTX )$, and hence also with $\MX $.  Notice that the
latter is generated by the $q(\tilde T _\mu )$, for $\mu $ not only in $\LX $, but in $\tilde \LX $.  Thus, to prove the desired coincidence of ranges
we need to show that $q(\tilde T _\varnothing )$, also known as the identity, is in the range of $\Phi $.  But this follows easily by noticing that $$
I_{\tilde \H }-\sum _{a\in \Sigma } \iota (T_a)\iota (T_a)^* $$ is a compact operator (not necessarily of rank one, it fixes $\delta _w$ for $w\in \LX
$ with $|w|\leq 1$ and annihilates all other $\delta _w$) and hence $$ 1=q(I_{\tilde \H }) = \sum _{a\in \Sigma } \Phi (T_a)\Phi (T_a)^*.  $$ We may
therefore view $\Phi $ as a surjective $*$-homomorphism $$ \Phi :\TX + K(\H ) \to \MX .  $$ Its kernel evidently contains $K(\H )$, and we claim that
$\Ker (\Phi )$ is in fact exactly equal to $K(\H )$.  To see this, let $U\in \TX + K(\H )$, and suppose that $\Phi (U)=0$.  Therefore $q(\iota (U))=0$,
so $\iota (U)$ is compact, but since $U$ coincides with the restriction of $\iota (U)$ to $\H $, we have that $U$ is also compact.  This proves the
claim, so $\Phi $ factors through the quotient of $\TX + K(\H )$ by $K(\H )$, providing the isomorphism sought.  \endProof

Unlike the representation $\tilde T $ discussed above, the representation $T$ of \ref {MatsuRepNoUnit} easily extends to $\hull (\SX ) $, as we shall
now prove.

\state Proposition \label ExistMatsuRep There exists a unique inverse semigroup representation $$ \pi : \hull (\SX ) \to B(\H ), $$ such that $\pi
(\theta _\mu ) = T_\mu $, for every $\mu $ in $\SX $.

\Proof It is easy to see that, given any $f$ in the symmetric inverse semigroup $\I (\LX )$, there exists a unique partial isometry $\tau (f)$ on $\H $
such that $$ \tau (f)\delta _\nu = \left \{\matrix { \delta _{f(\nu )}, & \hbox { if $\nu $ lies in the domain of } f, \cr \pilar {12pt} 0, & \hbox {
otherwise.}\hfill }\right .  $$ Moreover the correspondence $f\mapsto \tau (f)$ is clearly a representation of $\I (\LX )$ on $\H $.  The desired map
$\pi $ may then be obtained as the result of the composition $$ \hull (\SX ) \longhookrightarrow \I (\LX ) \arw \tau B(\H ).  $$ The uniqueness of $\pi $
follows from the fact that the $\theta _\mu $ generate $\hull (\SX ) $.  \endProof

Given our interest in the idempotent {\sla } of $\hull (\SX ) $, a detailed description of the restriction of $\pi $ to $\ehull (\SX ) $ will be useful.

\state Proposition \label PiOnConstruc Given any $Y$ in $\ehull (\SX ) $, one has that $\pi (Y)$ coincides with the orthogonal projection onto the
subspace $\ell ^2(Y)\subseteq \H $.

\Proof According to \first {\FPDefIdempSLA }, we are identifying $Y$ with $\id _Y$, so if we employ the representation $\tau $ used in the above proof,
we have that $$ \pi (Y)\delta _\nu = \tau (\id _Y)\delta _\nu = \bool {\nu {\in }Y}\delta _\nu , $$ so the range of $\pi (Y)$ is indeed $\ell ^2(Y)$,
as required.  \endProof

It is evident that the range of the representation $T$ of \ref {MatsuRepNoUnit} is contained in $\TX $, so the same clearly applies to the range of
$\pi $, whence the composition $q\circ \pi $ takes $\hull (\SX ) $ into $\MX $.

\definition \label DefMatsRep We will refer to the map $$ \rho :\hull (\SX ) \to \MX , $$ defined by $\rho =q\circ \pi $, as the Matsumoto representation
of $\hull (\SX ) $.  For the record, we notice that $\rho (\theta _\mu ) = T_\mu +K(\H )$, for every $\mu $ in $\SX $.

Strictly speaking $\rho $ is not a representation in the sense of \ref {HSRep.iii}, since it takes values in a C*-algebra, rather than in $B(\H )$.
However, every C*-algebra may be represented as an algebra of operators on some Hilbert space, so all of the concepts introduced in \ref {HSRep}
naturally extend to maps taking values in a C*-algebra.

We now wish to discuss the notion of the \emph {support} of a representation, so let us take a moment to consider an abstract inverse semigroup $\S $
(always assumed to have a zero element) and let $\pi $ be a representation of $\S $ in a C*-algebra $A$.  Letting $\E $ be the idempotent {\sla } of $\S
$, we have by \cite [10.6]{actions} that there exists a unique $*$-homomorphism $$ \Psi _\pi :C_0(\hat \E )\to A $$ such that $$ \Psi _\pi (1_e) = \pi
(e), \for e\in \E .  $$ The null space of $\Psi _\pi $ is evidently a closed 2-sided ideal of $C_0(\hat \E )$, and hence may be expressed as $C_0(U)$,
for some open subset $U\subseteq \hat \E $.  Following \cite [10.11]{actions}, the \emph {support} of $\pi $ is defined to be the subset of $\hat \E $
given by $$ \supp (\pi ) = \hat \E \setminus U.  $$ It is then easy to see that $\supp (\pi )$ consists precisely of those points $x$ in $\hat \E $
such that the evaluation character $\delta _x$ of $C_0(\hat \E )$ (given by $\delta _x(f)=f(x)$) vanishes on $\Ker (\Psi _\pi )$.  In symbols $$ \supp
(\pi ) = \big \{x\in \hat \E : \delta _x = 0 \text { on } \Ker (\Psi _\pi )\big \}.  \equationmark SpecInSymbols $$

Observe that the support of $\pi $ concerns only the idempotent {\sla } of $\S $. In particular the support of $\pi $ coincides with the support of
its restriction to $\E $.

We will soon describe the support of the representation $\rho $ of $\hull (\SX ) $ introduced in \ref {DefMatsRep}.  In doing so we will be aided by
the following general result:

\state Proposition \label SupportRelTigh Let $\pi $ be a representation of the semilattice $\E $ in a C*-algebra $A$.  Then the support of $\pi $
coincides with the set of all $\pi $-tight characters on $\E $.

\Proof This follows directly from \first {\FPFactorlemmaTwo }, by considering the Boolean algebra generated by the range of $\pi $ within $A$.
An alternative proof not using Boolean algebras is as follows.

Recall that the topology of $\hat \E $ is the topology of pointwise convergence, so the \emph {cylinder sets} $$ Z_{e_1,e_2,\ldots , e_n;f_1,f_2,\ldots
, f_m} := \Big \{\varphi \in \hat \E : \varphi (e_1)=\cdots =\varphi (e_n)=1,\ \varphi (f_1)=\cdots =\varphi (f_m)=0\Big \}, $$ where $e_1,e_2,\ldots ,
e_n, f_1,f_2,\ldots , f_m\in \E $, form a basis for the topology of $\hat \E $.

The above cylinder set is well defined even if $n$ (or $m$) vanishes but, should we take only those for which $n$ is nonzero, it is easy to see that
we still get a basis for the topology of $\hat \E $.

Fix, for the time being, a cylinder set as above, with $n>0$, and put $$ e'=\medprod _{i=1}^ne_i, \and f'_i = ef_i, \for i=1,\ldots ,m.  $$ We leave
it for the reader to verify that the cylinder set we fixed above coincides with $$ Z_{e';f'_1,f'_2,\ldots , f'_m}.  $$ The conclusion is that the
collection of all cylinder sets of the form $Z_{e;f_1,f_2,\ldots , f_m}$, whith $f_i\leq e$, also form a basis for the topology of $\hat \E $.

Changing subjects, let us now consider the set of all $(n+1)$-tuples $$ (e, f_1,f_2,\ldots ,f_n) \in \E ^{n+1} $$ such that $n\in {\bf N}$, $f_i\leq e$,
for all $i$, and $$ \pi (e) - \bigvee _{i=1}^n \pi (f_i) = 0.  \equationmark DataForGenerator $$ For each such $(n+1)$-tuple, one clearly has that the
element $b$ in $C_0(\hat \E )$ defined by $$ b:= 1_e - \bigvee _{i=1}^n 1_{f_i} \equationmark Generator $$ lies in $\Ker (\Psi _\pi )$ and we claim that
$\Ker (\Psi _\pi )$ coincides with the closed 2-sided ideal of $C_0(\hat \E )$ generated by the elements of the form $b$, as in \ref {Generator}, where
$e$ and the $f_i$ satisfy \ref {DataForGenerator}. In order to prove the claim, let $J$ be the latter ideal, so we evidently have that $J\subseteq \Ker
(\Psi _\pi )$.

Writing $J=C_0(V)$, for some open subset $V\subseteq \hat \E $, we then have that $V\subseteq U$, where $U$ is as before, namely $U = \hat \E \setminus
\supp (\pi )$.  Proving the claim therefore amounts to proving that $U\subseteq V$, so let us pick any $x$ in $U$.

Since $U$ is open, it is the union of basic open sets, so the first part of the present proof yields $(e, f_1,f_2,\ldots ,f_n)$ in $\E ^{n+1}$, with
$f_i\leq e$, such that $$ x\in Z_{e;f_1,f_2,\ldots , f_m}\subseteq U.  \equationmark HereIsX $$

With $b$ as in \ref {Generator}, one then checks that the support of $b$ coincides with $Z_{e;f_1,f_2,\ldots , f_m}$, so $$ b \in C_0(U) = \Ker (\Psi
_\pi ), $$ from where we deduce that $\Psi _\pi (b)=0$, whence \ref {DataForGenerator} holds.  It follows that $$ b \in J = C_0(V), $$ so the support of
$b$ is a subset of $V$, and we see that $x$ lies in $V$ by \ref {HereIsX}.  This concludes the proof of the claim that $J$ coincides with $\Ker (\Psi _\pi )$.

For each $(e, f_1,f_2,\ldots ,f_n)$ in $\E ^{n+1}$, with $f_i\leq e$, and for each $\varphi $ in $\hat \E $, define $b$ to be the right-hand-side of
\ref {Generator} and, considering the evaluation character $\delta _\varphi $ on $C_0(\hat \E )$, notice that $$ \delta _\varphi (b) = \delta _\varphi
\Big (1_e - \bigvee _{i=1}^n 1_{f_i}\Big ) = 1_e(\varphi ) - \bigvee _{i=1}^n 1_{f_i}(\varphi ) = \varphi (e) - \bigvee _{i=1}^n \varphi (f_i).  $$
As already mentioned, we have that $\varphi $ lies in $\supp (\pi )$ if and only $\delta _\varphi $ vanishes on $\Ker (\Psi _\pi )$, which is the case
if and only if the left-hand-side above vanishes whenever \ref {DataForGenerator} holds.  On the other hand, $\varphi $ is $\pi $-tight if and only if
the right-hand-side vanishes, so the proof is concluded.  \endProof

\state Proposition \label SupportEss The support of $\rho $ coincides with $\Ehat \sub {ess}(\SX ) $.

\Proof Using \ref {SupportRelTigh} it is enough to prove that $\Ehat \sub {ess}(\SX ) $ is exactly the set of all $\rho $-tight characters.  Recall that
a character $\varphi $ on $\ehull (\SX ) $ is $\rho $-tight if and only if $$ \varphi (X)=\bigvee _{i=1}^n\varphi (Y_i), $$ whenever $X, Y_1,\ldots
,Y_n\in \ehull (\SX ) $ are such that $$ \rho (X)=\bigvee _{i=1}^n\rho (Y_i).  $$ Since $\rho =q\circ \pi $, the equation displayed above is equivalent
to saying that $$ \pi (X) - \bigvee _{i=1}^n\pi (Y_i) \equationmark DifCpct $$ is a compact operator.  By \ref {PiOnConstruc} we have that $\pi (X)$
is the projection onto $\ell ^2(X)$, while $\bigvee _{i=1}^n\pi (Y_i)$ is the projection onto $$ \sum _{i=1}^n\ell ^2(Y_i) = \ell ^2(Y), $$ where
$Y=\bigcup _{i=1}^nY_i$.  Therefore the operator appearing in \ref {DifCpct} coincides with the projection onto $\ell ^2(X{\setminus } Y)\pilar {10pt}$
minus the projection onto $\ell ^2(Y{\setminus }X)$, which is a compact operator if and only if $X\Delta Y$ is finite.

This said, it is evident that $\varphi $ is $\rho $-tight if and only if $\varphi $ is essentially tight relative to the inclusion representation $\iota
$ of $\ehull (\SX ) $, which in turn is to say that $\varphi $ is in $\Ehat \sub {ess}(\SX ) $.  This concludes the proof.  \endProof

\def \F {{\bf F}} \def \TTensor {T^{\scriptscriptstyle \otimes }} 

\section The Carlsen-Matsumoto C*-algebra

As before, we fix a finite alphabet $\Sigma $ and a subshift $\X \subseteq \Sigma ^{\bf N}$.  We would now like to describe a variation of Definition \ref
{MatsuDefinition}, leading up to another C*-algebra closely related to $\X $, first introduced in \cite [Definition 6.4]{DokuchaExel}, and motivated by
\cite [Lemma 4.1]{MatsuAuto}\fn {We shoud however warn the reader that the isomorphism claimed in \cite [Lemma 4.1]{MatsuAuto} is incorrect.  See \cite
{MatsuCarl} for ways to correct this result.}  and \cite [Definition 2.1]{MatsuCarl}.

Since the present section will be largely independent of the previous one, we shall not make any attempt to distinguish the notations used here and there.
We therefore warn the reader that this section will share many symbols with section \ref {MatsuSection}, although they will often refer to completely
different objects.

While the description of the Matsumoto C*-algebra $\MX $ given in the previous section started out by considering the Hilbert space $\ell ^2(\tilde
\LX )$, and later $\ell ^2(\LX )$, here we will mostly work with operators on the Hilbert space $$ \H =\ell ^2(\X ).  $$ Since $\X $ is likely to be
an uncountable set, $\ell ^2(\X )$ might be a non-separable Hilbert space.  Nevertheless we will only consider separable C*-algebras of operators on
this oversized Hilbert space.

For each $\mu $ in $\SX $, consider the bounded linear operator $T_\mu $ on $\ell ^2(\X )$ defined on the standard orthonormal basis $\{\delta _\omega
\}_{\omega \in \X }$ by $$ T_\mu (\delta _\omega ) = \left \{\matrix { \delta _{\mu \omega }, & \hbox { if } \mu \omega \in \X , \cr \pilar {12pt} 0,
& \hbox { otherwise.} }\right .  $$ In the special case that $\mu $ is the zero element of $\SX $, we interpret the above in a somewhat ad-hoc way as
saying that $T_\mu =0$.

The term $\mu \omega $ employed above is meant to refer to the concatenation of the finite word $\mu $ with the infinite word $\omega $, evidently
resulting in an infinite word, namely $\mu \omega $.  The reader should notice the similarity of the above with \ref {DefineMatsuRep}, while we insist
that the present meaning of $T_\mu $ should not be confused with the one given by \ref {DefineMatsuRep}.

It is then evident that the correspondence $\mu \mapsto T_\mu $ defines a representation of $\SX $ on $\ell ^2(\X )$.

In our next concept we will employ the free group on the alphabet $\Sigma $, denoted by $\F $, and we will let $$ \lambda : \F \to B\big (\ell ^2(\F
)\big ) $$ be the left regular representation of $\F $.

\definition \cite [6.4]{DokuchaExel}.  The Carlsen-Matsumoto C*-algebra of $\X $, here denoted by $\OX $, is the closed $*$-subalgebra of operators on
$\ell ^2(\X )\otimes \ell ^2(\F )$ generated by the set $$ \{T_\mu \otimes \lambda _\mu : \mu \in \LX \}.  $$

We should perhaps clarify the meaning of the term $\lambda _\mu $, above: there is an obvious interpretation of finite words in $\Sigma $ as elements of
the group $\F $, and hence one may sensibly plug members of $\LX $ as parameters for the left regular representation.  However one should notice that,
unless $\X $ is the full shift, the correspondence $\mu \mapsto \lambda _\mu $ is not multiplicative if the multiplication operation considered in its
domain is that of $\SX $.  Nevertheless, we have:

\state Proposition The map $$ \TTensor :\SX \to B\big (\ell ^2(\X )\otimes \ell ^2(\F )\big ), $$ given by $\TTensor (\mu ) = T_\mu \otimes \lambda
_\mu $, for all $\mu $ in $\SX $, is a representation of $\SX $.

\Proof When $\mu $ is the zero element of $\SX $, the symbol $\lambda _\mu $, which has not yet been defined, will be understood\fn {In fact $\lambda
_0$ could be defined to be any operator whatsoever without affecting the term where it appears above.}  to mean the zero operator on $\ell ^2(\F )$.
Given $\mu $ and $\nu $ in $\SX $, we need to prove that $$ \TTensor (\mu \nu )=\TTensor (\mu )\TTensor (\nu ).  $$ The right-hand-side equals $$ \TTensor
(\mu )\TTensor (\nu ) = (T_\mu \otimes \lambda _\mu )(T_\nu \otimes \lambda _\nu ) = T_\mu T_\nu \otimes \lambda _\mu \lambda _\nu = T_{\mu \nu }\otimes
\lambda _{\mu \cdot \nu }, $$ where by $\mu \nu $ we denote the product of $\mu $ and $\nu $ in $\LX $, while $\mu \cdot \nu $ denotes their product
in $\F $.  Our task therefore consists in proving that $$ T_{\mu \nu }\otimes \lambda _{\mu \nu } = T_{\mu \nu }\otimes \lambda _{\mu \cdot \nu }.  $$
If $\mu \nu =0$, then $T_{\mu \nu }=0$, so both sides above vanish.  On the other hand, if $\mu \nu \neq 0$, then evidently $\mu \nu =\mu \cdot \nu $,
concluding the proof.  \endProof

We would now like to study a representation similar to the one given by \ref {DefMatsRep}, but with $\OX $ taking place of $\MX $.  In preparation for
this we will spend some time focusing on $T$ alone.

\state Proposition \label ExistCalMatsuRepOne There exists a unique inverse semigroup representation $$ \pi : \hull (\SX ) \to B(\ell ^2(\X )), $$
such that $\pi (\theta _\mu ) = T_\mu $, for every $\mu $ in $\SX $.

\Proof Let $\theta ^\star $ be the representation of $\SX $ on the set $\SX ^\star $ of all strings introduced in \first {10.9.ii}.  By \first
{\FPForwardInvariance } one has that $\SX ^\infty $ is invariant under (the forward action of) $\theta ^\star $, while \first {\FPBackInvarCorol }
implies that $\SX ^\infty $ is fully invariant under $\theta ^\star $, meaning that $\SX ^\infty $ is also invariant under $(\theta ^\star _\mu )\inv $,
for every $\mu $ in $\SX $.

Since $\hull (\SX ) $ is generated by the canonical image of $\SX $ under $\theta $, we conclude that $\SX ^\infty $ is also invariant under the
representation $\rho $ of $\hull (\SX ) $ on $\SX ^\star $ introduced in \first {\FPMapForISG }.

Recall from \ref {DescribeStringsAsWords} that the correpondence $\omega \mapsto \sigma _\omega $ defines a bijection between the set of all infinite
admissible words, also known as $\X $, and the set $\SX ^\infty $ of all maximal strings.  By identifying $\X $ with $\SX ^\infty $ under this
correspondence, we then get a representation $\psi $ of $\hull (\SX ) $ on $\X $, such that $$ \rho _\alpha (\sigma _\omega ) = \sigma _{\psi _\alpha
(\omega )}, \for \alpha \in \hull (\SX ) , \for \omega \in \X .  \equationmark RepHullX $$ The desired map $\pi $ may then be obtained as the result
of the composition $$ \hull (\SX ) \mathrel {\buildrel \psi \over \longhookrightarrow } \I (\X ) \arw \tau B(\ell ^2(\X ) ), $$ where $\tau $ is as in
the proof of \ref {ExistMatsuRep}, once $\X $ takes the place of $\LX $.  \endProof

Let us now present an abstract result designed to help determine the support of the above representation.

\state Proposition \label GeneralSupport Let $\S $ be an inverse semigroup with idempotent {\sla } $\E $.  Given a subset $Y\subseteq \hat \E $,
which is invariant in the sense of\/ \cite [10.7]{actions}, consider the representation $\pi $ of $\S $ on $\ell ^2(Y)$ given on the cannonical basis
$\{\delta _x\}_{x\in Y}$ by $$ \pi _s(\delta _x) = \clauses { \cl \delta _{\theta _s(x)} if x\in D_{s^*s}, \cl 0 otherwise. {} } $$ (In the above we
are using the notation employed in \cite [10.7]{actions}).  Then the support of $\pi $ coincides with the closure of\/ $Y$ in $\hat \E $.

\Proof Letting $\Psi _\pi $ be the representation of $C_0(\hat \E )$ as in \cite [10.6]{actions}, we claim that $$ \big \langle \Psi _\pi (f)\delta
_x,\delta _x\big \rangle = f(x), \equationmark DiagChar $$ for every $f$ in $C_0(\hat \E )$, and every $x$ in $Y$.  In order to prove this, assume
first that $f=1_e$, for some $e$ in $\E $.  Observing that $\pi _e$ is a diagonal idempotent operator whose range is $\ell ^2(D_e\cap Y)$, we have
that $$ \big \langle \Psi _\pi (1_e)\delta _x,\delta _x\big \rangle = \big \langle \pi (e) \delta _x,\delta _x\big \rangle = \bool {x\in D_e\cap Y}, $$
where the brackets denote boolean value.  On the other hand, we have $$ 1_e(x) \explica {[\withfont {cmbx6}{\possundef {actions}}:10.5]}= x(e) = \bool
{x\in D_e\cap Y}, $$ proving claim \ref {DiagChar} for the special case in which $f=1_e$.  Since the $1_e$ span a dense subspace of $C_0(\hat \E )$,
the claim follows.

From \ref {DiagChar} it then immediately follows that, for every $x$ in $Y$, the character on $C_0(\hat \E )$, given by evaluation at $x$, vanishes on
$\Ker (\Psi _\pi )$.  Therefore \ref {SpecInSymbols} implies that $Y\subseteq \supp (\pi )$.

In order to prove the opposite inclusion, observe that $\Psi _\pi (f)$ is a diagonal operator for every $f$ in $C_0(\hat \E )$.  Should $f$ vanish on $Y$,
then all of the diagonal entries of $\Psi _\pi (f)$ also vanish by \ref {DiagChar}, which is to say that $\Psi _\pi (f)=0$, and hence $f$ vanishes on
$\supp (\pi )$.  This proves that $$ f|_Y=0 \IMPLY f|_{\supp (\pi )}=0, $$ for every $f$ in $C_0(\hat \E )$ and this, together with Tietze's extension
Theorem, imply that $Y$ is dense in $\supp (\pi )$.  This concludes the proof.  \endProof

We are now in a position to determine the support of the representation $\pi $ introduced in \ref {ExistCalMatsuRepOne}.

\state Proposition \label SupportPi The support of $\pi $ coincides with the closure of the set $\Ehat \sub {max}(\SX )$ in $\Ehat (\SX ) $.

\Proof There are three representations of $\hull (\SX ) $ of interest to us here, namely: \iaitem \aitem the representation $\psi $ on $\X $ given in
\ref {RepHullX}, \aitem the representation $\rho $ on $\SX ^\star $ introduced in \first {\FPMapForISG }, and mentioned in \ref {RepHullX}, and \aitem
the canonical representation $\hat \theta $ on $\Ehat (\SX ) $ introduced in \cite [10.3.iv]{actions}.

\begingroup \noindent \hfill \beginpicture \setcoordinatesystem units <0.025truecm, -0.020truecm> \setplotarea x from -20 to 250, y from -60 to 180
\put {\null } at -20 -60 \put {\null } at 250 180 \put {$\X $} at 0 140 \put {$\buildrel \psi \over \curvearrowright $} at 0 105 \put {$\SX ^\infty
$} at 100 140 \put {$\SX ^\star $} at 100 0 \put {$\buildrel \rho \over \curvearrowright $} at 100 -35 \put {$(\SX ^\star )_\sharp $} at 100 70 \put
{$\Ehat (\SX )$} at 200 70 \put {$\buildrel \hat \theta \over \curvearrowright $} at 200 35 \put {$\Ehat \sub {max}(\SX )$} at 200 140 \arrow <0.15cm>
[0.25, 0.75] from 15.6 140 to 82 140 \put {$\scriptstyle \omega \mapsto \sigma _\omega $} at 48.8 127 \arrow <0.15cm> [0.25, 0.75] from 126.4 70 to
173.6 70 \put {$\scriptstyle \Phi $} at 150 57 \arrow <0.15cm> [0.25, 0.75] from 118 140 to 164 140 \put {\rotatebox {90}{$\subseteq $}} at 100 35 \put
{\rotatebox {90}{$\subseteq $}} at 100 105 \put {\rotatebox {90}{$\subseteq $}} at 200 105 \endpicture \hfill \null \endgroup

By \ref {RepHullX} we have that $\psi $ is equivalent to the restriction of $\rho $ to the invariant subset $\SX ^\infty $.  Furthermore the mapping $\Phi $
of \first {\FPCovarPhi }, defined on the $\rho $-invariant set $(\SX ^\star )_\sharp $ formed by the nondegenerate strings \first {\FPIntroPhiSigmaNonDeg
}, is covariant relative to the respective representations.  Although $\Phi $ is not necessarily injective, it has this property when restricted to the
open strings by \first {16.5}, and hence also when restricted to $\SX ^\infty $ by \first {\MaxOpenOrDeadEnd } because there are no nonzero elements
$\mu $ in $\SX $ such that $\mu \SX =\{0\}$.

Consequently the restriction of $\rho $ to $\SX ^\infty $ is equivalent to the restriction of $\hat \theta $ to the image of $\SX ^\infty $ under $\Phi
$, which turns out to be equal to $\Ehat \sub {max}(\SX )$ by definition.

The grand conclusion is then that $(\psi ,\X )$ and $\big (\hat \theta , \Ehat \sub {max}(\SX )\big )$ are equivalent systems and hence $\pi $, having
been built out of $\psi $ in \ref {ExistCalMatsuRepOne}, is unitarily equivalent to the representation arising from \ref {GeneralSupport} under the
choice of $Y=\Ehat \sub {max}(\SX )$.

Since unitarily equivalent representations obviously have the same support, the statement follows from \ref {GeneralSupport}.  \endProof

Even though the Carlsen-Matsumoto C*-algebra consists of operators on $\ell ^2(\X )\otimes \ell ^2(\F )$, since \ref {ExistCalMatsuRepOne} we have only
been discussing operators on $\ell ^2(\X )$.  In order to bring back the role of the space $\ell ^2(\F )$, we will now introduce an important map from
$\hull (\SX ) $ to the free group.

\state Proposition \label PartialHomomorphism There exists a map $$ d:\hull (\SX ) \to \F \cup \{0\}, $$ such that, for all $\alpha $ in $\hull (\SX )
$, one has \izitem \zitem $d(\alpha )=0$, if and only if $\alpha =0$, \zitem $d$ is \emph {idempotent pure}, in the sense that $d(\alpha )=1$, if and
only if $\alpha $ is a nonzero idempotent, \zitem $d$ is a \emph {partial homomorphism}, in the sense that $d(\alpha \beta )=d(\alpha )d(\beta )$,
whenever $\alpha $ and $\beta $ are elements of $\hull (\SX ) $ such that $\alpha \beta \neq 0$, \zitem $d(\theta _\mu )=\mu $, for every $\mu $ in $\SX $.

\Proof By \first {\FPFormOfHull }, every nonzero element $\alpha $ in $\hull (\SX ) $ may be written in a so called \emph {normal form} $$ \alpha =\theta
_\mu f_\Lambda \theta _\nu \inv , $$ where $\Lambda \subseteq \SX \cup \{1\}$ is a finite set, $\Lambda \cap \SX $ is nonempty, and and $\mu , \nu \in
\Lambda $.  This form is unique in the sense of \ref {SubshiftUniqueness}, so we may unambiguously define $$ d(\alpha ) = \clauses { \cl 0 if \alpha =0,
\cl {\mu \nu \inv } if \alpha =\theta _\mu f_\Lambda \theta _\nu \inv \neq 0.  } $$ We leave it for the reader to verify the validity of (i)--(iv).  \endProof

Inverse semigroups admiting a map $d$ into $G\cup \{0\}$, where $G$ is a group, and satisfying (i)--(iii) of \ref {PartialHomomorphism}, are said to
be \emph {strongly $0$-$E$-unitary} \cite {BFFG}, \cite {Lawson}, \cite {MilanSteinberg}.

For future reference we note the following:

\state Proposition \label HullEUnitary For every subshift $\X $, one has that $\hull (\SX ) $ is strongly $0$-$E$-unitary.

The above map $d$ is the key ingredient of a representation of crucial importance.

\state Proposition \label RepMax There exists a unique inverse semigroup representation $$ \rho : \hull (\SX ) \to \OX , $$ such that $\rho (\theta
_\mu ) = T_\mu \otimes \lambda _\mu $, for every $\mu $ in $\SX $.

\Proof It is enough to take $$ \rho (\alpha ) = \pi (\alpha )\otimes \lambda _{d(\alpha )}, \for \alpha \in \hull (\SX ) , $$ where $\pi $ is as in
\ref {ExistCalMatsuRepOne}.  Note that if $\alpha \beta =0$, then $\rho (\alpha )\rho (\beta ) = \pi (\alpha )\pi (\beta )\otimes \lambda _{d(\alpha
)d(\beta )} = 0 = \pi _{\alpha \beta }\otimes \lambda _{d(\alpha \beta )}$ \endProof

Our next goal is to determine the support of the above representation.

\state Proposition The support of the representation $\rho $ introduced above coincides with the closure of\/ $\Ehat \sub {max}(\SX )$ in $\Ehat (\SX ) $.

\Proof As already noticed, the support of the representation of an inverse semigroup coincides with the support of the restriction of that representation
to the idempotent {\sla }.  Observing that $$ \rho (\alpha ) = \pi (\alpha )\otimes 1, \for \alpha \in \ehull (\SX ) , $$ one may easily check that
$\rho $ and $\pi $ share supports, so the conclusion follows from \ref {SupportPi}.  \endProof

\section Amenability of the universal groupoid of a subshift semigroup

\noindent The universal groupoid of a countable inverse semigroup $\S $ with semilattice of idempotents $\E $ is the groupoid of germs for the action
$\hat \theta $ of $\S $ on $\hat \E $.  We shall show that the universal groupoid $\G $ of the inverse hull $\hull (\SX )$ of the subshift semigroup
$\SX $ associated to a subshift $\X \subseteq \Sigma ^{\bf N}$ is  amenable by proving that it is isomorphic to a Deaconu-Renault groupoid (also called a
semi-direct product groupoid in \cite [Definition 2.4]{cuntzlike},  where this notion was first introduced) for a certain local homeomorphism of $\Ehat
(\SX )$.  We proceed by first using the results of Milan and the second author \cite {MilanSteinberg} to realize $\G $ as the groupoid associated to
a certain partial action of the free group $\bf F$ on $\Sigma $ on $\Ehat (\SX )$.

If $G$ is a discrete group acting partially on a space $X$, the associated partial transformation groupoid $G\ltimes X$, first studied by Abadie \cite
{Abadie}, consists of all triples $(y,g,x)\in X\times G\times X$ with $x\in X_{g^{-1}}$ and $\alpha _g(x)=y$, where $\alpha _g\colon X_{g^{-1}}\to
X_g$, with the subspace topology induced by the product topology.  The product is given by $(y,g,x)(y',g',x') = (y,gg',x')$ if $x=y'$ and is undefined
otherwise.  The inverse is given by $(y,g,x)^{-1} = (x,g^{-1},y)$.  The unit space consists of the triples $(x,1,x)$ with $x\in X$ and hence can be
identified with $X$.  The groupoid $G\ltimes X$ is Hausdorff and \'etale.

If $\S $ is an inverse semigroup with $0$, then the universal group $G$ of $\S $ is the group with generating set $\{\gamma _s\mid s\in \S \setminus
\{0\}\}$ in bijection with the non-zero elements of $\S $ and relations of the form $\gamma _s\gamma _t=\gamma _{st}$ whenever $st\neq 0$.  The mapping
$\gamma \colon \S \to G\cup \{0\}$ given by $\gamma (s)= [\gamma _s]$ for $s\neq 0$ and $\gamma (0)=0$ is a partial homomorphism and is the universal
partial homomorphism from $\S $ to a group. Here $[\gamma _s]$ denotes the class of $\gamma _s$ in the group $G$.  In particular, $\S $ is strongly
$0$-$E$-unitary if and only if $\gamma $ is idempotent pure.  Note that it is immediate from the defining relations that if $A\subseteq \S $ generates
$\S $ as an inverse semigroup, then $G$ is generated by $\{[\gamma _a]\mid a\in A\}$.

If $\S $ is strongly $0$-$E$-unitary with universal group $G$, then Milan and the second author proved that the universal groupoid $\G $ of $S$ is
isomorphic to $G\ltimes \hat \E $ for a certain partial action of $G$ on $\hat \E $, where $\E $ is the semilattice of idempotents of $\S $.  If $g\in
G$, then the elements of $\gamma ^{-1}(g)$ are pairwise compatible and we define $\alpha _g=\bigvee _{s\in \gamma ^{-1}(g)}\hat \theta _s$, which is
a well-defined partial homeomorphism of $\hat \E $.   The isomorphism takes a germ $[s,x]$ to $(\alpha _{\gamma (s)}(x),\gamma (s),x)$.  Notice that
since each $\alpha _g$ is a join of actions of elements of $\S $, it follows that any invariant subspace of $\hat \E $ in $\G $ is also $G$-invariant
and so invariant in $G\ltimes \hat \E $.

We saw earlier, in \ref {HullEUnitary}, that $\hull (\SX )$ has an idempotent pure partial homomorphism $d\colon \hull (\SX )\to {\bf F}\cup \{0\}$ given
by $d(0)=0$ and $d(\theta _{u_1}f_{\Lambda _1} \theta _{v_1}\inv ) = u_1v_1\inv $. In particular, for $a\in \Sigma $, we have that $d(\theta _a)=a$.
We show that $d$ can be identified with the universal group partial homomorphism.  This is a consequence of the following simple proposition.

\state Proposition \label universalgroup Let $\S $ be an $A$-generated inverse semigroup and ${\bf F}$ the free group on $A$.  Suppose that there is
a partial homomorphism $d\colon \S \to {\bf F}\cup \{0\}$ such that $d(a)=a$ for all $a\in A$.  Then ${\bf F}$ is isomorphic to the universal group of
$\S $ and $d$ can be identified with the universal partial homomorphism.

\Proof Let $\gamma \colon \S \to G\cup \{0\}$ be the universal partial group homomorphism.  We already observed that $G$ is generated by $\{[\gamma
_a]\mid a\in A\}$.  There is a homomorphism $\psi \colon G\to {\bf F}$ such that $d=\psi \gamma $.  In particular, $a=d(a) = \psi ([\gamma _a])$
for $a\in A$.  Since $\bf F$ is free we can find a homomorphism $\phi \colon {\bf F}\to G$ with $\phi (a) = [\gamma _a]$.  It follows that $\psi $
and $\phi $ are inverse isomorphisms.  \endProof

In particular, we can identify $d\colon \hull (\SX )\to {\bf F}\cup \{0\}$ with the universal group homomorphism, and hence the universal groupoid
$\G $ of $\hull (\SX )$ is isomorphic to ${\bf F}\ltimes \Ehat (\SX )$ and each reduction of $\G $ is a reduction of ${\bf F}\ltimes \Ehat (\SX )$.
Let us state a consequence of this result.

\state Corollary \label partialgroup  The universal groupoid of $\hull (\SX )$ is the partial transformation groupoid of a partial action of the free
group on the alphabet of $\X $, as is the reduction to any invariant subspace of the universal groupoid.

The action of ${\bf F}$ on $\Ehat (\SX )$ has two additional properties, studied by the first author \cite {ortho}, namely that of being semi-saturated
and orthogonal.

Let ${\bf F}$ be a free group with basis $\Sigma $.  Let $\alpha $ be a partial action of ${\bf F}$ on a space $X$ with $\alpha _g\colon X_{g\inv }\to
X_g$ for each $g\in {\bf F}$.

\iItemize \iItem The action is semi-saturated if whenever the concatenation $uv$ of $u,v\in {\bf F}$ is reduced as written, then $\alpha _u\alpha
_v=\alpha _{uv}$.  \iItem The action is orthogonal if $X_a\cap X_b=\emptyset $ for all $a,b\in \Sigma $ with $a\neq b$.

\state Proposition \label IsSemiSat The partial action of\/ ${\bf F}$ on $\Ehat (\SX )$ is semi-saturated and orthogonal.

\Proof From the definition of $d$, the only elements of ${\bf F}$ with non-empty domain are those of the form $uv\inv $ with $u,v\in \Sigma ^*$.
Moreover, we may assume that this product is reduced as written.  It follows from the normal form in \ref {SubshiftUniqueness} that $\theta _u\theta
_v\inv $ is the unique maximum element of $d\inv (uv\inv )$ (in the natural partial order) under the assumption $uv\inv $ is reduced and so $\alpha
_{uv\inv } = \hat \theta _{uv\inv }= \hat \theta _u\hat \theta _v\inv $.  First we check the semi-saturated condition.  Let $u_1,v_1,u_2,v_2\in \Sigma
^*$ with $u_1v_1\inv $ and $u_2v_2\inv $ in ${\bf F}$ reduced as written and assume that $u_1v_1\inv u_2v_2\inv $ is also a reduced word.  Assume first
that both $v_1$ and $u_2$ are non-empty.  Then the fact that $v_1\inv u_2$ is reduced means that the first letter of $v_1$ differs from the first letter
of $u_2$.  Then $\alpha _{u_1v_1\inv u_2v_2\inv }=0$ and $$ \alpha _{u_1v_1\inv }\alpha _{u_2v_2\inv } = \hat \theta _{u_1}\hat \theta _{v_1}\inv \hat
\theta _{u_2}\hat \theta _{v_2\inv } =0 $$ as required.  Next assume that at least one of $v_1$ and $u_2$ are empty.  By taking inverses if necessary we
may assume that $v_1$ is empty.  Then $u_1u_2v_2\inv $ is reduced as written and $$ \alpha _{u_1}\alpha _{u_2v_2\inv } = \hat \theta _{u_1}\hat \theta
_{u_2}\hat \theta {v_2}\inv = \hat \theta _{u_1u_2}\hat \theta {v_2}\inv = \alpha _{u_1u_2v_2\inv }.  $$ We conclude that the action is semi-saturated.

If $a,b\in \Sigma $ with $a\neq b$, then $\alpha _a=\hat \theta _a$ and $\alpha _b=\hat \theta _b$.  The range of $\hat \theta _a$ consists of those
characters $\phi $ with $\phi (E_a) = 1$ and the range of $\hat \theta _b$ consists of those characters $\phi $ with $\phi (E_b)=1$.  Since $E_a\cap
E_b=\emptyset $, we conclude that these partial homeomorphisms have disjoint range, whence $X_a\cap X_b=\emptyset $.  Thus the action is orthogonal.

\endProof

We aim to prove that groupoids of the from ${\bf F}\ltimes X$ for a semi-saturated and orthogonal partial action of ${\bf F}$ on a locally compact
Hausdorff spaces $X$ are isomorphic to Deaconu-Renault groupoids.  Our first step will be to show that the partial action looks very similar to the
way it did in the case of the inverse hull of a subshift semigroup.

\state Proposition \label DescribeSemisat Let $\alpha $ be a semi-saturated orthogonal action of\/ ${\bf F}$ on a locally compact Hausdorff space $X$.
Then every element of ${\bf F}$ with non-empty domain is of the form $uv\inv $ with $u,v\in \Sigma ^*$ and the product reduced as written.

\Proof Let $w=a_1^{\epsilon _1}\cdots a_n^{\epsilon _n}$ with $a_i\in \Sigma $ and $\epsilon _i\in \{\pm 1\}$ be a reduced word.  Then since the action
is semi-saturated, we have that $\alpha _w = \alpha _{a_1}^{\epsilon _1}\cdots \alpha _{a_n}^{\epsilon _n}$.  If there is an index $i$ with $\epsilon
_i=-1$ and $\epsilon _{i+1}=1$, then since $w$ is reduced $a_i\neq a_{i+1}$.  But orthogonality then implies that $\alpha _{a_i}^{\epsilon _i}\alpha
_{i+1}^{\epsilon _i}=\alpha _{a_i}\inv \alpha _{a_{i+1}} =0$.  Thus $w$ is of the form $uv\inv $ with $u,v\in \Sigma ^*$ and the product reduced
as written.  \endProof

If ${\bf F}$ has a semi-saturated and orthogonal action $\alpha $ on a locally compact space $X$, then the $\alpha _a$ with $a\in \Sigma $ have disjoint
domains and so $T = \bigvee _{a\in \Sigma } \alpha _a\inv = \coprod _{a\in \Sigma } \alpha _a\inv $ is a well-defined local homeomorphism $\coprod
_{a\in A} X_a\to X$.  Thus we can form the Deaconu-Renault groupoid $\G _{(T,X)}$.  Let us recall the definition.  The underlying set of $\G _{(T,X)}$
consists of all triples $(y,m-n,x)\in X\times {\bf Z}\times X$ such that $T^m(y)=T^n(x)$ with $m,n\geq 0$.  The product is defined by $(y,k,x)(y',k',x')
= (y,k+k',x')$ if $x=y'$ and is undefined otherwise.  The inverse is given by $(y,k,x)\inv = (x,-k,y)$.  The topology has as a basis all triples $$(U,
(m,n), V) = \{(y,m-n,x)\in \G _{T,X}\mid y\in U, x\in V\}$$ with $U,V$ open in $X$ and $T^m(U)=T^n(V)$ with $T^m|_U$ and $T^n|_V$ a homeomorphism.
The unit space of $\G _{(T,X)}$ consists of those triples $(x,0,x)$ with $x\in X$ and can be identified with $X$.

We aim to prove that ${\bf F}\ltimes X$ is isomorphic to $\G _{(T,X)}$.  First we prove some elementary properties of $T$.

\state Proposition \label TProperties Let ${\bf F}$ have a semi-saturated and orthogonal action on a locally compact Hausdorff space $X$.

\iaitem \aitem If $u,v\in \Sigma ^*$ with $|u|\leq |v|$, then $X_u\cap X_v\neq \emptyset $ implies that $u$ is a prefix of $v$, in which case
$X_v\subseteq X_u$.  \aitem $T^n =\coprod _{w\in \Sigma ^n} \alpha _w\inv $.

\Proof \itmproof (a) We prove that $u$ is a prefix of $v$ by induction on $|v|$.  If $|v|=0$, then $u=v$ is empty.  Assume true when the longer word has
length $n$ and assume $|v|=n+1$. Write $u=au'$ and $v=bv'$ with $a,b\in \Sigma $ and note $|u'|\leq |v'|=n$.  Then since the action is semi-saturated
and $u,v$ are reduced, we have that $\alpha _u = \alpha _a\alpha _{u'}$ and $\alpha _v=\alpha _b\alpha _{v'}$.  If $a\neq b$, then the assumption that
$X_a\cap X_b=\emptyset $ implies that $X_u\cap X_v\subseteq X_a\cap X_b=\emptyset $.  Thus we must have $a=b$.  Then if $x\in X_u\cap X_v$, we have
$\alpha _a\inv (x)\in X_{u'}\cap X_{v'}$  and so $u'$ is a prefix of $v'$ by induction.  Writing $v'=u'w$ we then have that $v=av'=au'w=uw$ and so $u$
is a prefix of $v$.    Also, since $uw$ is reduced and the action is semi-saturated, $\alpha _v= \alpha _u\alpha _w$ and so $X_v\subseteq X_u$.

\itmproof (b) Since the symmetric inverse monoid on $X$ is distributive, it follows that $T^n =\bigvee _{w\in \Sigma ^n}\alpha _w\inv $.  But by (i),
the collection $\{X_w\mid w\in \Sigma ^n\}$ consists of pairwise disjoint subsets and so $T^n = \coprod _{w\in \Sigma ^n} \alpha _w\inv $.

\endProof

We define a mapping $\Phi \colon {\bf F}\ltimes X\to \G _{(T,X)}$ by $$ \Phi (y,uv\inv ,x) = (y,|u|-|v|,x), $$ where $uv\inv $ is reduced with $u,v\in
\Sigma ^*$.  Let us check that this is well-defined.  Indeed, since $y=\alpha _{uv\inv }(x) = \alpha _u\alpha _v\inv (x)$, as the action is semi-saturated, we
deduce that $\alpha _u\inv (y) = \alpha _v\inv (x)$.  By \ref {TProperties} we deduce that $T^{|u|}(y) = \alpha _u\inv (y) = \alpha _v\inv (x) = T^{|v|}(x)$.

\state Theorem \label isowithDR Let $\alpha $ be a semi-saturated and orthogonal action of\/ ${\bf F}$ on a locally compact Hausdorff space and put
$T=\coprod _{a\in \Sigma }\alpha _a\inv $.  Then ${\bf F}\ltimes X\cong \G _{(T,X)}$.

\Proof We check first that $\Phi $ is a homomorphism.  Note that $\Phi (x,1,x) = (x,0,x)$ and so $\Phi $ sends the unit space of ${\bf F}\ltimes
X$ homeomorphically to the unit space of $\G _{(T,X)}$.  We next check that $\Phi $ is injective.  Suppose that $\Phi (y_1,u_1v_1\inv ,x_1) = \Phi
(y_2,u_2v_2\inv ,x_2)$ where $u_1v_1\inv $ and $u_2v_2\inv $ are reduced as written. First note that $y_1=y_2$ and $x_1=x_2$.  Also $|u_1|-|v_1|=|u_2|-|v_2|$.
Without loss of generality, assume that $|u_1|\leq |u_2|$.  Then $0\leq |u_2|-|u_1| = |v_2|-|v_1|$ and so $|v_1|\leq |v_2|$.  Note that since the action
is semi-saturated we have $y_1=\alpha _{u_i}\alpha _{v_i}\inv (x_1)$ for $i=1,2$ and so $\alpha _{u_i}\inv (y_1) = \alpha _{v_i}\inv (x_1)$ for $i=1,2$.
Thus $y_1\in X_{u_1}\cap X_{u_2}$ and $x_1\in X_{v_1}\cap X_{v_2}$.  We conclude that $u_2=u_1u$ and $v_2=v_1v$ for some words $u,v\in \Sigma ^*$ with
$|u|=|u_2|-|u_1|=|v_2|-|v_1| = |v|$. Then $$\alpha _{u_1}\alpha _{v_1}\inv (x_1) = \alpha _{u_2}\alpha _{v_2}\inv (x_1) =\alpha _{u_1}\alpha _u\alpha
_v\inv \alpha _{v_1}\inv (x_1)$$ and so $\alpha _{v_1}\inv (x_1) = \alpha _u\alpha _v\inv \alpha _{v_1}\inv (x_1)$.  There for $\alpha _{v_1}\inv
(x_1)\in X_u\cap X_v$ and so $u=v$ by (a) of \ref {TProperties} since $|u|=|v|$ and hence they are both prefixes of each other.  But then $u_2v_2\inv =
u_1uu\inv v_1\inv $ is not reduced as written unless $u$ is empty, that is, $u_1=u_2$ and $v_1=v_2$.  This establishes that $\Phi $ is injective.

Observe that a triple $(y,k,x)\in \G _{(T,X)}$ may have many representations of the form $(y,m-n,x)$ with $m,n\geq 0$ and $T^m(y)=T^n(x)$.  Let us
assume that $m$ is chosen to be minimum among all such representations.  Then note that $n$ is minimum as well since if $m-n=k=m'-n'$ with $m\leq m'$,
then $0\leq m'-m = n'-n$ and so $n\leq n'$.  We will say that $(y,m-n,x)$ is a minimum representative if $m$ and $n$ are minimum so that $T^m(y)=T^n(x)$.
If $(y,k,x)\in \G _{(T,X)}$ has minimum representative $(y,m-n,x)$, then by \ref {TProperties} we have unique words $u\in \Sigma ^m$ and $v\in \Sigma ^n$
such that $\alpha _u\inv (y)=T^m(y)=T^n(x) = \alpha _v\inv (x)$.  Then $y=\alpha _u\alpha _v\inv (x)$.  We claim that $uv\inv $ is reduced as written.
Otherwise, there is a non-empty word $w\in \Sigma ^+$ with $u=u'w$ and $v=v'w$.  But then $\alpha _u=\alpha _{u'}\alpha _w$ and $\alpha _v=\alpha
_{v'}\alpha _w$ and hence $$\alpha _w\inv \alpha _{u'}\inv (y)=\alpha _u\inv (y) = \alpha _v\inv (x) = \alpha _w\inv \alpha _{v'}\inv (x)$$ and so
$T^{|u'|}(y)=\alpha _{u'}\inv (y) = \alpha _{v'}\inv (x)=T^{|v'|}(x)$ with $|u'|<|u|$ and $|v'|<|v|$ contradicting that $(y,m-n,x)$ was a minimum
representative.  Since $uv\inv $ is reduced as written, we have that $y=\alpha _u\alpha _v\inv (x) = \alpha _{uv\inv }(x)$ by the semi-saturated property
of the action and so $(y,m-n,x) = (y,|u|-|v|,x) = \Phi (y,uv\inv ,x)$.  Thus $\Phi $ is surjective.

We check that $\Phi $ is a functor.  Since it is bijective on the unit space, it suffices to consider the effect of $\Phi $ on a product of the form
$(y,u_1v_1\inv ,z)(z, u_2v_2\inv x)$ with $u_iv_i\inv $ reduced as written, $u_i,v_i\in \Sigma ^*$ for $i=1,2$.   The fact that $u_1v_1\inv u_2v_2\inv
$ is defined at $x$ means that either $v_1$ is a prefix of $u_2$ or $u_2$ is a prefix of $v_1$.  We handle the first case, as the other is similar.
So write $u_2=v_1w$.  Then $u_1v_1\inv u_2v_2\inv = u_1v_1\inv v_1wv_2\inv = u_1wv_2\inv $ and the right hand side is reduced as written.  Thus we have
$\Phi ((y,u_1v_1\inv ,z)(z, u_2v_2\inv x))=\Phi (y,u_1wv_2\inv ,x) = (u, |u_1|+|w|-|v_2|)$. On the other hand, $\Phi (y,u_1v_1\inv ,z)\Phi (z,u_2v_2\inv
,x) = (y,|u_1|-|v_1|,z)(z, |v_1|+|w|-|v_2|,x) = (y, |u_1|+|w|-|v_2|,x)$ as required.

It is straightforward to verify that $\Phi $ is a homeomorphism.  For example, to see that $\Phi $ is open, a basic neighborhood of ${\bf F}\ltimes
X$ is of the form $ \alpha _g(U)\times \{g\}\times U$ where $U\subseteq X_{g\inv }$.  If $g=uv\inv $ with $u,v\in \Sigma ^*$ and the product reduced
as written, then the image of this neighborhood under $\Phi $ is the basic open set $(\alpha _g(U),(|u|,|v|), U)$ of $\G _{(T,X)}$.  If $(U,(m,n),V)$
is a basic neighborhood of $\G _{(T,X)}$, then its preimage under $\Phi $ is the union of all open sets of the form $\alpha _{uv\inv }(U\cap X_{uv\inv
})\times \{uv\inv \}\times U\cap X_{uv\inv }$ with $u\in \Sigma ^m$ and $v\in \Sigma ^n$ (we do not require $uv\inv $ to be reduced).  \endProof

Since Deaconu-Renault groupoids of local homeomorphisms are always amenable by \cite [Proposition 2.4]{cuntzlike}, we obtain the following corollaries.

\state Corollary \label semisatamenable Let $\alpha $ be a semi-saturated and orthogonal partial action of a free group ${\bf F}$ on a locally compact
Hausdorff space $X$.  Then the groupoid ${\bf F}\ltimes X$ is an amenable Hausdorff \'etale groupoid.

\state Corollary \label amenablesubshift Let $\X $ be a subshift.  Then the universal groupoid of $\hull (\SX )$ is an amenable Hausdorff groupoid and
hence so are all its reductions.

\section Groupoid models for Matsumoto's C*-algebras

\noindent Recall from \cite [Theorem 4.4.1]{Paterson} that, given a countable inverse semigroup $\S $, with idempotent {\sla } $\E $, the groupoid
of germs for the canonical action $\hat \theta $ of $\S $ on $\hat \E $ (called the universal groupoid), say $\G $, is a groupoid model for $C^*(\S
)$, the universal C*-algebra of $\S $, in the sense that $$ C^*(\S ) \simeq C^*(\G ), $$ where $C^*(\G )$ denotes the groupoid C*-algebra of $\G $.
Also the reduced $C^*$-algebra of $\G $ coincides with the reduced $C^*$-algebra of $\S $.

The unit space of $\G $ naturally identifies with $\hat \E $ so, given any closed, invariant subset $W\subseteq \hat \E $, one may consider the groupoid
of germs for the restriction of $\hat \theta $ to $W$ or, equivalently, the reduction of $\G $ to $W$, which we will henceforth denote by $\G ^W $.

A well known example of such an invariant subset is $W=\hat \E \sub {tight} $ \cite [Proposition 12.11]{actions}, and the groupoid obtained by reducing
$\G $ to $\hat \E \sub {tight} $ is a model in the above sense for $C^*\sub {tight} (\S )$, the tight C*-algebra of $\S $ \cite [Proposition 13.3]{actions}.

If $W\subseteq \hat \E $ is invariant, but not necessarily closed, then its closure is also invariant, so the above game may be played with $\overline W$.
A natural example is the set $\hat \E _\infty $ formed by all ultra-characters, whose closure is well known to be $\hat \E \sub {tight} $.

Here we wish to consider this circle of ideas applied to the inverse semigroup $\S =\hull (S) $, for a given $0$-left cancellative semigroup $S$
admitting least common multiples.  One of the distinctive features of this situation is that, besides the space of ultra-characters and the space of
tight characters, there are many other interesting invariant subsets of $\Ehat (S)$. Since $\G $ is amenable in this case by \ref {amenablesubshift}
we note that the $C^*(\G )$ is isomorphic to the reduced $C^*$-algebra of $\G $ and so we do not distinguish them.

\state Proposition \label ManyInvariant Let $S$ be a $0$-left-cancellative semigroup admitting least common multiples.  Then the following subsets of\/
$\Ehat (S)$ are invariant under the dual action of\/ $\hull (S) $: \iItemize \iItem $\Ehat \sub {ess}(S) $, \iItem the set formed by all open characters
\first {16.16}, \iItem $\{\varphi _\sigma :\sigma \hbox { is an open string}\}$.

\Proof Since $\hull (S) $ is generated by the $\theta _s$ and the $\theta _s\inv $, which are respectively represented as $\dualrep _s$ and $\dualrep
_s\inv $ on $\Ehat (S)$, in order to show that any given subset of $\Ehat (S)$ is invariant under $\hat \theta $, it suffices to prove that said set
is invariant under $\dualrep _s$ and $\dualrep _s\inv $, for every $s$ in $S$.

\itmproof (i) Given $s$ in $S$, pick $\varphi $ in $\Ehat \sub {ess}(S) $ lying in the domain of $\dualrep _s$, and let us prove that $\ac s\varphi
$ also belongs to $ \Ehat \sub {ess}(S) $. For this we choose $\theta $-constructible sets $X, Y_1,\ldots ,Y_n$ such that $ X\mathop {\Delta }\big
(\medcup _{i=1}^nY_i\big ) $ is finite, whence so is $$ \theta _s\inv \Big (E^\theta _s\cap \big ( X\mathop {\Delta }\big (\medcup _{i=1}^nY_i\big
)\big )\Big ) = X^s\mathop {\Delta }\big (\medcup _{i=1}^nY_i^s\big ), $$ where $X^s=\theta _s\inv (E^\theta _s\cap X)$, and $Y_i^s=\theta _s\inv
(E^\theta _s\cap Y_i)$.  We then have that $$ \ac s\varphi (X) = \varphi \big (\theta _s\inv (E^\theta _s\cap X)\big ) = \varphi \big (X^s) \quebra =
\bigvee _{i=1}^n\varphi (Y^s_i) = \bigvee _{i=1}^n \varphi \big (\theta _s\inv (E^\theta _s\cap Y_i)\big ) = \bigvee _{i=1}^n \ac s\varphi (Y_i), $$
as desired.  In a similar way one proves that $\acinv s\varphi $ is in $ \Ehat \sub {ess}(S) $, whenever $\varphi $ lies both in $ \Ehat \sub {ess}(S)
$ and in the range of $\dualrep _s$.

\itmproof (ii) Given $s$ in $S$, pick an open character $\varphi $ in the domain of $\dualrep _s$, and let us prove that $\ac s\varphi $ is open.  Since
$\varphi $ is open we have in particular that $\sigma _\varphi $ is nonempty, whence we may use \first {\FPBirthOfString .i} to conclude that $\sigma
_\varphi \in F^\star _s$, and $$ \sigma _{\ac s\varphi } = \theta ^\star _s(\sigma _\varphi ).  $$ It then follows from \first {\FPOpenInvarUnderTheta
.i} that $\theta ^\star _s(\sigma _\varphi )$ is an open string, whence $\ac s\varphi $ is an open character.

Next let us suppose that $\varphi $ is an open character in the range of $\dualrep _s$.  Then $\varphi (E^\theta _s)=1$, so $s\in \sigma _\varphi $.
Once we know that $\varphi $ is open, it follows that $s$ is in fact in the interior of $\sigma _\varphi $, so \first {\FPDualBackOnStrings .i} implies
that $\sigma _\varphi \in E^\star _s$, and $$ \sigma _{\acinv s\varphi } = \theta ^{\star -1}_s(\sigma _\varphi ).  $$ It then follows from \first
{\FPOpenInvarUnderTheta .ii} that $\theta ^{\star -1}_s(\sigma _\varphi )$ is an open string, whence $\acinv s\varphi $ is an open character.

\itmproof (iii) Recall from \first {\FPOpenStringsInvariant } that the subset of $S^\star $ formed by all open strings is invariant under the action
of $\hull (S) $, so the conclusion follows from \first {\FPCovarPhi }.  \endProof

Even though the question of invariance for the set $S^\infty $ of all maximal strings is a bit touchy (see \first {10.19} and the subsequent discussion),
for subshift semigroups we may add $ \Ehat \sub {max}(S) $ to the list above:

\state Proposition Let $\X $ be a subshift.  Regarding the associated semigroup $\SX $, one has that $ \Ehat \sub {max}(\SX ) $ is $\hat \theta $-invariant.

\Proof By \ref {DescribeStringsAsWords} we have that a string is maximal if and only if it is open.  Therefore $ \Ehat \sub {max}(\SX ) $ is invariant
because it coincides with the set appearing in \ref {ManyInvariant.iii}.  \endProof

The Matsumoto and the Carlsen-Matsumoto C*-algebras are among the C*-algebras of groupoids obtained as the reduction of the universal groupoid to
suitable invariant subsets, as we shall now show.

\def \Gess {\G _\X ^{^{\hbox {\sixrm ess}}}\pilar {10pt}} \def \Gmax {\G _\X ^{^{\hbox {\sixrm max}}}}

\state Theorem \label Main Given a subshift $\X $, let $\G _\X $ be the groupoid of germs for the canonical action $\hat \theta $ of\/ $\hull (\SX ) $
on $\Ehat (\SX ) $.  \izitem \zitem Letting $\Gess $ be the reduction of\/ $\G _\X $ to $ \Ehat \sub {ess}(\SX ) $, one has that $$ \MX \simeq C^*\big
(\Gess \big ).  $$ \zitem Letting $\Gmax $ be the reduction of\/ $\G _\X $ to $\overline {\Ehat \sub {max}(\SX )} $, one has that $$ \OX \simeq C^*\big
(\Gmax \big ).  $$

\Proof Focusing on (i), and plugging the representation $\rho $ of \ref {DefMatsRep} in \cite [10.14]{actions}, we obtain a $*$-homomor\-phism $\Phi $
from $C^*\big (\Gess \big )$ to $\MX $, such that the diagram

\begingroup \noindent \hfill \beginpicture \setcoordinatesystem units <0.025truecm, -0.020truecm> \setplotarea x from -60 to 200, y from -30 to 100
\put {\null } at -60 -30 \put {\null } at 200 100 \put {$\hull (\SX )$} at 0 0 \put {$\MX $} at 140 0 \put {$C^*\big (\Gess \big )$} at 70 80 \arrow
<0.15cm> [0.25, 0.75] from 28.8 0 to 111.2 0 \put {$\rho $} at 70 -13 \arrow <0.15cm> [0.25, 0.75] from 89.222 58.032 to 123.312 19.072 \put {$\Phi $}
at 116.051 47.112 \arrow <0.15cm> [0.25, 0.75] from 16.688 19.072 to 50.778 58.032 \put {$\iota $} at 23.949 47.112 \endpicture \hfill \null \endgroup

\noindent commutes, where $\iota $ is the map sending each $\alpha $ in $\hull (\SX ) $ to the element in $C^*\big (\Gess \big )$ represented by the
characteristic function on the bisection $\Omega _\alpha $ formed by all germs $[\alpha ,x]$, as $x$ range in $ \Ehat \sub {ess}(\SX ) $.  In symbols
$$ \Omega _\alpha = \big \{[\alpha ,x]: x \in \Ehat \sub {ess}(\SX ) \big \}.  $$ Notice that $\Phi $ is onto $\MX $ because its range is a closed
$*$-subalgebra containing the range of $\rho $.  We next set out to prove that $\Phi $ is one-to-one.  Since $\hull (\SX ) $ is $0$-$E$-unitary by \ref
{HullEUnitary}, we have that $\Gess $ is Hausdorff by \cite [10.9]{actions}, so there exists a canonical conditional expectation $$ P:C^*\big (\Gess
\big )\to C_0( \Ehat \sub {ess}(\SX ) ) $$ given by restricting functions to the unit space, which we claim satisfies $$ P\big (\iota (\alpha )\big )
= \clauses { \kl \iota (\alpha ) / if / \alpha \in \ehull (\SX ) , / \kl 0 / otherwise, / / } \equationmark ClaimConExpOne $$ for every $\alpha $ in
$\hull (\SX ) $.  In order to prove the claim, recall that when $\alpha $ is idempotent, every germ $[\alpha ,x]$ is a unit, so the support of $\iota
(\alpha )$, also known as $\Omega _\alpha $, is contained in the unit space of $\Gess $, and hence indeed $P\big (\iota (\alpha )\big ) = \iota (\alpha )$.

When $\alpha $ is not idempotent we will prove that $\Omega _\alpha $ has an empty intersection with the unit space.  Arguing by contradiction assume
that there is a unit of the form $[\alpha ,x]$, for some $x$ in $\Ehat (\SX ) $.  This means that $[\alpha ,x]=[\varepsilon ,x]$, for some idempotent
$\varepsilon $ \cite [4.11]{actions}, which in turn implies that $\alpha \zeta =\varepsilon \zeta \neq 0$, for some idempotent $\zeta $, according
to \cite [4.6]{actions}.  Consequently $\varepsilon \zeta $ is a nonzero idempotent dominated by $\alpha $, and hence $\alpha $ itself is idempotent
because $\hull (\SX ) $ is $0$-$E$-unitary, as seen in \ref {HullEUnitary}.  Reaching a contradiction, we have thus proved that $\Omega _\alpha $ is
indeed disjoint from the unit space.  So, restricting the characteristic function on $\Omega _\alpha $ to the unit space leads to the zero function,
meaning that $P\big (\iota (\alpha )\big )=0$, as desired.

Since $\Gess $ is amenable by \ref {amenablesubshift}, $P$ is faithful.

We shall also employ another conditional expectation, this time defined on $\MX $, which we will now describe.  We initially consider the canonical
conditional expectation $Q$ defined on $B(\H )$ (according to the convention adopted in section \ref {MatsuSection}, we let $\H =\ell ^2(\LX )$
here) onto the set of all diagonal operators.  Given any $\alpha $ in $\hull (\SX ) $, and recalling that $\rho =q\circ \pi $, with $\pi $ as in \ref
{ExistMatsuRep}, we claim that $$ Q\big (\pi (\alpha )\big ) = \clauses { \kl \pi (\alpha ) / if / \alpha \in \ehull (\SX ) , / \kl 0 / otherwise. / /
} \equationmark ClaimConExpTwo $$ In order to verify this we exclude the trivial case in which $\alpha =0$, and we use \first {\FPFormOfHull } to write
$\alpha = \theta _\mu f_\Lambda \theta _\nu \inv $, where $\Lambda \subseteq \LX \cup \{1\}$ is finite, $\Lambda \cap \SX $ is nonempty, and $\mu ,
\nu \in \Lambda $.  Therefore $$ \pi (\alpha ) = T_\mu \pi (f_\Lambda )T_\nu ^*.  $$

If $\alpha $ is in $\ehull (\SX ) $, we may assume that $\mu =\nu $, whence $\pi (\alpha )$ is a diagonal idempotent operator, so $Q\big (\pi (\alpha
)\big ) =\pi (\alpha )$.  On the other hand, if $\alpha $ is not idempotent, then all diagonal entries of $\pi (\alpha )$ must vanish, because otherwise
there exists some basis vector $\delta _\xi $ such that $ \langle \pi (\alpha )\delta _\xi , \delta _\xi \rangle $ is nonzero.  This implies that $\pi
(\alpha )\delta _\xi $ is nonzero, so $\xi =\nu \xi '$, for some $\xi '$, and $$ \pi (\alpha )\delta _\xi = \big (T_\mu \pi (f_\Lambda )T_\nu ^*\big
)\delta _{\nu \xi '} = \delta _{\mu \xi '}, $$ whence $$ 0\neq \langle \pi (\alpha )\delta _\xi , \delta _\xi \rangle = \langle \delta _{\mu \xi '},
\delta _{\nu \xi '}\rangle .  $$ Therefore $\nu \xi '=\mu \xi '$, and we deduce that $\nu =\mu $, so $\alpha $ lies in $\ehull (\SX ) $, a contradiction,
hence proving \ref {ClaimConExpTwo}.

Observing that $Q$ leaves $K(\H )$ invariant, it factors through the quotient providing a linear operator $\tilde Q$ on $B(\H )/K(\H )$ such that $$
\tilde Q\big (U+K(\H ) \big ) = Q(U)+K(\H ) , $$ for every $U$ in $B(\H )$.  We then have that the diagram

\begingroup \noindent \hfill \beginpicture \setcoordinatesystem units <0.03truecm, -0.020truecm> \setplotarea x from -60 to 200, y from -30 to 130 \put
{\null } at -60 -30 \put {\null } at 200 130 \put {$C^*\big (\Gess \big )$} at 0 0 \put {$B(\H )/K(\H )$} at 120 0 \put {$C_0( \Ehat \sub {ess}(\SX ))$}
at 0 100 \put {$B(\H )/K(\H )$} at 120 100 \arrow <0.15cm> [0.25, 0.75] from 28.6 0 to 81 0 \put {$\Phi $} at 54.8 -13 \arrow <0.15cm> [0.25, 0.75]
from 39 100 to 81 100 \put {$\Phi $} at 60 87 \arrow <0.15cm> [0.25, 0.75] from 0 16.9 to 0 83.1 \put {$P$} at -13 50 \arrow <0.15cm> [0.25, 0.75]
from 120 16.9 to 120 83.1 \put {$\tilde Q$} at 133 50 \endpicture \hfill \null \endgroup

\noindent commutes, as this can be seen by checking on the elements of the form $\iota (\alpha )$, and using \ref {ClaimConExpOne} and \ref {ClaimConExpTwo}.

The proof of injectivity of $\Phi $ may now be given as follows: if $\Phi (a)=0$, for some $a$ in $C^*\big (\Gess \big )$, then $$ 0 = \tilde Q\big (\Phi
(a^*a)\big ) = \Phi \big (P(a^*a)\big ), $$ whence $P(a^*a)=0$, because $\Phi $ is injective on $C_0( \Ehat \sub {ess}(\SX ) )$ by \ref {SupportEss}.
Since $\Gess $ is amenable by \ref {amenablesubshift}, we have that $P$ is faithful, so $a=0$, as desired.  This shows that $\Phi $ is injective, proving (i).

The proof of (ii) follows essentially the same lines, the punch line coming after showing that the diagram

\begingroup \noindent \hfill \beginpicture \setcoordinatesystem units <0.035truecm, -0.020truecm> \setplotarea x from -60 to 200, y from -30 to 130
\put {\null } at -60 -30 \put {\null } at 200 130 \put {$C^*\big (\Gmax \big )$} at 0 0 \put {$B\big (\ell ^2(\X )\otimes \ell ^2(\F )\big )$} at
120 0 \put {$C_0( \Ehat \sub {max}(\SX ) )$} at 0 100 \put {$B\big (\ell ^2(\X )\otimes \ell ^2(\F )\big )$} at 120 100 \arrow <0.15cm> [0.25, 0.75]
from 27 0 to 75 0 \put {$\Phi $} at 51 -13 \arrow <0.15cm> [0.25, 0.75] from 36 100 to 75 100 \put {$\Phi $} at 55.5 87 \arrow <0.15cm> [0.25, 0.75]
from 0 19.5 to 0 80.5 \put {$P$} at -13 50 \arrow <0.15cm> [0.25, 0.75] from 120 19.5 to 120 80.5 \put {$Q$} at 133 50 \endpicture \hfill \null \endgroup

\noindent commutes, where $\Phi $ is now obtained by applying \cite [10.14]{actions} to the representation $\rho $ of \ref {RepMax}, while $P$ and $Q$
are similarly defined.  The key point here is that if $\alpha \in \hull (\SX )$ is not idempotent, then $d(\alpha )\neq 1$ and so $\rho (\alpha )=\pi
(\alpha )\otimes \lambda _{d(\alpha )}$ has no non-zero diagonal coefficients.  We leave the details for the reader.  \endProof

Let us conclude by noticing the following variation of \cite [Theorem 3.4]{MatsuCarl}:

\state Corollary If the subshift $\X $ is such that $\SX $ satisfies Carlsen and Matsumoto's condition $(*)$  (see \ref {CarlMatsuCond}), then $\MX $
is naturally isomorphic to $\OX $.

\Proof Follows immediately from \ref {CMConditionForSX} and the result above.  \endProof

Another consequence of Theorem \ref {Main}, in light of Corollary \ref {partialgroup} and \cite {Abadie}, is that both the Matsumoto and the
Carlsen-Matsumoto algebras are partial crossed products with respect to free group actions.  This point of view for the Carlsen-Matsumoto algebra was
taken in \cite {DokuchaExel}, but it may be new for the Matsumoto algebra.

\state Corollary Given a subshift $\X $, the Matsumoto algebra $\MX $ and the Carlsen-Matsumoto algebra $\OX $ can be realized as partial crossed
products of the free group on the alphabet of the shift with a commutative $C^*$-algebra.

\references

\def \ignore #1\par {}

\Article Abadie F. Abadie; On partial actions and groupoids; Proc. Amer. Math. Soc., 132 (2004), 1037-1047

\ignore \Bibitem AielloContiRossiStammeier V. Aiello, R. Conti, S. Rossi and N. Stammeier; The inner structure of boundary quotients of right LCM
semigroups; arXiv:1709.08839, 2017

\Bibitem BedosSpielberg E. B\'edos, S. Kaliszewski, J. Quigg and J. Spielberg; On finitely aligned left cancellative small categories, Zappa-Sz\'ep
products and Exel-Pardo algebras; arXiv:1712.09432, 2017

\Article BCEIsotopy M. Boyle, T. M. Carlsen and  S. Eilers; Flow equivalence and isotopy for subshifts; Dyn. Syst.,  32  (2017),  no. 3, 305-325

\Article BCESofic M. Boyle, T. M. Carlsen and  S. Eilers; Flow equivalence of sofic shifts; Israel J. Math.,  225  (2018),  no. 1, 111-146

\Bibitem BrixCarl K. A. Brix and T. M. Carlsen; Cuntz-Krieger algebras and one-sided conjugacy of shifts of finite type and their groupoids; arXiv:1712.00179
[math.OA]

\ignore \Bibitem BrownloweLarsenStammeier N. Brownlowe, N. S. Larsen and N. Stammeier; C*-Algebras of algebraic dynamical systems and right LCM semigroups;
to appear in Indiana Univ. Math. J., arXiv:1503.01599, 2017

\Article BFFG S. Bulman-Fleming, J. Fountain, and V. Gould; Inverse semigroups with zero: covers and their structure; J. Austral. Math. Soc. Ser. A,
67 (1999), 15-30

\Article CarlsenSofic T. M. Carlsen; On C*-algebras associated with sofic shifts; J. Operator Theory,  49  (2003),  no. 1, 203-212

\Bibitem CarlsenIntro T. M. Carlsen; An introduction to the C*-algebra of a one-sided shift space; Operator algebra and dynamics, 63--88, Springer
Proc. Math. Stat., 58, Springer, Heidelberg,  2013

\Bibitem CarlsenRigidity T. M. Carlsen; C*-rigidity of dynamical systems and \'etale groupoids; arXiv:1803.05326 [math.OA]

\Article CEOR T. M. Carlsen,  S. Eilers, E. Ortega and G. Restorff; Flow equivalence and orbit equivalence for shifts of finite type and isomorphism
of their groupoids; J. Math. Anal. Appl.,  469  (2019),  no. 2, 1088-1110

\Article MatsuCarl T. M. Carlsen and K. Matsumoto; Some remarks on the C*-algebras associated with subshifts; Math. Scand., 95 (2004), 145-160

\Bibitem CRST T. M. Carlsen, E. Ruiz, A. Sims and M. Tomforde; Reconstruction of groupoids and C*-rigidity of dynamical systems; arXiv:1711.01052 [math.OA]

\Article CarlsenSilvestrov T. M. Carlsen and S. Silvestrov; C*-crossed products and shift spaces; Expo. Math., 25 (2007), 275-307

\Article CarlsenSilvestrovKTheory T. M. Carlsen and S. Silvestrov; On the K-theory of the C*-algebra associated with a one-sided shift space;
Proc. Est. Acad. Sci.,  59  (2010),  no. 4, 272-279

\Article CarlsenThomsen T. M. Carlsen and K. Thomsen; The structure of the C*-algebra of a locally injective surjection; Ergodic Theory Dynam. Systems,
32  (2012),  no. 4, 1226-1248

\ignore \Article Cherubini A. Cherubini and M. Petrich; The Inverse Hull of Right Cancellative Semigroups; J. Algebra, 111 (1987), 74-113

\ignore \Bibitem CP A.~H. Clifford and G.~B. Preston; The algebraic theory of semigroups. Vol. I; Mathematical Surveys, No. 7.; American Mathematical
Society, Providence, R.I., 1961

\ignore \Article CoOne L. A. Coburn; The C*-algebra generated by an isometry I; Bull. Amer. Math. Soc., 73 (1967), 722-726

\ignore \Article CoTwo L. A. Coburn; The C*-algebra generated by an isometry II; Trans. Amer. Math. Soc., 137 (1969), 211-217

\Bibitem CELY J. Cuntz, S. Echterhoff,  X. Li and G. Yu; K-theory for group C*-algebras and semigroup C*-algebras; Oberwolfach Seminars, vol. 47, 2017

\Article CK J. Cuntz and W. Krieger; A class of C*-algebras and topological Markov chains; Invent. Math., 56 (1980), no. 3, 251-268

\Article DokuchaExel M. Dokuchaev and R. Exel; Partial actions and subshifts; J. Funct. Analysis, 272 (2017), 5038-5106

\Article ortho R. Exel; Partial representations and amenable Fell bundles over free groups; Pacific J. Math, 192 (2000), 39-63

\Article actions R. Exel; Inverse semigroups and combinatorial C*-algebras; Bull. Braz. Math. Soc. (N.S.), 39 (2008), 191-313

\ignore \Article semigpds R. Exel; Semigroupoid C*-Algebras; J. Math. Anal. Appl., 377 (2011), 303-318

\ignore \Article infinoa R. Exel and M. Laca; Cuntz-Krieger algebras for infinite matrices; J. reine angew. Math., 512 (1999), 119-172

\ignore \Article ExelPardo R. Exel and E. Pardo; The tight groupoid of an inverse semigroup; Semigroup Forum, 92 (2016), 274-303

\Article announce R. Exel and B. Steinberg; The inverse hull of $0$-left-cancellative semigroups; Proc. Int. Cong. of Math.,  1 (2018), 21-30.
arXiv:1710.04722 [math.OA]

\Bibitem ESOne R. Exel and B. Steinberg; Representations of the inverse hull of a $0$-left-cancellative semigroup; arXiv:1802.06281 [math.OA]

\Article HeadlundMorse G. A. Hedlund and M. Morse; Unending chess, symbolic dynamics and a problem in semigroups; Duke Math. J., 11 (1944), 1-7

\ignore \Article Howson A. G. Howson; On the intersection of finitely generated free groups; J. London Math. Soc., 29 (1954), 428-434

\ignore \Article Keimel K. Keimel; Alg\`ebres commutatives engendr\'ees par leurs \'el\'ements idempotents; Canad. J. Math.,  22 (1970), 1071-1078

\ignore \Article KumjianPask A. Kumjian and  D. Pask; Higher rank graph C*-algebras; New York J. Math., 6 (2000), 1-20

\ignore \Bibitem KwasniewskiLarsen B. K. Kwasniewski and N. S. Larsen; Nica-Toeplitz algebras associated with right tensor C*-precategories over right
LCM semigroups: part II examples; arXiv:1706.04951, 2017

\Article Lawson M. V. Lawson; The structure of $0$-$E$-unitary inverse semigroups. I. The monoid case; Proc. Edinburgh Math. Soc. (2), 42 (1999), 497-520,

\ignore \Article Li X. Li; Semigroup C*-algebras and amenability of semigroups; J. Funct. Anal., 262 (2012), 4302-4340

\Article MatsuOri K. Matsumoto; On C*-algebras associated with subshifts; Internat. J. Math., 8 (1997), 357-374

\Article MatsumotoDimension K. Matsumoto; Dimension groups for subshifts and simplicity of the associated C*-algebras; J. Math. Soc. Japan, 51 (1999), 679-698

\Article MatsuAuto K. Matsumoto; On automorphisms of C*-algebras associated with subshifts; J. Operator Theory,  44 (2000), 91-112

\Article MatsumotoStabilized K. Matsumoto; Stabilized C*-algebras constructed from symbolic dynamical systems; Ergodic Theory Dynam. Systems, 20
(2000), 821-841

\Article MatsOrbitEquivalenceTop K. Matsumoto; Orbit equivalence of topological Markov shifts and Cuntz-Krieger algebras; Pacific J. Math., 246  (2010),
no. 1, 199-225

\Article MatsOrbitEquivalenceOneSid K. Matsumoto; Orbit equivalence of one-sided subshifts and the associated C∗-algebras; Yokohama Math. J., 56
(2010),  no. 1-2, 59-85

\Article MatsAClass K. Matsumoto; A class of simple C∗-algebras arising from certain non-sofic subshifts; Ergodic Theory Dynam. Systems, 31  (2011),
no. 2, 459-482

\Article MatsSomeRemarks K. Matsumoto; Some remarks on orbit equivalence of topological Markov shifts and Cuntz-Krieger algebras; Yokohama Math. J.,
58  (2012), 41-52

\Article MatsClassificationOf K. Matsumoto; Classification of Cuntz-Krieger algebras by orbit equivalence of topological Markov shifts; Proc. Amer. Math. Soc.,
141  (2013),  no. 7, 2329-2342

\Article MatsStronglyContinuous K. Matsumoto; Strongly continuous orbit equivalence of one-sided topological Markov shifts; J. Operator Theory, 74
(2015),  no. 2, 457-483

\Article MatsOnFlow K. Matsumoto; On flow equivalence of one-sided topological Markov shifts; Proc. Amer. Math. Soc., 144  (2016),  no. 7, 2923-2937

\Article MatsUniformlyContinuous K. Matsumoto; Uniformly continuous orbit equivalence of Markov shifts and gauge actions on Cuntz-Krieger algebras;
Proc. Amer. Math. Soc., 145  (2017),  no. 3, 1131-1140

\Article MatsContinuousOrbit K. Matsumoto; Continuous orbit equivalence, flow equivalence of Markov shifts and circle actions on Cuntz-Krieger algebras;
Math. Z., 285  (2017),  no. 1-2, 121-141

\Article MatsTopologicalConjugacy K. Matsumoto; Topological conjugacy of topological Markov shifts and Cuntz-Krieger algebras; Doc. Math., 22  (2017), 873-915

\Article MatsImprimitivityBimodules K. Matsumoto; Imprimitivity bimodules of Cuntz-Krieger algebras and strong shift equivalences of matrices; Dyn. Syst.,
33  (2018),  no. 2, 253-274

\Article MatsRelativeMorita K. Matsumoto; Relative Morita equivalence of Cuntz-Krieger algebras and flow equivalence of topological Markov shifts;
Trans. Amer. Math. Soc., 370  (2018),  no. 10, 7011-7050

\Article MatsAShort K. Matsumoto; A short note on Cuntz splice from a viewpoint of continuous orbit equivalence of topological Markov shifts; Math. Scand.,
123  (2018),  no. 1, 91-100

\Article MatsStateSplitting K. Matsumoto; State splitting, strong shift equivalence and stable isomorphism of Cuntz-Krieger algebras; Dyn. Syst., 34
(2019),  no. 1, 93-112

\Article MatsMatuiContinuousOrbitMark K. Matsumoto and H. Matui; Continuous orbit equivalence of topological Markov shifts and Cuntz-Krieger algebras;
Kyoto J. Math., 54  (2014),  no. 4, 863-877

\Article MatsMatuiContinuousOrbitZeta K. Matsumoto and H. Matui; Continuous orbit equivalence of topological Markov shifts and dynamical zeta functions;
Ergodic Theory Dynam. Systems, 36  (2016),  no. 5, 1557-1581

\Article MilanSteinberg D. Milan and B. Steinberg; On inverse semigroup C*-algebras and crossed products; Groups Geom. Dyn., 8 (2014), 485-512

\ignore \Article Munn W. D. Munn; Brandt congruences on inverse semigroups; Proc. London Math. Soc., {\rm (3)} 14 (1964), 154-164.

\ignore \Article MurOne G. J. Murphy; Ordered groups and Toeplitz algebras; J. Operator Theory,   18 (1987), 303-326

\ignore \Article MurTwo G. J. Murphy; Ordered groups and crossed products of C*-algebras; Pacific J. Math., 2 (1991), 319-349

\ignore \Article MurThree G. J. Murphy; Crossed products of C*-algebras by semigroups of automorphisms; Proc.  London Math. Soc., 3 (1994), 423-448

\ignore \Article Nica A. Nica; C*-algebras generated by isometries and Wiener-Hopf operators; J. Operator Theory, 27 (1992), 17-52

\Bibitem Paterson A. L. T. Paterson; Groupoids, inverse semigroups, and their operator algebras; Birkh\umlaut auser, 1999

\ignore \Article Petrichcong M. Petrich; Congruences on inverse semigroups; J. Algebra, 55{\rm (2)} (1978), 231-256

\ignore \Bibitem petrich M. Petrich; Inverse semigroups; Pure and Applied Mathematics (New York). John Wiley \& Sons, Inc., New York, 1984.
A Wiley-Interscience Publication

\Bibitem cuntzlike J. Renault; Cuntz-like algebras; Proceedings of the 17th International Conference on Operator Theory (Timisoara 98), The Theta
Fondation, 2000

\Bibitem SpielbergA J. Spielberg; Groupoids and C*-algebras for categories of paths; arXiv:1111.6924v4, 2014

\Bibitem SpielbergB J. Spielberg; Groupoids and $C^*$-algebras for left cancellative small categories; arXiv:1712.07720, 2017

\ignore \Bibitem Stammeier N. Stammeier; A boundary quotient diagram for right LCM semigroups; Semigroup Forum, pp. 1-16 (2017, first online), arXiv:1604.03172

\ignore \Bibitem Starling C. Starling; Boundary quotients of C*-algebras of right LCM semigroups; arXiv:1409.1549, 2017

\ignore \Article SteinbergPrimitive B. Steinberg; Simplicity, primitivity and semiprimitivity of \'etale groupoid algebras with applications to inverse
semigroup algebras; J. Pure Appl. Algebra, 220 (2016), 1035-1054

\Article Thomsen K. Thomsen; Semi-\'etale groupoids and applications; Ann. Inst. Fourier (Grenoble), 60 (2010), 759-800

\endgroup

\close

\close